\documentclass[lettersize,onecolumn]{IEEEtran}

\usepackage{amsmath,amsfonts}
\usepackage{algorithmic}
\usepackage{algorithm}
\usepackage{array}
\usepackage[caption=false,font=normalsize,labelfont=sf,textfont=sf]{subfig}
\usepackage{textcomp}
\usepackage{stfloats}
\usepackage{url}
\usepackage{verbatim}
\usepackage{graphicx}
\usepackage{cite}

\usepackage{amssymb,amsmath}
\usepackage{graphics}
\usepackage{float}
\usepackage{color}
\usepackage{dsfont}
\usepackage{bm}
\usepackage{upgreek}
\usepackage{arydshln}
\usepackage{enumerate}
\usepackage{amsthm}
\usepackage{graphicx}
\usepackage{threeparttable}
\usepackage{caption}
\usepackage{multirow}
\usepackage{color}
\usepackage{xfrac}
\usepackage{array}
\usepackage{nomencl}
\usepackage{hyphenat}
\usepackage{algorithm,algorithmic}
\usepackage{url}
\usepackage[font=small]{caption}
\usepackage{hyphenat}
\usepackage{setspace}
\usepackage{cancel}
\usepackage[dvipsnames]{xcolor}

\def\QED{\mbox{\rule[0pt]{1.4ex}{1.4ex}}} 

\usepackage{soul}  

\usepackage{chngcntr}
\usepackage{apptools}
\usepackage{enumitem}

\newcommand{\innerprod}[2]{ \left\langle #1 , #2\right\rangle}
\newcommand{\parder}[2]{ \frac{\partial #1}{\partial #2} }
\newcommand{\norm}[1]{ || #1 ||}
\newcommand{\trace}[1]{ \mathrm{tr}( #1 )}
\newcommand{\mb}[1]{\mathbf{#1}}
\newcommand{\bs}[1]{\boldsymbol{#1}}

\newcommand{\virg}[1]{\textquotedblleft#1\textquotedblright}

\newcommand{\kronsigmatinvm}{\bs{\Sigma}_0^{-1/2}\otimes\bs{\Sigma}_0^{-1/2}}

\newcommand\e[1]{E\{ #1 \} }
\newcommand\eg[1]{E\left\lbrace  #1 \right\rbrace  }

\newcommand\ev[1]{{E_0\{ #1 \} }}

\newcommand{\cvec}[1]{ \mathrm{vec}(  #1 ) }
\newcommand{\vecs}[1]{\mathrm{vecs}(#1)}
\newcommand{\ovec}[1]{\underline{\mathrm{vec}}(#1)}
\newcommand{\ovecs}[1]{\underline{\mathrm{vecs}}(#1)}
\newcommand{\tonde}[1]{\left( #1 \right)  }
\newcommand{\quadre}[1]{\left[  #1 \right]  }
\newcommand{\graffe}[1]{\left\lbrace   #1 \right\rbrace   }

\newtheorem{theorem}{Theorem}[section]
\newtheorem{corollary}{Corollary}[section]
\newtheorem{lemma}{Lemma}[section]

\newtheorem{proposition}{Proposition}[section]

\newcommand{\vsVso}{\underline{\mathrm{vecs}}(\mb{V}_{S,0})}
\newcommand{\vsSs}{\mathrm{vecs}(\bs{\Sigma}_0)}

\newcommand{\vVsoinv}{\mathrm{vec}(\mb{V}_{S,0}^{-1})}

\newcommand{\kronVsoinvm}{\mb{V}_{S,0}^{-1/2}\otimes\mb{V}_{S,0}^{-1/2}}

\begin{document}

\title{Nuisance parameters and elliptically symmetric distributions: a geometric approach to parametric and semiparametric efficiency}

\author{Stefano Fortunati, Jean-Pierre Delmas~\IEEEmembership{Senior Member,~IEEE,} and Esa Ollila~\IEEEmembership{Senior Member,~IEEE.}
\thanks{This paper was produced by the IEEE Publication Technology Group. They are in Piscataway, NJ.}
\thanks{Manuscript received XXX, 2025; revised YYY, 2025.}}

\markboth{Journal of \LaTeX\ Class Files,~Vol.~1, No.~2, ZZZ~2025}%
{Shell \MakeLowercase{\textit{et al.}}: A Sample Article Using IEEEtran.cls for IEEE Journals}

\IEEEpubid{0000--0000~\copyright~2023 IEEE}

\maketitle

\begin{abstract}
	Elliptically symmetric distributions are a classic example of a semiparametric model where the location vector and the scatter matrix (or a parameterization of them) are the two finite-dimensional parameters of interest, while the density generator represents an \textit{infinite-dimensional nuisance} term. This basic representation of the elliptic model can be made more accurate, rich, and flexible by considering additional \textit{finite-dimensional nuisance} parameters. Our aim is therefore to investigate the deep and counter-intuitive links between statistical efficiency in estimating the parameters of interest in the presence of both finite and infinite-dimensional nuisance parameters. Previous seminal works have addressed this problem by leveraging a general result: if the statistical model has a specific group invariance, then the projection operator onto the semiparametric nuisance tangent space can be asymptotically expressed as a conditional expectation with respect to the maximal invariant sub-$\sigma$ algebra. In this article, we show that, for the statistical model of elliptical distributions, the projection operator can be explicitly computed without relying on the above-mentioned asymptotic approximation. This allows us to obtain original results also for the case in which the location vector and the scatter matrix are parameterized by a finite-dimensional vector that can be partitioned in two sub-vectors: one containing the parameters of interest and the other containing the nuisance parameters. As an example, we illustrate how the obtained results can be applied to the well-known \virg{low-rank} parameterization. Furthermore, while the theoretical analysis will be developed for Real Elliptically Symmetric (RES) distributions, we show how to extend our results to the case of Circular and Non-Circular Complex Elliptically Symmetric (C-CES and NC-CES) distributions.
\end{abstract}


\section{Introduction}
\IEEEPARstart{A}{}semiparametric model is a statistical model that involves not only a finite-dimensional parameter vector of interest $\bs{\theta} \in \Theta \subseteq \mathbb{R}^d$ but also an infinite-dimensional parameter, i.e. a function $g$, that often represents a nuisance parameter. This characterization is general enough to include many well-known examples: the symmetric location model, linear and logistic regression, errors in variables models, missing data and censoring models, copula models and even times series models such as ARMA or ARCH models. We refer the reader to \cite[Sect. 4]{BKRW} and \cite{Hallin_Werker} for a discussion of how the semiparametric formalism applies to the above-mentioned examples. The clear advantage of adopting a semiparametric model is in the potential gain in term of robustness with respect to some missing knowledge about the random experiment at hand that are indeed required when we use a parametric model. On the other hand, the fact that the function $g$ is left fully unspecified may lead to some efficiency losses in the estimation of the parameter vector of interest $\bs{\theta} \in \Theta \subseteq \mathbb{R}^d$ in respect of the parametric case. Nevertheless, there are cases in which parametric and semi-parametric efficiency coincide. In other words, in these specific cases, we can gain robustness without losing efficiency. This property of the statistical model is known as \textit{adaptivity} \cite{Bickel_paper,BKRW,Hallin_Werker}. 

In this work, we will analyze the parametric and the semiparametric efficiency for the estimation of the parameter vector of interest when the considered semiparametric model contains nuisance parameters of \textit{both finite and infinite dimension}. Specifically, we will focus on the statistical model of the elliptically symmetric distributions. 

Elliptically symmetric distributions have established themselves as a statistical model capable of capturing the heterogeneous nature of data in a wide range of applications: from remote sensing and communication to seismology and finance (see e.g. \cite{RES_rew, fang2004elliptically, gupta2013elliptically,Esa,Book_ES} for a complete list of references and examples). As for a Gaussian distribution, an elliptical distribution has the advantage of being fully characterized by the location vector $\bs{\mu}$ and the covariance/scatter matrix $\bs{\Sigma}$, while its flexibility with respect to (w.r.t.) the latter is provided by the \textit{density generator}, i.e. a function $g$, that is able to characterizes the \textit{heavy or light-tailed} behavior of the data. To this end, the density generator $g$ may depend on additional parameters that control the \virg{level of the tails} of the resulting distribution. Among the most popular and widely-used elliptical distributions, we may cite the Student $t$-distribution (characterized by the, so called, \textit{degree of freedom} $\nu$) and the Generalized Gaussian distribution (characterized by the \textit{scale} and \textit{shape} parameters) \cite{Chap_back_JP}.

Inference in elliptical distributions generally requires the estimation of $\bs{\mu}$ and $\bs{\Sigma}$ as main \textit{parameters of interest} in the eventual presence of additional \textit{nuisance} terms. For clarity, we recall that a nuisance term is an unknown parameter whose estimation is not strictly required but which can affect the performance of the estimator of the parameter of interest, i.e. the \textit{efficiency}, of the parameters of interest. 

To formalize this inference problem, we need to specify the statistical model we are considering. Three options, ranging from the most to the least restrictive, can be considered:
\begin{itemize}
	\item[M1] \textit{Parametric modeling with fully specified density generator $g$}. This is the least general case in which a full knowledge (i.e. both its functional form and its parameters) of the specific density generator $g$ is assumed to be a-priori known. To fix the idea, in order to derive estimation procedure for $\bs{\mu}$ and $\bs{\Sigma}$, a practitioner may assume \textit{to know a-priori} that the data follow a $t$-distribution with an \textit{a-priori known} degree of freedom $\nu$. This full knowledge of the density generator is not realistic in practice and need to be relaxed.
	\item[M2] \textit{Parametric modeling with specified density generator $g$ up to its parameters}. This is the classical approach in parametric elliptical inference and most of the literature deals with this case. Returning once again to the example of the $t$-distribution, in this case a practitioner may \textit{suppose a-priori} that the data generating process is a $t$-distribution characterized by the unknown degree of freedom $\nu$. Consequently, unlike in M1, the parameters of interest $\bs{\mu}$ and $\bs{\Sigma}$ need to be estimated together with the nuisance term $\nu$. In general this lack of knowledge on the true value of $\nu$ may lead to a performance degradation, i.e. \textit{efficiency losses}, in the estimation of $\bs{\mu}$ and $\bs{\Sigma}$. Even if more flexible w.r.t. the previous modeling approach, the requirement of the a-priori knowledge of the functional form of the density generator $g$ may be questionable. This leads to the semiparametric modeling.
	\item[M3] \textit{Semiparametric modeling where the functional form of the density generator $g$ is left unspecified}. This is the most realistic case in which the practitioner only supposes that the data are elliptically distributed, without assuming any specific density generator $g$, while estimating the parameters of interest $\bs{\mu}$ and $\bs{\Sigma}$. As a consequence, the density generator has to be considered as an infinite-dimensional nuisance term. Intuitively, one might expect efficiency losses in this case to be greater than those incurred in M2. After all, in model M3, it is the entire functional form of $g$ that is unknown, not just the value of its characterizing parameters.    
\end{itemize}

As this discussion suggests, it would be useful to carefully analyze the efficiency losses in the estimation of the parameters of interest when moving from the more to the less stringent modeling. In particular, the following question may arise: \textit{is it possible to relax unrealistic assumptions on the adopted statistical model (for example, moving from M2 to M3, or even from M1 to M3) without losing efficiency?} 

As showed by Hallin and Paindaveine in two seminal papers \cite{Hallin_P_2006,PAINDAVEINE}, the answer may be surprising and counter-intuitive. They proved that a decisive role is played by the additional finite-dimensional nuisance parameters that are involved, in an implicit or explicit manner, in the semiparametric modeling of elliptical distributions. Specifically, in their works, Hallin and Paidaveine built upon a general result, proposed in \cite[Prop. 2.6]{Hallin_Werker}, to bring to light the fundamental role of a finite-dimensional nuisance parameter hidden in the semiparametric elliptic model: the \textit{scale parameter}. In short, it is well-known that, in an elliptical model, the scatter matrix $\bs{\Sigma}$ and the density generator $g$ are not jointly identifiable due to a scale ambiguity. In order to remove this ambiguity, $\bs{\Sigma}$ must be rewritten as the product of a scale parameter $s \triangleq S(\bs{\Sigma})$ and of a \textit{shape matrix} $\mb{V}_S \triangleq \bs{\Sigma}/S(\bs{\Sigma})$, where $S(\cdot)$ is a given homogeneous functional of order one. In the resulting inference problem, $\bs{\mu}$ and $\mb{V}_S$ are to be considered as parameters of interest, while $s$ and $g$ are finite and infinite-dimensional nuisance parameters, respectively. Then, in \cite{Hallin_P_2006,PAINDAVEINE}, it has been shown that, if the scale parameter $s$ is considered as a nuisance term, then: \textit{i}) not knowing the functional form of $g$ does not lead to any efficiency loss on the estimation of $\bs{\mu}$ and $\mb{V}_S$ w.r.t. the case in which $g$ is a-priori known (so-called \textit{restricted adaptivity} \cite{Hallin_P_2006}), and \textit{ii}) if a determinant-based scale function $S(\cdot)$ is adopted, then not knowing $g$ and $s$ does not lead to efficiency losses on $\bs{\mu}$ and $\mb{V}_S$ w.r.t. the case in which both $g$ and $s$ are perfectly known (the \textit{full adaptivity}).

Inspired by this surprising result, in the first part of this paper, we focus on the three models M1, M2 and M3 discussed above. Specifically, our goal here is to recover the adaptivity conditions found in \cite{Hallin_P_2006,PAINDAVEINE}, but using a more direct approach. In fact, the approach proposed in \cite{Hallin_P_2006,PAINDAVEINE} is to derive this condition of adaptivity from a more general result proposed in \cite[Prop. 2.6]{Hallin_Werker}. Without going into detail, this general result shows that, for certain models, the projection operator can be asymptotically approximated as a conditional expectation with respect to the maximal invariant (with respect to the group that generates the statistical model) sub-$\sigma$ algebra. In particular, the generating group of the RES distribution, along with its maximal invariant statistics, has been introduced and discussed in \cite{Hallin_P_Annals}. In contrast, in this article we show how to derive, for the elliptical model, this projection operator can be derived directly without resorting to its asymptotic approximation.

To make this work as self-contained as possible, in Sect. \ref{sec_lem1_2} we present the geometrical tools (i.e. projections and tangent spaces) that will be at the core of our analysis of the semiparametric elliptical model. Moreover, two fundamental lemmas are provided and their proofs detailed in the Appendix. We move then to Sect. \ref{RES_models} where, after briefly recalling of the main definitions and properties of the Real Elliptically Symmetric (RES) distributions, we provide an extensive discussion of score vectors and related Fisher Information Matrix (FIM) for the parametric RES model. In doing this, we will make use of some fundamental outcomes obtained in \cite{Hallin_P_2006,PAINDAVEINE} for the inference of the location $\bs{\mu}$ and the shape $\mb{V}_S$ in the presence of the nuisance scale parameter $s$. In Sect. \ref{sec_geo_V_g} the geometrical tools introduced in Sect. \ref{sec_lem1_2} will be specified for the semiparametric RES elliptical models under considerations, while the main results of the first part of the paper are collected in Sect. \ref{sec:CRB}, in particular in Proposition \ref{Prop_chains}, where the closed-form expressions of parametric and semiparametric FIM and the related information bounds for the models in M1, M2 and M3. Proposition \ref{Prop_chains} plays a crucial role in formalizing and fully understanding the relationship between the (lack of) efficiency losses and the finite and infinite-dimensional nuisance terms involved in the considered elliptical parametric and semiparametric models.    

The second and last part of our paper, basically concentrated in Sect. \ref{sec_para}, deals with parametrized elliptical distributions. In particular, we suppose that the location vector $\bs{\mu}(\bs{\theta})$ and the covariance/scatter matrix $\bs{\Sigma}(\bs{\theta})$ are parameterized by a finite dimensional vector $\bs{\theta} = (\bs{\gamma}^T, \bs{\xi}^T)^T$, where the sub-vector $\bs{\gamma}$ contains the parameters of interest, while $\bs{\xi}$ collects the finite-dimensional nuisance terms. We aim then to investigate the efficiency losses on the following two scenarios:
\begin{itemize}
	\item[E1] \textit{Parametric modeling: estimation of $\bs{\gamma}$ in the presence of the finite-dimensional nuisance $\bs{\xi}$ with fully specified density generator $g$}. Again, this is the classical parametric context adopted in the vast majority of the applications. Nevertheless, as said before, the assumption of a perfect a-priori knowledge of the density generator $g$ may be unrealistic in practice.
	\item[E2] \textit{Semiparametric modeling: estimation of $\bs{\gamma}$ in the presence of the finite-dimensional nuisance $\bs{\xi}$ and the infinite-dimensional nuisance $g$}. This is the more realistic case in which we only need to assume that the data generating process follows an elliptical distribution while its density generator is considered as an infinite-dimensional unknown term. 
\end{itemize}
So we ask ourselves the same question as before: \textit{is it possible to relax the parametric assumption, moving from E1 to E2, without losing efficiency?} The answer to this question depends on the parameterization $\bs{\theta} \mapsto (\bs{\mu}(\bs{\theta}),\bs{\Sigma}(\bs{\theta}))$ at hand. Proposition \ref{prop_corr_1_RES}, which represents the main result of Sect. \ref{sec_para}, provides an \textit{adaptivity} condition that allows one to determine whether the given parameterization leads to efficiency losses or not. As an example, this condition will be applied to two cases of practical importance: \textit{i}) the parameterization $(\bs{\gamma}^T, \bs{\xi}^T)^T \mapsto (\bs{\mu}(\bs{\gamma}),\bs{\Sigma}(\bs{\xi}))$, i.e. when the location and the scatter matrix have no parameters in common and the \textit{ii}) \virg{low-rank} parameterization.  
Then, Sect. \ref{Applications} generalizes all the results previously obtained in the context of RES distributions to the case of Circular Complex Elliptically Symmetric (C-CES) and Non Circular CES (NC-CES) distributions. Finally Sect. \ref{Conclusion} concludes the paper and the technical proofs are reported in the Appendix.

\textit{Notation}: Italics indicates scalar quantities ($a$), lower case and upper case boldface indicate column vectors ($\mathbf{a}$) and matrices ($\mathbf{A}$). The superscripts $T$, $*$, $H$ and ${\#}$ indicate the transpose, the complex conjugate, the Hermitian and the Moore-Penrose inverse operators such that $\mb{A}^{- T} \triangleq (\mb{A}^{-1})^T = (\mb{A}^T)^{-1}$, $\mb{A}^{-*} \triangleq (\mb{A}^{-1})^* = (\mb{A}^*)^{-1}$, $\mb{A}^{-H} \triangleq (\mb{A}^{-1})^H = (\mb{A}^H)^{-1}$ and $\mb{A}^{\# T} \triangleq (\mb{A}^{\#})^T= (\mb{A}^T)^{\#}$. Moreover, $\mb{A}^{-1/2} \triangleq (\mb{A}^{-1})^{1/2}$ denotes any square root of the inverse of the symmetric positive definite matrix $\mb{A}$.
Each entry of a vector $\mb{a}$ and of a matrix $\mb{A}$ is indicated as $a_{i}\triangleq [\mb{a}]_{i}$ and $a_{ij}\triangleq [\mb{A}]_{ij}$, respectively. The symbol $\mathrm{vec}$ indicates the standard vectorization operator that maps column-wise the entry of an $m \times m$ matrix $\mb{A}$ in an $m^2$-dimensional column vector $\cvec{\mb{A}}$. The Hadamard product  $\mb{B}\odot\mb{C}$ is the matrix whose $(i,j)$-th element is $[\mb{B}]_{ij}[\mb{C}]_{ij}$. The Kronecker product ${\mb B}\otimes{\mb C}$ denotes the block matrix whose $(i,j)$ block element is $b_{ij}{\mb C}$, and the commutation matrix ${\mb K}_m$ is such that  $\cvec{\mb{A}^T}={\bf K}_m\cvec{\mb{A}}$.
The operator $\ovec{\mb{A}}$ defines the $m^2-1$-dimensional vector obtained from $\cvec{\mb{A}}$ by deleting its first element, i.e. 
$\cvec{\mb{A}} \triangleq [a_{11},\ovec{\mb{A}}^T]^T$. For any $m \times m$ symmetric matrix $\mb{A}$, $\vecs{\mb{A}}$ indicates the $m(m+1)/2$-dimensional vector of the entries of the lower triangular part of $\mb{A}$. The \textit{duplication matrix} $\mb{D}_m$ is implicitly defined as the unique $m^2 \times m(m+1)/2$ matrix satisfying $\mb{D}_m \vecs{\mb{A}}=\cvec{\mb{A}}$ for any symmetric matrix $\mb{A}$ \cite{Magnus1, Magnus2}. Let us now implicitly define the operator $\ovecs{\cdot}$ as $\vecs{\mb{A}} \triangleq [a_{11},\ovecs{\mb{A}}^T]^T$.
The identity matrix of dimension $m$ is indicated as $\mb{I}_m$ and $\mb{e}_{k,l}$ indicates the $k$th vector of the canonical basis of $\mathbb{R}^{l}$. Moreover, let $\underline{\mb{I}}_m$ be the operator such that $\ovecs{\mb{A}}=\underline{\mb{I}}_m \vecs{\mb{A}}$ that can be obtained from $\mb{I}_{m(m+1)/2}$ by removing its first row.

Let $\mb{A}(\bs{\theta})$ be a matrix (or possibly a vector or even a scalar) function of the \textit{real} vector $\bs{\theta} \in \Theta \subseteq \mathbb{R}^d$, then $\mb{A}_0 \triangleq \mb{A}(\bs{\theta}_0)$ while $\mb{A}_i^{0}\triangleq {\frac{\partial \mb{A}(\bs{\theta})}{\partial \theta_i}}|_{\bs{\theta}=\bs{\theta}_0}$, where $\bs{\theta}_0$ is a particular (or \textit{true}) value of $\bs{\theta} \in \Theta$. Similarly, the gradient of a function $f(\bs{\theta})$ evaluated at $\bs{\theta}_0$, i.e. $\nabla_{\bs{\theta}}f(\bs{\theta})|_{\bs{\theta} = \bs{\theta}_0}$, will be indicated as $\nabla_{\bs{\theta}}f(\bs{\theta}_0)$. 

Let $(\mathcal{X},\mathfrak{B}(\mathcal{X}),P_0)$ be a probability space where the sample space $\mathcal{X}$ is a subset of $\mathbb{R}^m$, $\mathfrak{B}(\mathcal{X})$ is the Borel $\sigma$-algebra on $\mathcal{X}$ and $P_0$ is a probability measure. Moreover, $P_0$ is assumed to be \textit{absolutely continuous} with probability density function (pdf), associated to the Lebesgue measure on $\mathbb{R}^m$, given by $dP_0(\mb{x})=p_0(\mb{x})d\mb{x}$. Let $f:\mathcal{X}\rightarrow \mathbb{R}$ be an $\mathfrak{B}(\mathcal{X})$-measurable function, then $E_0\{f\} \triangleq \int f(\mb{x})dP_0(\mb{x})$ indicates its expectation w.r.t. $P_0$. Let $(\Omega, \mathcal{A},P)$ be a probability space and let $(\Psi, \mathcal{E})$ be a measurable space. Then, two random variables (that is, two measurable functions), $X,\;Y: (\Omega, \mathcal{A})  \rightarrow (\Psi,\mathcal{E})$ are said to be \textit{almost sure equal}, denoted as $X = Y$, if and only if $X(\omega) = Y(\omega)$ for $P$-almost all $\omega \in \Omega$.

Let us now introduce the Hilbert space $(\mathcal{H},\innerprod{\cdot}{\cdot}_{\mathcal{H}})$ as the (infinite-dimensional) linear space of the $\mathfrak{B}(\mathcal{X})$-measurable scalar functions with zero-mean and finite variance:
\begin{equation}\label{H_set}
	\mathcal{H} \triangleq \graffe{h \in L_2(P_0)|\ev{h} = 0} = \graffe{h:\mathcal{X}\rightarrow \mathbb{R}|\ev{h} = 0, \ev{h^2}< +\infty},
\end{equation}
endowed with the canonical inner product 
\begin{equation}
	\label{H_inner_prod}
	\innerprod{h_1}{h_2}_{\mathcal{H}}\triangleq \ev{h_1h_2} = \int_\mathcal{X} h_1(\mb{x})h_2(\mb{x})dP_0(\mb{x}), \quad \forall h_1,h_2 \in \mathcal{H}.
\end{equation}
We note that the norm associated to the inner product in \eqref{H_inner_prod} is $\norm{h}_{\mathcal{H}}=\sqrt{\ev{h^2}}$ that is the standard deviation of $h \in \mathcal{H}$.

Let us now introduce the $q$-\textit{replicating} Hilbert space $\mathcal{H}_q = \mathcal{H} \times \cdots \times \mathcal{H}$ as the linear space of the $\mathfrak{B}(\mathcal{X})$-measurable, $q$-variate functions $\mb{h}: \mathcal{X} \rightarrow \mathbb{R}^q$. This set can clearly be obtained as $q$ Cartesian products of $\mathcal{H}$ in \eqref{H_set} and this explain the name of $q$-replicating space. 


\textit{Remark}: It is worth noticing that the $q$-replicating space $\mathcal{H}_q$ has been introduced only to simplify some notation and to improve the readability of some result proposed in the following sections. However, it does not introduce any additional geometrical structure since it is just composed of $q$ copies of $(\mathcal{H},\innerprod{\cdot}{\cdot}_{\mathcal{H}})$. In particular, when we work on (subspaces of) $\mathcal{H}$, instead of $\mathcal{H}_q$, all the operators, such as expectations and projections, has to be interpreted component-wise (see e.g. \cite[Sect. 2.4, Remark 2]{Tsiatis} or \cite[Sect. 2.4]{BKRW}) as is proved in Lemma \ref{lem_p} and Theorem \ref{Pyt_theo} reported in Appendix \ref{App1}. Specifically, for each $\mb{h} = (h_1,\ldots,h_q)^T \in \mathcal{H}_q$ s.t. $h_i \in \mathcal{H}$ and $u \in \mathcal{U} \subseteq \mathcal{H}$, we have:
\begin{equation}\label{pr_comp}
	[\Pi(\mb{h}|\mathcal{U})]_i \triangleq \Pi([\mb{h}]_i|\mathcal{U}) = \Pi(h_i|\mathcal{U}),
\end{equation}
\begin{equation}
	[\ev{\mb{h}u}]_i \triangleq \ev{[\mb{h}]_iu} = \innerprod{h_i}{u}_{\mathcal{H}},
\end{equation}
for $i=1,\ldots,q$. 

\section{Models, score vectors and tangent spaces} \label{sec_lem1_2}
A semiparametric model is a set of probability density functions (pdf) parameterized by a finite dimensional vector of interest and by an infinite-dimensional parameter, i.e. a function, that generally plays the role of a nuisance term. The study of the estimation efficiency in semiparametric models, along with the derivation of the relevant information bounds, is a well-established topic in the statistical literature (see e.g. \cite{BKRW, Tsiatis, Newey} and the reference therein). We refer the reader to \cite{For_EUSIPCO_18, Chap_sem} for a tutorial introduction on this subject. In this work, we focus on the case in which the considered semiparametric model involves an additional finite-dimensional nuisance vector, along with the infinite-dimensional nuisance term. More formally, let us consider the general semiparametric model:
\begin{equation}\label{sp1}
	\mathcal{P} = \graffe{p_X(\mb{x}|\bs{\gamma},\bs{\xi},g): \bs{\gamma}\in \Gamma, \bs{\xi}\in \Psi,  
		g \in \mathcal{G} },
\end{equation}
where $\Gamma \subseteq \mathbb{R}^q$ is the set of the (finite-dimensional) parameter vectors $\bs{\gamma}$ of interests, $\Psi \subseteq \mathbb{R}^r$ is the set of (finite-dimensional) nuisance parameter vectors $\bs{\xi}$ and $\mathcal{G}$ is the set of the (infinite-dimensional) nuisance functions $g$. 
We indicate as $\bs{\gamma}_0 \in \Gamma$, $\bs{\xi}_0 \in \Psi$ and $g_0 \in \mathcal{G}$ the true, but unknown, related quantities and with $E_0\{f\} \triangleq \int f(\mb{x})dP_X(\mb{x}|\bs{\gamma}_0,\bs{\xi}_0,g_0)$ the expectation of a $\mathfrak{B}(\mathcal{X})$-measurable function $f$ w.r.t. the true distribution $P_0(\mb{x})=P_X(\mb{x}|\bs{\gamma}_0,\bs{\xi}_0,g_0)$. 

By using a self-explanatory notation, we now introduce three \textit{parametric} sub-models of $\mathcal{P}$ as:
\begin{equation}\label{P1}
	\mathcal{P}_{1} = \graffe{p_X(\mb{x}|\bs{\gamma},\bs{\xi}_0,g_0): \bs{\gamma}\in \Gamma },
\end{equation}
\begin{equation}\label{P2}
	\mathcal{P}_{2} = \graffe{p_X(\mb{x}|\bs{\gamma}_0,\bs{\xi},g_0): \bs{\xi}\in \Psi},
\end{equation}
\begin{equation}\label{P12}
	\mathcal{P}_{1,2} = \graffe{p_X(\mb{x}|\bs{\gamma},\bs{\xi},g_0): \bs{\gamma}\in \Gamma,  \bs{\xi}\in \Psi},
\end{equation}
along with the \textit{non-parametric} model
\begin{equation}\label{P3}
	\mathcal{P}_{3} = \graffe{p_X(\mb{x}|\bs{\gamma}_0,\bs{\xi}_0,g):  g \in \mathcal{G}}.
\end{equation}

Let us define the score vector $\mb{s}_{\bs{\gamma}_0}$ in $\mathcal{P}_{1}$ as:
\begin{equation}
	\label{score_vect_nu}
	\left[ \mb{s}_{\bs{\gamma}_0} \right]_i \triangleq  \left[ \mb{s}_{\bs{\gamma}_0}(\mb{x}) \right]_i = \left. \partial\ln p_X(\mb{x}|\bs{\gamma}, \bs{\xi}_0, g_0)/\partial \gamma_i \right|_{\bs{\gamma}=\bs{\gamma}_0},
\end{equation}
for $i=1,\ldots,q$, that represents the score vector of the parameters of interest. Similarly, the score vector $\mb{s}_{\bs{\xi}_0}$ of the finite-dimensional nuisance parameters $\bs{\xi}_0$ in $\mathcal{P}_{2}$ is given by:
\begin{equation}
	\label{score_vect_zeta}
	\left[ \mb{s}_{\bs{\xi}_0} \right]_j  \triangleq  \left[ \mb{s}_{\bs{\xi}_0}(\mb{x}) \right]_j =  \left. \partial \ln p_X(\mb{x}|\bs{\gamma}_0, \bs{\xi}, g_0)/\partial \xi_j \right|_{\bs{\xi}=\bs{\xi}_0}, 
\end{equation}
for $j=1,\ldots,r$. Under the regularity conditions discussed in \cite[Sects. 6.2, 6.3]{Lehmann}, it is immediate to verify that \mbox{$[ \mb{s}_{\bs{\gamma}_0} ]_i,[ \mb{s}_{\bs{\xi}_0} ]_j \in \mathcal{H}$}, $\forall i,j$, i.e. they have zero-mean and finite variance.

We can now introduce the (finite-dimensional) tangent space of the parametric sub-models $\mathcal{P}_{2}$ as the linear span of $\mb{s}_{\bs{\xi}_0}$ in $\mathcal{H}$ \cite[Sect. 2.3]{Tsiatis}, \cite[eq. (5.24)]{Chap_sem}:
\begin{equation}
	\label{ts_2}
	\mathcal{H}  \supseteq  \mathcal{T}_{2} \triangleq \mathrm{Span}\{[ \mb{s}_{\bs{\xi}_0}]_1, \ldots, [ \mb{s}_{\bs{\xi}_0}]_r\}.
\end{equation}

Finally, the nuisance tangent space $\mathcal{T}_3 \subseteq \mathcal{H}$ of the non-parametric model $\mathcal{P}_3$ is defined as in \cite[Sect. 3.2, Def. 2]{BKRW}, \cite[Sect. 4.4]{Tsiatis}. Note that, by construction, $\mathcal{T}_{2}$ and $\mathcal{T}_{3}$ are finite- and infinite-dimensional closed subspaces of $\mathcal{H}$.

According to the previous definition, let us define the \textit{efficient score vector} $\bar{\mb{t}}_{\bs{\gamma}_0}$ for the vector $\bs{\gamma}_0$ of the parameter of interest in the parametric sub-model $\mathcal{P}_{1,2}$ in \eqref{P12} as \cite[Sect. 3.4]{Tsiatis}, \cite[Sect. 2.4]{BKRW} and \cite[Def. 4]{Chap_sem}:
\begin{equation}\label{t_def}
	\bar{\mb{t}}_{\bs{\gamma}_0} \triangleq \mb{s}_{\bs{\gamma}_0} -\Pi(\mb{s}_{\bs{\gamma}_0}|\mathcal{T}_{2}),
\end{equation}
where the projection $\Pi(\mb{s}_{\bs{\gamma}_0}|\mathcal{T}_{2})$ is to be interpreted \textit{component-wise} as indicated in \eqref{pr_comp}.
Since $\mathcal{T}_{2}$ is a finite-dimensional subspace of $\mathcal{H}$, the projection operator $\Pi(\cdot|\mathcal{T}_{2})$ can be derived in closed form as \cite[eq. (5.7)]{Chap_sem} \cite[Sect. 2.4, Ex. 1]{Tsiatis}:
\begin{equation}\label{proj_T2}
	\Pi(h|\mathcal{T}_2) = E_0\{h\mb{s}_{\bs{\xi}_0}^T\}\mb{I}_{\bs{\xi}_0}^{-1}\mb{s}_{\bs{\xi}_0}, \; h \in \mathcal{H},
\end{equation}
where:
\begin{equation}
	\mb{I}_{\bs{\xi}_0} =E_0\{\mb{s}_{\bs{\xi}_0}\mb{s}_{\bs{\xi}_0}^T\},
\end{equation}
is the Fisher Information Matrix (FIM) for $\bs{\xi}_0$ in the parametric sub-model $\mathcal{P}_{2}$ in \eqref{P2}. Then, we have that $\bar{\mb{t}}_{\bs{\gamma}_0}$ can be explicitly expressed as:
\begin{equation}\label{t_def2}
	\bar{\mb{t}}_{\bs{\gamma}_0} = \mb{s}_{\bs{\gamma}_0} - E_0\{\mb{s}_{\bs{\gamma}_0}\mb{s}_{\bs{\xi}_0}^T\}\mb{I}_{\bs{\xi}_0}^{-1}\mb{s}_{\bs{\xi}_0} = \mb{s}_{\bs{\gamma}_0} -  \mb{I}_{\bs{\gamma}_0\bs{\xi}_0}\mb{I}_{\bs{\xi}_0}^{-1}\mb{s}_{\bs{\xi}_0},
\end{equation} 
where 
\begin{equation}
	\mb{I}_{\bs{\gamma}_0\bs{\xi}_0} \triangleq E_0\{\mb{s}_{\bs{\gamma}_0}\mb{s}_{\bs{\xi}_0}^T\}
\end{equation}
is the matrix of the cross-information terms in the parametric sub-model $\mathcal{P}_{1,2}$ in \eqref{P12}.

Using the same geometrical approach, the \textit{semiparametric efficient score vector} $\bar{\mb{s}}_{\bs{\gamma}_0}$ for $\bs{\gamma}_0$ in the semiparametric model $\mathcal{P}$ in \eqref{sp1} is given by \cite[Sect. 3.4, eq. (18)]{BKRW}:
\begin{equation}\label{bar_s_nu}
	\bar{\mb{s}}_{\bs{\gamma}_0} \triangleq \mb{s}_{\bs{\gamma}_0} - \Pi(\mb{s}_{\bs{\gamma}_0}|\mathcal{T}_{2} + \mathcal{T}_{3}),
\end{equation}
where:
\begin{equation}\label{H}
	\mathcal{H} \supseteq \mathcal{T}_{2} + \mathcal{T}_{3} \triangleq \graffe{h \in \mathcal{H} | h = o + l, o \in \mathcal{T}_{2},  l \in \mathcal{T}_{3}},
\end{equation}
and since $\mathcal{T}_{2}$ is a (closed) finite-dimensional subspace and $\mathcal{T}_{3}$ is closed, then  $\mathcal{T}_{2} + \mathcal{T}_{3}$ is closed. In general, since $\mathcal{T}_{2} + \mathcal{T}_{3}$ is infinite-dimensional, a closed form for the projection operator $\Pi(\cdot|\mathcal{T}_{2} + \mathcal{T}_{3})$ does not exist. Fortunately, some further manipulation is still possible. In fact, let us first recall that, for two orthogonal closed subspaces $\mathcal{A}$ and $\mathcal{B}$ of $\mathcal{H}$ we have the following property:
\begin{equation}\label{orthogonal projector}
\Pi(h| \mathcal{A}+\mathcal{B})
=
\Pi(h| \mathcal{A})
+
\Pi(h| \mathcal{B}),\; \forall h \in \mathcal{H} .
\end{equation}
Moreover, it can be noted that $\mathcal{T}_{2} + \mathcal{T}_{3}$ can be expressed as the orthogonal direct sum of the two orthogonal subpsaces $\mathcal{T}_{2}$ and \hbox{$((\mathcal{T}_{2} + \mathcal{T}_{3})\cap \mathcal{T}_{2}^{\bot})$}. Then, from \eqref{orthogonal projector}, we immediately have that:
\begin{eqnarray}
\nonumber
\bar{\mb{s}}_{\bs{\gamma}_0}
& \triangleq&
\mb{s}_{\bs{\gamma}_0} - \Pi(\mb{s}_{\bs{\gamma}_0}|\mathcal{T}_{2} + \mathcal{T}_{3})
\\
\nonumber
& =& 
\mb{s}_{\bs{\gamma}_0} - \Pi(\mb{s}_{\bs{\gamma}_0}|\mathcal{T}_{2}) - \Pi(\mb{s}_{\bs{\gamma}_0}|(\mathcal{T}_{2} + \mathcal{T}_{3}) \cap \mathcal{T}_{2}^\perp) 
\\
\label{eff_score_nu}	
&  =& 
\bar{\mb{t}}_{\bs{\gamma}_0} - \Pi(\bar{\mb{t}}_{\bs{\gamma}_0}|(\mathcal{T}_{2} + \mathcal{T}_{3}) \cap \mathcal{T}_{2}^\perp),
\end{eqnarray}
where the last equality comes from \eqref{t_def} and from the fact that $\mathcal{T}_{2} \supset \Pi(\mb{s}_{\bs{\gamma}_0}|\mathcal{T}_{2}) \perp (\mathcal{T}_{2} + \mathcal{T}_{3}) \cap \mathcal{T}_{2}^\perp$.

The relation \eqref{eff_score_nu} between the efficient score vectors $\bar{\mb{s}}_{\bs{\gamma}_0}$ and $\bar{\mb{t}}_{\bs{\gamma}_0}$ is the key tool to compare the efficient Semiparametric FIM (SFIM) 
\begin{equation}\label{eff_sem_FIM}
	\bar{\mb{I}}(\bs{\gamma}_0|\bs{\xi}_0,g_0)\triangleq E_0\graffe{\bar{\mb{s}}_{\bs{\gamma}_0}\bar{\mb{s}}_{\bs{\gamma}_0}^T}
\end{equation}
for the parameter of interest $\bs{\gamma}_0$ in the presence of both the finite- and infinite-dimensional nuisance parameters $\bs{\xi}_0$ and $g_0$ with the efficient FIM
\begin{equation}\label{eff_FIM}
	\bar{\mb{I}}(\bs{\gamma}_0|\bs{\xi}_0)\triangleq E_0\graffe{\bar{\mb{t}}_{\bs{\gamma}_0}\bar{\mb{t}}_{\bs{\gamma}_0}^T}
\end{equation} 
for $\bs{\gamma}_0$ in the presence of only the finite-dimensional nuisance $\bs{\xi}_0$, while $g_0$ is known. Specifically, for this comparison, we will make use of the following lemma, stated without proof in \cite[Sect. 3.4, Prop. 3]{BKRW}, for which we provide a full proof in the Appendix \ref{App2} of this paper.
\begin{lemma}\label{Lem1}
The efficient SFIM $\bar{\mb{I}}(\bs{\gamma}_0|\bs{\xi}_0,g_0)$ and the efficient FIM $\bar{\mb{I}}(\bs{\gamma}_0|\bs{\xi}_0)$ for $\bs{\gamma}_0$ in the presence of respectively, the finite- and infinite-dimensional nuisance terms $\bs{\xi}_0$ and $g_0$, and only the finite-dimensional nuisance term $\bs{\xi}_0$, are connected through the relation:
\begin{equation}\label{Sem_FIM_eff}
		\bar{\mb{I}}(\bs{\gamma}_0|\bs{\xi}_0,g_0) = \bar{\mb{I}}(\bs{\gamma}_0|\bs{\xi}_0) -E_0\graffe{\mb{p}\mb{p}^T},
\end{equation}
where 
\begin{equation}\label{def_w}
		\mb{p} \triangleq \Pi\tonde{\bar{\mb{t}}_{\bs{\gamma}_0}|(\mathcal{T}_{2} + \mathcal{T}_{3}) \cap \mathcal{T}_{2}^\perp}.
\end{equation}
\end{lemma}

It is worth noticing that the matrix $\bar{\mb{I}}(\bs{\gamma}_0|\bs{\xi}_0)$ is the efficient FIM for $\bs{\gamma}_0$ in the presence of the finite-dimensional nuisance vector $\bs{\xi}_0$ in the parametric sub-model $\mathcal{P}_{1,2}$ in \eqref{P12}. Consequently the Cram\'{e}r-Rao bound (CRB) for $\bs{\gamma}_0$ in the presence of $\bs{\xi}_0$ is given by:
\begin{equation}\label{CRB_12}
	\mathrm{CRB}(\bs{\gamma}_0|\bs{\xi}_0) = \bar{\mb{I}}(\bs{\gamma}_0|\bs{\xi}_0)^{-1} = \quadre{\mb{I}_{\bs{\gamma}_0} - \mb{I}_{\bs{\gamma}_0\bs{\xi}_0}\mb{I}_{\bs{\xi}_0}^{-1}\mb{I}_{\bs{\gamma}_0\bs{\xi}_0}^T}^{-1},
\end{equation}  
where:
\begin{equation}
	\mb{I}_{\bs{\gamma}_0} \triangleq E_0\{\mb{s}_{\bs{\gamma}_0}\mb{s}_{\bs{\gamma}_0}^T\},
\end{equation}
is the FIM for $\bs{\gamma}_0$ in the parametric sub-model $\mathcal{P}_{1}$ in \eqref{P1}. It is immediate to recognize in \eqref{CRB_12} the well-known expression of the CRB for parametric estimation in the presence of a finite-dimensional nuisance vector.

Moreover, we note that this lemma implies that $ \bar{\mb{I}}(\bs{\gamma}_0|\bs{\xi}_0,g_0)\leq \bar{\mb{I}}(\bs{\gamma}_0|\bs{\xi}_0)$. This means that the information about the parameter $\bs{\gamma}_0$ is reduced or remains the same in the presence of the infinite-dimensional nuisance function $g_0$. It is clear from Lemma \ref{Lem1} that there is no loss of information iff $\mb{p}=\mb{0}$. The following lemma (also stated without proof, in \cite[Sect. 3.4, Prop. 3]{BKRW}) and proven in the Appendix \ref{App31} of this work, specifies the condition under which the presence of the nuisance function $g_0$ will not bring any loss of information.

\begin{lemma}\label{Lem2}
	The efficient SFIM $\bar{\mb{I}}(\bs{\gamma}_0|\bs{\xi}_0,g_0)$ for the model $\mathcal{P}$ in \eqref{sp1} is equal to the parametric efficient FIM $\bar{\mb{I}}(\bs{\gamma}_0|\bs{\xi}_0)$ in $\mathcal{P}_{1,2}$, i.e. $\mb{p}= \mb{0}$ if and only if (iff) the following condition is satisfied
	\begin{equation}\label{main_cond}
	\mb{s}_{\bs{\gamma}_0} - \Pi(\mb{s}_{\bs{\gamma}_0}|\mathcal{T}_2)\triangleq \bar{\mb{t}}_{\bs{\gamma}_0}   \perp \mathcal{T}_3.
	\end{equation}
\end{lemma}

In Sect. \ref{sec_geo_V_g} and in particular in Sect. \ref{sec_para}, we will make extensive use of this framework and of Lemma \ref{Lem2} to bring to light the sometimes surprising and counter-intuitive relationships between parametric and semiparametric efficiency in the family of elliptical distributions.

\section{A short introduction to the RES distributions}
\label{RES_models}
In this section, we briefly recall the definition and the main properties of the RES distributions distributions. Many different, yet equivalent, representation of the RES family can be found in the literature \cite{CAMBANIS1981,RES_Fang,Hallin_P_2006,PAINDAVEINE,Chap_back_JP}. Here, we adopt the approach and the notation introduced in \cite{Chap_back_JP}. Moreover, in the following, we will consider only the \textit{absolutely continuous case}, i.e. we suppose that each distribution admits a density w.r.t. the Lebesgue measure on $\mathbb{R}^m$.  
Before moving on, we would like to underline that some of the outcomes discussed in the following sections have already been presented in \cite{Hallin_P_2006,PAINDAVEINE} by using a \virg{semiparametric generalization} of the Le Cam theory on Local Asymptotically Normal (LAN) families of distributions \cite{LeCam60}, \cite[Ch.6]{LeCam_Yang}. Our main goal here is to recast the problem of the statistical inference in RES distributions in the framework of the Hilbert spaces and present a systematic analysis based on purely geometrical concepts, such as those used in Lemmas \ref{Lem1} and \ref{Lem2}.        
\subsection{Essentials on RES distributions}
\label{RES_models_essentials}
A real-valued, random observation vector $\mb{x} \in \mathcal{X} \subseteq \mathbb{R}^m$ is said to be elliptically symmetric distributed if its probability density function (pdf) can be expressed (in the absolutely continuous case) as: 
\begin{equation}
\label{RES_pdf}
p_X(\mb{x}|\bs{\mu},\bs{\Sigma},g)=|\bs{\Sigma}|^{-1/2} g \left((\mb{x}-\bs{\mu})^T\bs{\Sigma}^{-1}(\mb{x}-\bs{\mu}) \right),
\end{equation}
where $\bs{\mu} \in \mathbb{R}^m$ is a location vector, $\bs{\Sigma} \in \mathcal{S}_m^\mathbb{R}$ is an $m \times m$, positive definite, \textit{scatter} matrix in the set $\mathcal{S}_m^\mathbb{R}$ of the symmetric real matrices. \footnote{In this article we will limit ourselves to considering scatter, covariance and shape matrices as elements of the linear subspace of symmetric matrices and not as elements of the manifold of positive definite matrices.} The \textit{density generator} $g \in \mathcal{G}$ is a function belonging to a set $\mathcal{G}$ such that:
\begin{equation}
\label{set_G}
\mathcal{G} = \graffe{ g: \mathbb{R}^{+} \rightarrow \mathbb{R}_0^{+} \left| \delta_m \triangleq \int_{0}^{\infty}t^{m/2-1}g(t)dt = \pi^{-m/2}\Gamma(m/2) \right. }.
\end{equation}
where the value of $\delta_m$ is such that \eqref{RES_pdf} is a proper density that integrates to 1. In the following, the notation $\mb{x} \sim RES_m(\bs{\mu},\bs{\Sigma},g)$ indicates that a random vector $\mb{x} \in \mathcal{X}$ has the density given in \eqref{RES_pdf}.

A fundamental result for RES distributed vectors is the \textit{Stochastic Representation Theorem}. Specifically, a RES distributed vector $\mb{x} \sim RES_m(\bs{\mu},\bs{\Sigma},g)$ can be expressed as:
\begin{equation}
\label{SRT_dec}
\mb{x} = \bs{\mu} + \sqrt{\mathcal{Q}}\bs{\Sigma}^{1/2}\mb{u},
\end{equation}
where the random vector $\mb{u} \sim \mathcal{U}(S_{\mathbb{R}}^{m-1})$ is uniformly distributed on the unit sphere $S_{\mathbb{R}}^{m-1} \triangleq \{\mb{u}\in \mathbb{R}^m|\norm{\mb{u}}=1\}$ and consequently satisfies $E\{\mb{u}\}=\mb{0}$ and $E\{\mb{u}\mb{u}^T\}=m^{-1}\mb{I}_m$. The positive random variable $\mathcal{Q}$, called \textit{2nd-order modular variate}, is such that (s.t.)
\begin{equation}
\label{Q_RES}
\mathcal{Q} = Q_{\bs{\mu},\bs{\Sigma}}(\mb{x})\triangleq (\mb{x}-\bs{\mu})^T\bs{\Sigma}^{-1}(\mb{x}-\bs{\mu}),\; \mb{x} \in \mathcal{X}
\end{equation}
and it is independent of $\mb{u} \sim \mathcal{U}(S_{\mathbb{R}}^{m-1})$. \footnote{Formally, given an abstract probability space $(\Omega, \mathcal{A},P)$, the three random variables/vectors $\mb{x}$, $\mathcal{Q}$ and $\mb{u}$ are defined as the following three measurable functions: $\mb{x}: (\Omega, \mathcal{A})\rightarrow (\mathbb{R}^m,\mathfrak{B}(\mathbb{R}^m))$, $\mathcal{Q}: (\Omega, \mathcal{A})\rightarrow (\mathbb{R}^+,\mathfrak{B}(\mathbb{R}^+))$, $\mb{u}: (\Omega, \mathcal{A})\rightarrow (S_{\mathbb{R}}^{m-1},\mathfrak{B}(S_{\mathbb{R}}^{m-1}))$. Note that $\mathfrak{B}(S_{\mathbb{R}}^{m-1})$ is the Borel $\sigma$-algebra on the unit sphere $S_{\mathbb{R}}^{m-1}$ induced by the subspace topology from $\mathbb{R}^m$.} Moreover, $\mathcal{Q}$ has pdf given by:
\begin{equation}
\label{Q_pdf}
p_{\mathcal{Q}}(q) = \delta_m^{-1} q^{m/2-1} g (q).
\end{equation}

It is immediate to verify that the definition of elliptical density suffers from a lack of indentifiability for the couple $(\bs{\Sigma},g)$. Specifically, we can easily note that $RES_{m}(\bs{\mu},\bs{\Sigma},g(t)) \equiv RES_{m}(\bs{\mu},c\bs{\Sigma},c^{m/2}g(ct)), \forall c>0$. To avoid this ambiguity, we may decide to put a constraint on the \virg{functional form} of the density generator $g$. In particular, we force the density generator to belong to the following set:  
\begin{equation}
	\label{set_G_c}
	\overline{\mathcal{G}} = \graffe{ \bar{g} \in \mathcal{G} \left| \delta_m^{-1} \int_{0}^{\infty}q^{m/2}\bar{g}(q)dq= E\{\mathcal{Q}\} = m \right. }.
\end{equation}
Note that, from the stochastic representation \eqref{SRT_dec}, the properties of $\mb{u}\sim \mathcal{U}(S_{\mathbb{R}}^{m-1})$ and the fact that $\mathcal{Q}$ is independent of $\mb{u}$, the constraint $E\{\mathcal{Q}\} = m$ is equivalent to:
\begin{equation}\label{scatter_cov}
	E\{(\mb{x}-\bs{\mu})(\mb{x}-\bs{\mu})^T\} = E\{\mathcal{Q}\} \bs{\Sigma}^{1/2} E\{\mb{u}\mb{u}^T\}\bs{\Sigma}^{T/2} = \bs{\Sigma},
\end{equation}
(where $\bs{\Sigma}^{T/2} = (\bs{\Sigma}^{1/2})^T = \bs{\Sigma}^{1/2}$) i.e.\ the scatter matrix can be directly interpreted as the usual covariance matrix. It is important to highlight that, as shown in \cite{Hallin_P_2006,PAINDAVEINE}, it is possible to choose more general constraints that do not require the finiteness of any moment of the 2nd-order modular variate $\mathcal{Q}$, unlike the one in \eqref{set_G_c} used in this work.

Let us now introduce the matrix scale function:
\begin{equation}\label{scale_f}
	\begin{split}
		S :\;& \mathcal{S}_m^\mathbb{R} \rightarrow \mathbb{R}^{+}\\
		& \bs{\Sigma} \mapsto S(\bs{\Sigma}) = s
	\end{split}
\end{equation}
satisfying the following assumptions \cite{Hallin_P_2006,PAINDAVEINE}:
\begin{itemize}
	\item[A1] Homogeneity of order one: $S(c \cdot \bs{\Sigma}) = c \cdot S(\bs{\Sigma}), \forall c>0$,
	\item[A2] Differentiability over $\mathcal{S}_m^\mathbb{R}$ with $\parder{S(\bs{\Sigma})}{[\bs{\Sigma}]_{11}} \ne 0$,
	\item[A3] $S(\mb{I}_m) = 1$.
\end{itemize}

Then we define the \textit{shape} matrix $\mb{V}_S$ as:
\begin{equation}
	\label{shape_m}
	\mb{V}_S \triangleq \bs{\Sigma}/S(\bs{\Sigma} ) \in \mathcal{S}_{m,S}^\mathbb{R},
\end{equation}
where $\mathcal{S}_{m,S}^\mathbb{R}$ is a (non-linear) differentiable manifold on dimension $m(m+1)/2-1$ such that:
\begin{equation}\label{man_S}
	\mathcal{S}_{m,S}^\mathbb{R} \triangleq \{\mb{V}_S \in \mathcal{S}_{m}^\mathbb{R}|S(\mb{V}_S)=1\}. 
\end{equation}
We note in passing that the most popular choices for the \textit{scale} function $S(\cdot)$ are $S(\bs{\Sigma}) = [\bs{\Sigma}]_{11}$, $S(\bs{\Sigma}) = \trace{\bs{\Sigma}}/m$ and $S(\bs{\Sigma}) = |\bs{\Sigma}|^{1/m}$. It is worth emphasizing that, under A1, A2 and A3, and as a direct consequence of the implicit function theorem, the first top-left entry of $\mb{V}_S$, i.e. $[\mb{V}_{S}]_{11}$, can (locally at $\mb{V}_S$) be expressed as function of the other entries. For example, if we choose $S(\bs{\Sigma}) = [\bs{\Sigma}]_{11}$, then $[\mb{V}_{S}]_{11}$ is trivially given by $[\mb{V}_{S}]_{11} = 1$. For $S(\bs{\Sigma}) = \trace{\bs{\Sigma}}/m$, we have that $[\mb{V}_{S}]_{11} = m-\sum\nolimits_{i=2}^m[\mb{V}_{S}]_{ii}$. Lastly, for $S(\bs{\Sigma}) = |\bs{\Sigma}|^{1/m}$ it can be easily shown that, using the Laplace's expansion of the determinant along the first row of $\mb{V}_S$, $[\mb{V}_{S}]_{11}$ can be recovered as $[\mb{V}_{S}]_{11} = \frac{1}{C_{11}}\tonde{1-\sum\nolimits_{i=2}^m(-1)^{1+i}[\mb{V}_{S}]_{1i}C_{1i}}$ where $C_{ij}$ indicates the cofactor of $[\mb{V}_{S}]_{ij}$.

As discussed in \cite{Hallin_P_2006}, the RES model can be parameterized in two different, yet equivalent, ways:
\begin{equation}\label{RES_full_par_model_Sigma}
	\mathcal{P}_{\bs{\gamma},g} = \graffe{p_X(\mb{x}|\bs{\gamma},g) = |\mb{V}_{S}|^{-1/2}  g \left((\mb{x}-\bs{\mu})^T\mb{V}_{S}^{-1}(\mb{x}-\bs{\mu}) \right); \bs{\gamma} \in \Gamma, g \in \mathcal{G} },
\end{equation}
and
\begin{equation}\label{RES_full_par_model_V}
	\mathcal{P}_{\bs{\eta},\bar{g}} = \graffe{p_X(\mb{x}|\bs{\eta},\bar{g}) = s^{-m/2}|\mb{V}_S|^{-1/2}  \bar{g} \left(s^{-1}(\mb{x}-\bs{\mu})^T\mb{V}_S^{-1}(\mb{x}-\bs{\mu}) \right); \bs{\eta} \in \Phi, \bar{g} \in \overline{\mathcal{G}} }
\end{equation}
where the finite-dimensional parameter of the first parameterization \eqref{RES_full_par_model_Sigma} makes use of the scatter matrix:
\begin{equation}\label{nu_Sigma}
	\bs{\gamma} \triangleq (\bs{\mu}^T,\ovecs{\mb{V}_S}^T)^T \in \Gamma \subseteq \mathbb{R}^m \times \ovecs{\mathcal{S}_{m,S}^\mathbb{R}},
\end{equation}
while the one of the second parameterization \eqref{RES_full_par_model_V} is based on the shape matrix and on the scale:
\begin{equation}\label{eta_Sigma}
	\bs{\eta} \triangleq (\bs{\mu}^T,\ovecs{\mb{V}_S}^T,s)^T \in \Phi \subseteq \mathbb{R}^m  \times \ovecs{\mathcal{S}_{m,S}^\mathbb{R}} \times \mathbb{R}^+.
\end{equation}
In both parameterizations, the infinite-dimensional parameter is the density generator, specifically its unconstrained version $g \in \mathcal{G}$ in \eqref{RES_full_par_model_Sigma} and the constrained one $\bar{g} \in \overline{\mathcal{G}}$ in \eqref{RES_full_par_model_V}. Moreover, for any scale function $S$, we have that and $g(t)$ and $\bar{g}(t)$ are related by $g(t) = s^{-m/2}\bar{g}(s^{-1}t)$, $t \in \mathbb{R}^{+}$.

It is important to note that the two parameterizations in terms of $\ovecs{\mb{V}_S}$ are of interest for both practical and theoretical reasons. From a practical standpoint, in many applications such as Principal Component Analysis (PCA), Canonical Correlation Analysis (CCA) and subspace-based methods, only a scaled version of the covariance matrix, i.e. $\mb{V}_S$, is of interest while the scale term $s$ can be treated as a \textit{nuisance} parameter. Moreover, from a theoretical viewpoint, the parameterization in \eqref{RES_full_par_model_V} allows us to investigate the hidden and counter-intuitive relationships between the parametric and semiparametric efficiency on the shape matrix $\mb{V}_S$ and the scale $s$.

\subsection{Parametric and non-parametric sub-models of $\mathcal{P}_{\bs{\gamma},g}$ and $\mathcal{P}_{\bs{\eta},\bar{g}}$} 

For further reference, we summarize below the parametric and non-parametric sub-models of the two RES-models in \eqref{RES_full_par_model_Sigma} and \eqref{RES_full_par_model_V} that will be extensively used in the rest of this paper. For the sake of clarity, we will indicate with a superscript $u$ the \virg{unconstrained} sub-models, i.e. the ones involving an unconstrained density generator $g \in \mathcal{G}$, and with a superscript $c$ the \virg{constrained} sub-models, i.e. the ones involving a constrained $\bar{g} \in \overline{\mathcal{G}}$.

For the model $\mathcal{P}_{\bs{\eta},\bar{g}}$ in \eqref{RES_full_par_model_V}, we may identify three parametric sub-models as: 
\begin{equation}\label{P1_eta}
	\mathcal{P}_{\bs{\gamma}}^c = \graffe{p_X(\mb{x}|\bs{\gamma},s_0,\bar{g}_0): \bs{\gamma} \in \Gamma},
\end{equation}
\begin{equation}\label{P2_eta}
	\mathcal{P}^c_{s} = \graffe{p_X(\mb{x}|\bs{\gamma}_0 ,s,\bar{g}_0): s\in \mathbb{R}^+},
\end{equation}
\begin{equation}\label{P12_eta}
	\mathcal{P}^c_{\bs{\gamma},s} = \graffe{p_X(\mb{x}|\bs{\gamma},s,\bar{g}_0): \bs{\gamma} \in \Gamma,  s\in \mathbb{R}^+},
\end{equation}
and a non-parametric sub-model as:
\begin{equation}\label{con_non_par} 
	\mathcal{P}^c_{\bar{g}} = \graffe{p_X(\mb{x}|\bs{\gamma}_0,s_0,\bar{g}): \bar{g} \in \overline{\mathcal{G}}}.
\end{equation}

We may also define a parametric sub-model of $\mathcal{P}_{\bs{\eta},\bar{g}}$ that is largely exploited in applications. In fact, since the constrained density generator $\bar{g} \in \overline{\mathcal{G}}$ is generally parameterized by a set of $p$ parameters, say $\bs{\zeta} \in \Delta \subseteq \mathbb{R}^p$ (for example the degrees of freedom of the $t$-distribution or the shape parameter of the Generalized Gaussian distribution),
one may assume that functional form of $\bar{g} \in \overline{\mathcal{G}}$ is known a-priori, up to its parameter vector $\bs{\zeta} \in \Delta$. This leads to the definition of the following parametric sub-model:   
\begin{equation}\label{P_124}
	\mathcal{P}^c_{\bs{\gamma},s,\bs{\zeta}} = \graffe{ p_X(\mb{x}|\bs{\gamma},s,\bs{\zeta}) : \bs{\gamma} \in \Gamma,  s\in \mathbb{R}^+, \bs{\zeta} \in \Delta }
\end{equation} 
where 
\begin{equation}
	p_X(\mb{x}|\bs{\gamma},s,\bs{\zeta}) = s^{-m/2}|\mb{V}_S|^{-1/2}  \bar{g}_{\bs{\zeta}} \tonde{s^{-1}(\mb{x}-\bs{\mu})^T\mb{V}_S^{-1}(\mb{x}-\bs{\mu})},
\end{equation}
$\bar{g}_{\bs{\zeta}} \in \overline{\mathcal{G}}$ and $\bs{\zeta} \in \Delta$ represents a second finite-dimensional nuisance term. 

For model $\mathcal{P}_{\bs{\gamma},g}$ in \eqref{RES_full_par_model_Sigma}, we can identify the parametric sub-model 
\begin{equation}\label{P1u_eta}
	\mathcal{P}^u_{\bs{\gamma}} = \graffe{p_X(\mb{x}|\bs{\gamma},g_0): \bs{\gamma} \in \Gamma}
\end{equation}
and an \textit{unconstrained} non-parametric sub-model as:
\begin{equation}\label{unc_non_par}
	\mathcal{P}^u_{g} = \graffe{p_X(\mb{x}|\bs{\gamma}_0,g): g \in \mathcal{G}}.
\end{equation}

This framework will allow us to apply Lemmas \ref{Lem1} and \ref{Lem2} to reveal the counter-intuitive links (in terms of efficiency) between the different RES sub-models. The first step, which enables us to analyze parametric and semiparametric efficiency in the RES model, is the derivation of the score vectors of the finite-dimensional parameters, i.e. the location vector $\bs{\mu} \in \mathbb{R}^m$, the shape matrix $\mb{V}_S \in \mathcal{S}_{m,S}^\mathbb{R}$ (or of its vectorized version) and the scale $s \in \mathbb{R}^+$.

\subsection{The score vectors  $\mb{s}_{\bs{\mu}_0}$, $\mb{s}_{\vsVso}$ and $s_{s_0}$ in $\mathcal{P}^c_{\bs{\gamma},s}$}

Calculating the score vectors $\mb{s}_{\vsVso}$ and $s_{s_0}$ is not an easy task because the differentiable manifold $\mathcal{S}_{m,S}^\mathbb{R}$ in \eqref{man_S} is not a linear subspace of the set of real matrices $\mathbb{R}^{m \times m}$. On the contrary, since $\mathcal{S}_{m}^\mathbb{R}$ is indeed a linear subspace of $\mathbb{R}^{m \times m}$, the score vector $\mb{s}_{\vsSs}$ for the scatter matrix $\bs{\Sigma}$ is easy to derive and well-known in the literature (see e.g. \cite{For_SCRB, Chap_sem}). Fortunately, due to the differentability (under A1, A2 and A3) of $\mathcal{S}_{m,S}^\mathbb{R}$, there exists a diffeomorphism, say $\mb{w}$, between $(\mb{V}_S,s)$ and $\bs{\Sigma}$. Let us first provide an explicit expression for such diffeomorphism. It is important to highlight that most of the calculations that follow have already been derived in \cite{Hallin_P_2006,PAINDAVEINE}. In reporting them here, our aim is to provide a treatment consistent with the notations used in this article. We believe this will greatly facilitate reading and understanding.

By recalling that $\bs{\Sigma}=s\mb{V}_S$, under A1, A2 and A3, we can move from the parameterization based on $(\bs{\mu},\mb{V}_S,s)$ to the one based on $(\bs{\mu},\bs{\Sigma})$ and vice versa by means of
\begin{equation}\label{diff_par}
	\begin{split}
		\mb{w}:\;& \Phi \rightarrow \Omega \subseteq \mathbb{R}^m \times \vecs{\mathcal{S}_{m}^\mathbb{R}}\\
		& \bs{\eta} = (\bs{\mu}^T,\ovecs{\mb{V}_S}^T,s)^T \mapsto \mb{w}(\bs{\eta}) = (\bs{\mu}^T,s\cdot[\mb{V}_S]_{11}, s\cdot\ovecs{\mb{V}_S}^T)^T,
	\end{split}  
\end{equation}
whose Jacobian matrix $\mb{J}[\mb{w}](\bs{\eta}_0)$, evaluated at $\bs{\eta}_0 \in \Phi$, is given by:
\begin{equation}\label{Jac_d}
	\mb{J}[\mb{w}](\bs{\eta}_0) = \tonde{
		\begin{array}{ccc}
			\mb{I}_m & \mb{0}  & \mb{0} \\
			 \mb{0} & s_0\bs{\nabla}_{\ovecs{\mb{V}_S}}^T [\mb{V}_{S,0}]_{11} &  [\mb{V}_{S,0}]_{11} \\
			 \mb{0} & s_0\mb{I}_{m(m+1)/2-1} &  \ovecs{\mb{V}_{S,0}} \\
		\end{array}
	}\triangleq
	\tonde{
		\begin{array}{cc}
			\mb{I}_m & \mb{0}  \\
			\mb{0} & \mb{J}[\mb{w}](\ovecs{\mb{V}_{S,0}},s_0)  \\
		\end{array}
	}
\end{equation}
where $s_0 = S(\bs{\Sigma}_0)$ and $\mb{V}_{S,0} = \bs{\Sigma}_0/s_0$. Furthermore, the 2nd diagonal component in \eqref{Jac_d}, as shown in \cite[Sect. 4]{PAINDAVEINE} and, by using our notation in the Appendix \ref{App4}, of this work is given by
\begin{equation}
	\begin{split}
		\bs{\nabla}_{\ovecs{\mb{V}_S}} [\mb{V}_{S,0}]_{11} = - \left. \frac{\bs{\nabla}_{\ovecs{\mb{V}_S}}S([\mb{V}_{S}]_{11},\ovecs{\mb{V}_S})}{\partial S([\mb{V}_{S}]_{11},\ovecs{\mb{V}_S})/\partial [\mb{V}_{S}]_{11}} \right|_{\mb{V}_S = \mb{V}_{S,0}} ,
	\end{split}
\end{equation}
where $\bs{\nabla}_{\ovecs{\mb{V}_S}}$ indicates the \virg{constrained} gradient w.r.t. $\ovecs{\mb{V}_S} \in \ovecs{\mathcal{S}_{m,S}^\mathbb{R}}$ and $S([\mb{V}_{S}]_{11},\ovecs{\mb{V}_S})$ denotes the scale function applied to $\mb{V}_{S}$ reconstructed from $[\mb{V}_{S}]_{11}$ and $\ovecs{\mb{V}_S}$. For further reference, let us recast $\mb{J}[\mb{w}](\ovecs{\mb{V}_{S,0}},s_0)$ as:
\begin{equation}\label{Jac_Vs}
	\mb{J}[\mb{w}](\ovecs{\mb{V}_{S,0}},s_0) = \quadre{s_0\mb{K}_{\mb{V}_{S,0}} \; \vecs{\mb{V}_{S,0}} },
\end{equation}
where $\mb{K}_{\mb{V}_{S}}$ is a block-matrix defined as:
\begin{equation}\label{K_S}
	\mb{K}_{\mb{V}_{S}} \triangleq \quadre{
		\begin{array}{c}
		\nabla_{\ovecs{\mb{V}_S}}^T [\mb{V}_{S}]_{11} \\
		\mb{I}_{m(m+1)/2-1}
	\end{array}},
\end{equation}
evaluated at $\mb{V}_{S,0}$. 
As an example, we show in Appendix \ref{App5} that the term $\nabla_{\ovecs{\mb{V}_S}}^T [\mb{V}_{S}]_{11}$ can be explicitly obtained for the three above-mentioned scale functions as:
\begin{itemize}
	\item $\mb{0}^T_{m(m+1)/2-1}$ for $S(\bs{\Sigma}) = [\bs{\Sigma}]_{11}$,
	\item $-\frac{\cvec{\mb{I}_m}^T\mb{D}_m\underline{\mb{I}}_m^T}{(\cvec{\mb{I}_m}^T\mb{D}_m)_1}$ for $S(\bs{\Sigma}) = \trace{\bs{\Sigma}}/m$,
	\item $-\frac{\cvec{\mb{V}_S^{-1}}^T\mb{D}_m\underline{\mb{I}}_m^T}{(\cvec{\mb{V}_S^{-1}}^T\mb{D}_m)_1}$ for $S(\bs{\Sigma}) = |\bs{\Sigma}|^{1/m}$.
\end{itemize} 

The derivation of the gradient $\mb{J}[\mb{w}](\bs{\eta}_0)$ in \eqref{Jac_d} can be also found in \cite{PAINDAVEINE,Hallin_P_2006} as function of the matrix $\mb{M}_S^{\mb{V}_S}$ which can be linked to our results through the relation:
\begin{equation}\label{M_mat}
	\mb{M}_S^{\mb{V}_{S,0}} = \mb{K}_{\mb{V}_{S,0}}^T\mb{D}_m^T.
\end{equation}
Due to the crucial importance of this matrix in the subsequent sections, we collect in the Appendix \ref{App6} both its known properties (from \cite{PAINDAVEINE,Hallin_P_2006}) and some new ones, which will be extensively used in the reminder of the paper. 

Let us now introduce the score vectors $\mb{s}_{\bs{\mu}_0}$ and $\mb{s}_{\vsSs}$. By applying the rules of the differential matrix calculus detailed in \cite[Ch. 8]{Magnus}, the Kronecker product and the vec operator \cite[Ch. 3]{Magnus}, and the stochastic representation theorem \eqref{SRT_dec}, we get the following expressions of the scores (\cite[eqs. (34) and (37)]{For_SCRB},\cite[Sec. 5.3.2.1]{Chap_sem}):
\begin{equation}\label{s_mu}
	\mb{s}_{\bs{\mu}_0} =  \sqrt{\mathcal{Q}} \bar{\varphi}_0(\mathcal{Q}) \bs{\Sigma}^{-1/2}_{0} \mb{u},
\end{equation}
\begin{equation}\label{s_Sigma}
	\mb{s}_{\vsSs} = 2^{-1}\mb{D}_m^T[\kronsigmatinvm]\cvec{\mathcal{Q}\bar{\varphi}_0(\mathcal{Q})  \mb{u}\mb{u}^T - \mb{I}_m},
\end{equation}
where, following the notation adopted in \cite{Chap_back_JP}, we defined: \footnote{To avoid confusion, it is important to note that the function $\bar{\varphi}_0(t)$ is linked to the equivalent function $\psi_0(t)$ used in \cite{For_SCRB,For_SCRB_complex} by the constant $-2$: $\bar{\varphi}_0(t) = -2\psi_0(t)$.}
\begin{equation}\label{psi}
	\bar{\varphi}_0(t) \triangleq \frac{-2}{\bar{g}_0(t)}\frac{d\bar{g}_0(t)}{dt}.
\end{equation}

By using the diffeomorphism in \eqref{diff_par} and its Jacobian in \eqref{Jac_d}, the score vector for $\bs{\eta}_0 \triangleq (\bs{\gamma}_0^T,s_0)^T \in \Phi$ in $\mathcal{P}_{\bs{\eta},\bar{g}}$ in \eqref{RES_full_par_model_V} can be obtained as:
\begin{equation}
	\mb{s}_{\bs{\eta}_0} = \tonde{\begin{array}{c}
			\mb{s}_{\bs{\mu}_0}\\
			\quadre{\mb{J}[\mb{w}](\ovecs{\mb{V}_{S,0}},s_0)}^T \mb{s}_{\vsSs}
	\end{array}}.
\end{equation}

From the block-diagonal structure $\mb{J}[\mb{w}](\bs{\eta}_0)$ in \eqref{Jac_d}, we immediately have that the score vector $\mb{s}_{\bs{\mu}_0}$ for the location $\bs{\mu}_0$ is equal to \eqref{s_mu}. Moreover, from the block structure of $\mb{J}[\mb{w}](\ovecs{\mb{V}_{S,0}},s_0)$, by the definition of the matrix $\mb{M}_S^{\mb{V}_{S,0}}$ in \eqref{M_mat}, and by using the fact that:
\begin{equation}
	\vecs{\mb{V}_{S,0}}^T\mb{D}_m^T[\kronVsoinvm]\cvec{\mb{A}} = \trace{\mb{A}}
\end{equation}
for any $m \times m$ symmetric matrix $\mb{A}$, we have that:
\begin{equation}\label{s_V}
	\mb{s}_{\vsVso} = 2^{-1}\mb{M}_S^{\mb{V}_{S,0}}[\kronVsoinvm]\cvec{\mathcal{Q}\bar{\varphi}_0(\mathcal{Q}) \mb{u}\mb{u}^T - \mb{I}_m}
\end{equation}
and
\begin{equation}\label{s_s}
	s_{s_0} = (2s_0)^{-1} \tonde{\mathcal{Q}\bar{\varphi}_0(\mathcal{Q}) - m}.
\end{equation}
We note in passing that the expression of $\mb{s}_{\vsVso}$ in \eqref{s_V} was already derived in \cite[eq. (5)]{PAINDAVEINE} whereas the expression for $s_{s_0}$ is new.
Now that we have the explicit expression for the score vector $\mb{s}_{\bs{\eta}_0}$ in the model $\mathcal{P}_{\bs{\eta},\bar{g}}$ we can easily evaluate the related FIM $\mb{I}_{\bs{\eta}_0} \triangleq E_0\{\mb{s}_{\bs{\eta}_0}\mb{s}_{\bs{\eta}_0}^T\}$. Before computing FIM, we define two functionals of the true density generator $\bar{g}_0$:
\begin{equation}\label{def alpha}
	\alpha(\bar{g}_0) \triangleq \frac{E\{\mathcal{Q}^2\bar{\varphi}^2_0(\mathcal{Q})\}}{m(m+2)},
\end{equation}
\begin{equation}\label{def beta}
	\beta(\bar{g}_0) \triangleq \frac{E\{\mathcal{Q}\bar{\varphi}^2_0(\mathcal{Q})\}}{m}.
\end{equation}
In the rest of this work, we will always assume that the density generator $\bar{g}$ is a $C^1$-function satisfying the regularity condition $\int_{0}^{\infty}\bar{g}(t)^{-1}(\bar{g}(t)/dt)^2t^{m/2+1}dt  <\infty$. This guarantees that $\alpha(\bar{g}_0)$ and $\beta(\bar{g}_0)$ are both finite. It is important to note that this condition corresponds to the Assumption 2.2 in \cite{Hallin_P_2006} which gives sufficient conditions to assure that the elliptical distributions are a locally asymptotically normal (LAN) family.

Then, from standard calculations and by using the independence between $\mathcal{Q}$ and $\mb{u}$ (along with the properties of $\mb{u}$), the FIM for $\bs{\eta}_0 \in \Phi$ is given by:
\begin{equation}\label{Par_model_eta_block}
	\mb{I}_{\bs{\eta}_0} =\tonde{
		\begin{array}{cc}
			\mb{I}_{\bs{\gamma}_0} & \mb{I}_{\bs{\gamma}_0,s_0} \\
			\mb{I}_{\bs{\gamma}_0,s_0}^T & I_{s_0}
		\end{array}
	}= \tonde{
	\begin{array}{ccc}
		\mb{I}_{\bs{\mu}_0} & \mb{I}_{\bs{\mu}_0,\vsVso} & \mb{I}_{\bs{\mu}_0,s_0} \\
					\mb{I}_{\bs{\mu}_0,\vsVso}^T	& \mb{I}_{\vsVso} & \mb{I}_{\vsVso, s_0} \\
					\mb{I}_{\bs{\mu}_0,s_0}^T & \mb{I}^T_{\vsVso, s_0} & I_{s_0}
	\end{array}
	},
\end{equation}
where \footnote{Note that, in the calculations, we used the equality $E\{\mathcal{Q}\bar{\varphi}_0(\mathcal{Q})\}=m$. In fact, from \eqref{psi}, \eqref{Q_pdf} and \eqref{set_G}, we have $E\{\mathcal{Q}\bar{\varphi}_0(\mathcal{Q})\} = -2 \delta_m^{-1} \int_{0}^{\infty}q^{m/2}d\bar{g}
	=  -2 \delta_m^{-1}[q^{m/2}\bar{g}(q)]_0^{\infty} + m\delta_m^{-1}\int_{0}^{\infty}q^{m/2-1}\bar{g}(q)dq=m$. Similar calculation gives $E\{\mathcal{Q}^2\bar{\varphi}_0(\mathcal{Q})\}=m(m+2)$ under the additional condition $\lim_{q  \rightarrow \infty}q^{m/2+1}\bar{g}(q)=0$. 
	}
\begin{equation}\label{I_mu}
	\mb{I}_{\bs{\mu}_0} \triangleq E_0\{\mb{s}_{\bs{\mu}_0}\mb{s}_{\bs{\mu}_0}^T\} =  \beta(\bar{g}_0) [s_0\mb{V}_{S,0}]^{-1},
\end{equation}
\begin{equation}\label{I_mu_V}
	\mb{I}_{\bs{\mu}_0,\vsVso} \triangleq E_0\{\mb{s}_{\bs{\mu}_0}\mb{s}_{\vsVso}^T\} = \mb{0},
\end{equation}
\begin{equation}\label{I_mu_s}
\mb{I}_{\bs{\mu}_0,s_0} \triangleq E_0\{\mb{s}_{\bs{\mu}_0}s_{s_0}\} = \mb{0},
\end{equation}
\begin{equation}
	\begin{split}\label{I_V}
		\mb{I}_{\vsVso} & \triangleq E_0\{\mb{s}_{\vsVso}\mb{s}_{\vsVso}^T\} =\\
		&= \frac{1}{4}\mb{M}_S^{\mb{V}_{S,0}}[2\alpha(\bar{g}_0)(\mb{V}_{S,0}^{-1}\otimes\mb{V}_{S,0}^{-1})+(\alpha(\bar{g}_0) -1){\rm vec}({\bf V}_{S,0}^{-1}){\rm vec}^T({\bf V}_{S,0}^{-1})][\mb{M}_S^{\mb{V}_{S,0}}]^T,
	\end{split}
\end{equation}
\begin{equation}\label{I_s}
	\begin{split}
		I_{s_0} \triangleq E_0\{(s_{s_0})^2\} =\frac{m(m+2)\alpha(\bar{g}_0) -m^2}{4s_0^2 },
	\end{split}
\end{equation}
\begin{equation}\label{I_Vs}
	\mb{I}_{\vsVso, s_0} \triangleq E_0\{\mb{s}_{\vsVso}s_{s_0}\} = \frac{(m+2)\alpha(\bar{g}_0) -m}{ 4  s_0}\mb{M}_S^{\mb{V}_{S,0}}\vVsoinv.
\end{equation}

We note that these expressions coincide with those already given in \cite[eqs. (6) and (7)]{PAINDAVEINE}.

\subsection{The score vectors  $\mb{s}_{\bs{\mu}_0}$, $\mb{s}_{\vsVso}$  in $\mathcal{P}^u_{\bs{\gamma}}$}\label{sec_score_u}
Since this sub-model involves an unconstrained $g_0 \in \mathcal{G}$, we need to define an \virg{unconstrained} version of the function $\bar{\varphi}_0(t)$ in \eqref{psi} as:
\begin{equation}\label{psi_u}
	\varphi_0(t) = s_0^{-1}\bar{\varphi}_0(s_0^{-1}t).
\end{equation}
Consequently, from \eqref{Q_RES} and from the definition of the shape matrix in \eqref{shape_m}, it is immediate to verify that the score vectors $\mb{s}_{\bs{\mu}_0}$ and $\mb{s}_{\vsVso}$ in $\mathcal{P}^u_{\bs{\gamma}}$ are equal in distribution to the ones already derived in $\mathcal{P}^c_{\bs{\gamma},s}$ whose expressions are given in \eqref{s_mu} and \eqref{s_V}. Moreover, we have that:
\begin{equation}
	\alpha(g_0) = \alpha(\bar{g}_0),\quad \beta(g_0)=s_0^{-1}\beta(\bar{g}_0).
\end{equation}

\section{The geometry of $\mathcal{P}_{\bs{\gamma},g}$ and $\mathcal{P}_{\bs{\eta},\bar{g}}$ and some results on parametric and semiparametric information matrices}\label{sec_geo_V_g}
We now state the key result that allows us to study the impact of finite and infinite-dimensional nuisance parameters on the semiparametric efficiency in $\mathcal{P}_{\bs{\gamma},g}$ and $\mathcal{P}_{\bs{\eta},\bar{g}}$ for the parameter vector of interest $\bs{\gamma}_0 = (\bs{\mu}_0^T,\vsVso^T)^T$. To this end, we have the following tools at our disposal:
\begin{enumerate}
	\item The score vector of $\bs{\gamma}_0$ in $\mathcal{P}^c_{\bs{\gamma},s}$ and $\mathcal{P}^u_{\bs{\gamma}}$
	\begin{equation}\label{score_gamma}
		\mb{s}_{\bs{\gamma}_0} = \Big(\mb{s}_{\bs{\mu}_0}^T,\mb{s}_{\vsVso}^T\Big)^T,
	\end{equation}
	where $\mb{s}_{\bs{\mu}_0}$ and $\mb{s}_{\vsVso}$ are given in \eqref{s_mu} and \eqref{s_V}, respectively. Moreover, we recall here that $[\mb{s}_{\bs{\gamma}_0}]_i \in \mathcal{H}$.
	\item The (one-dimensional) nuisance tangent space $\mathcal{T}_{s_0}$ of the parametric sub-model $\mathcal{P}^c_{s}$ in \eqref{P2_eta} at $s_0$
	\begin{equation}
		\mathcal{T}_{s_0} = \mathrm{Span} \{s_{s_0}\} \subset \mathcal{H}
	\end{equation}
	and the relevant orthogonal projection derived form \eqref{proj_T2} as:
	\begin{equation}\label{proj_T2_s}
		\Pi(h|\mathcal{T}_{s_0}) =  E\{h s_{s_0}\}I_{s_0}^{-1}s_{s_0}, \; h \in \mathcal{H},
	\end{equation}
	where $I_{s_0}$ has been derived in \eqref{I_s}.  
	\item The (infinite-dimensional) nuisance tangent space $\mathcal{T}^u_{g_0}$ of the unconstrained non-parametric model $\mathcal{P}^u_{g}$ in \eqref{unc_non_par} at $g_0 \in \mathcal{G}$:
	\begin{equation}\label{Tu_g0_1}
		\begin{split}
			\mathcal{T}^u_{g_0} = \graffe{h \in \mathcal{H}| h\;\text{is $\sigma(\mathcal{Q})$-measurable}},
		\end{split}
	\end{equation}
	 where $\sigma(\mathcal{Q}) \subset \mathfrak{B}(\mathcal{X})$ is the sub-$\sigma$-algebra generated by the random variable $\mathcal{Q}$ in \eqref{Q_RES}. A detailed derivation of this result is given in Appendix \ref{proof_tan_spa}. From e.g. \cite[Ch. 23, Def. 4]{Jacod}, we have that the orthogonal projection of a generic element $h \in \mathcal{H}$  onto $\mathcal{T}_{g_0}^u$ can be obtained as: \footnote{It is important to note that for a generic $f \in L_2(P_0)$, this projection is $\Pi(f|\mathcal{T}^u_{g_0}) = E\{f|\mathcal{Q}\} - E\{f\}$.}
	\begin{equation}\label{proj_ug}
		\Pi(h|\mathcal{T}^u_{g_0}) = E\{h|\mathcal{Q}\}, \; \forall h \in \mathcal{H}. 
	\end{equation}
	\item The (infinite-dimensional) nuisance tangent space $\mathcal{T}^c_{\bar{g}_0}$ of the constrained non-parametric model $\mathcal{P}^c_{\bar{g}}$ in \eqref{con_non_par} at $\bar{g}_0 \in \overline{\mathcal{G}}$:
	\begin{equation}\label{Tc_g0_1}
		\begin{split}
			\mathcal{T}^c_{\bar{g}_0} = \graffe{l \in \mathcal{T}^u_{g_0}| \e{(\mathcal{Q}-m)l(\mathcal{Q})} = 0} = \graffe{l \in \mathcal{T}^u_{g_0}| \e{\mathcal{Q}l(\mathcal{Q})} = 0},
		\end{split}
	\end{equation}
	and then, as a direct consequence, we have that:
	\begin{equation}\label{direct_sum_Tu}
		\mathcal{H} \supseteq \mathcal{T}^u_{g_0} = \mathcal{T}^c_{\bar{g}_0} \oplus \mathrm{Span}\{\mathcal{Q}-m\}.
	\end{equation}
	A detailed proof of this result is given in Appendix \ref{proof_tan_spa} along with the derivation of the relevant orthogonal projection operator that can be expressed as:
	\begin{equation}\label{proj_cg}
		\Pi(h|\mathcal{T}^c_{\bar{g}_0})  = E\{h|\mathcal{Q}\} - \eg{\mathcal{Q}h}\sigma_{\mathcal{Q}}^{-2}(\mathcal{Q}-m), \quad \forall h\in \mathcal{H}.
	\end{equation}
	where $\sigma_{\mathcal{Q}} \triangleq \sqrt{\eg{(\mathcal{Q}-m)^2}}$. It is important to note here that to ensure that the decomposition in \eqref{direct_sum_Tu} and the projection operator in \eqref{proj_cg} are well-defined, we need to impose that $E\{\mathcal{Q}^2\} < \infty$ (along with, of course $E\{\mathcal{Q}\} =m$).
	\item From its expression given in \eqref{s_s} and the footnote 3, we verify that, when the score $s_{s_0}$ has zero-mean and is $\sigma(\mathcal{Q})$-measurable, then $\mathcal{T}_{s_0} \subset \mathcal{T}^u_{g_0}$. However, through some simple calculation, we can verify that $\e{\mathcal{Q}s_{s_0}(\mathcal{Q})} \neq 0$, then $\mathcal{T}_{s_0} \nsubseteq \mathcal{T}^c_{g_0}$. Finally, from \cite[Theo. 5, A.4]{BKRW} (see also \eqref{proj_sum} and the detailed explanation in Appendix \ref{App_Theo_sum}), we have that:
	\begin{equation}\label{proj_s_Tc}
		\Pi(h|\mathcal{T}_{s_0} + \mathcal{T}^c_{\bar{g}_0}) = \Pi(h|\mathcal{T}^c_{\bar{g}_0}) + \e{hr_0}(\e{r_0^2})^{-1}r_0 = \Pi(h|\mathcal{T}^u_{g_0}), \; \forall h \in \mathcal{H},
	\end{equation}
	where $r_0 \triangleq s_{s_0} - \Pi(s_{s_0}|\mathcal{T}^c_{\bar{g}_0})$ and where the second equality holds because $r_0 \in \mathrm{Span}\{\mathcal{Q}-m\}$. \footnote{Note that, from \eqref{s_s}, \eqref{proj_cg} and footnote 3, we have $\Pi(s_{s_0}|\mathcal{T}^c_{\bar{g}_0})  = E\{s_{s_0}|\mathcal{Q}\} - \eg{\mathcal{Q}s_{s_0}}\sigma_{\mathcal{Q}}^{-2}(\mathcal{Q}-m) = s_{s_0} - ms_0^{-1}\sigma_{\mathcal{Q}}^{-2}(\mathcal{Q}-m)$. Consequently, $r_0 = ms_0^{-1}\sigma_{\mathcal{Q}}^{-2}(\mathcal{Q}-m) \in \mathrm{Span}\{\mathcal{Q}-m\}$.}	
\end{enumerate}

We are now ready to establish an interesting result on the \textit{adaptivity property} of $\bs{\gamma}_0$ with respect to the infinite-dimensional nuisance $g_0$ (or $\bar{g}_0$) in the presence of the finite-dimensional nuisance $s_0$. This surprising adaptivity result was already discovered in \cite{PAINDAVEINE,Hallin_P_2006} as an application to the elliptical model of the general result in \cite[Prop. 2.6]{Hallin_Werker}. Here we will provide a more direct proof, based on the explicit form (as opposed to the asymptotic approximation used in \cite{PAINDAVEINE,Hallin_P_2006}) of the efficient semiparametric score vector $\bar{\mb{s}}_{\bs{\gamma}_0}$.

\begin{proposition}{(Restricted adaptivity property):}\label{Prop_sg} 
	Let us assume that $\mathcal{Q}$ in \eqref{Q_RES} has a finite second order moment, $\e{\mathcal{Q}^2}<\infty$, and that $\alpha(\bar{g}_0)$ in \eqref{def alpha} is finite as well. Then, the two efficient semiparametric FIMs $\bar{\mb{I}}(\bs{\gamma}_0|g_0)$ for the \virg{unconstrained} model $\mathcal{P}_{\bs{\gamma},g}$ in \eqref{RES_full_par_model_Sigma} and $\bar{\mb{I}}(\bs{\gamma}_0|s_0,\bar{g}_0)$ for the \virg{constrained} model $\mathcal{P}_{\bs{\eta},\bar{g}}$ in \eqref{RES_full_par_model_V} are both equal to the efficient parametric FIM $\bar{\mb{I}}(\bs{\gamma}_0|s_0)$ in $\mathcal{P}^c_{\bs{\gamma},s}$ in \eqref{P12_eta} that is a block diagonal matrix of the form
	\begin{equation}
		\bar{\mb{I}}(\bs{\gamma}_0|g_0) = \bar{\mb{I}}(\bs{\gamma}_0|s_0,\bar{g}_0) = \bar{\mb{I}}(\bs{\gamma}_0|s_0) = \mathrm{blkdiag}(\mb{I}_{\bs{\mu}_0},\bar{\mb{I}}_{\vsVso}),
	\end{equation}
	where $\mb{I}_{\bs{\mu}_0}$ is given in \eqref{I_mu} and 
		\begin{equation}\label{bar_I_V}
	\bar{\mb{I}}_{\vsVso} = 2^{-1}\alpha(\bar{g}_0)\mb{M}_S^{\mb{V}_{S,0}}[({\bf V}_{S,0}^{-1}\otimes{\bf V}_{S,0}^{-1})-m^{-1}{\rm vec}({\bf V}_{S,0}^{-1}){\rm vec}^T({\bf V}_{S,0}^{-1})][\mb{M}_S^{\mb{V}_{S,0}}]^T.
	\end{equation}
\end{proposition}

\begin{IEEEproof}
Let us start by proving that $\bar{\mb{I}}(\bs{\gamma}_0|g_0) = \bar{\mb{I}}(\bs{\gamma}_0|s_0)$. To this end, we just need to show that $\Pi(\mb{s}_{\bs{\gamma}_0}|\mathcal{T}^u_{g_0}) = \Pi(\mb{s}_{\bs{\gamma}_0}|\mathcal{T}_{s_0})$. For the location parameter $\bs{\mu}_0$, from the expression of $\Pi(\cdot|\mathcal{T}_{s_0})$ in \eqref{proj_T2_s} and from the fact that $\mb{I}_{\bs{\mu}_0,s_0} \triangleq E_0\{\mb{s}_{\bs{\mu}_0}s_{s_0}\} = \mb{0}$ as indicated in \eqref{I_mu_s}, we immediately have that:
\begin{equation}\label{proj_appo}
	\Pi(\mb{s}_{\bs{\mu}_0}|\mathcal{T}_{s_0}) = \mb{I}_{\bs{\mu}_0,s_0} I_{s_0}^{-1}s_{s_0} =\mb{0},
\end{equation} 
then $\mb{s}_{\bs{\mu}_0} \perp \mathcal{T}_{s_0}$. Moreover, from the expression of the score $\mb{s}_{\bs{\mu}_0}$ given by \eqref{s_mu}, the independence between $\mathcal{Q}$ and $\mb{u} \sim \mathcal{U}(S_{\mathbb{R}}^{m-1})$ and the property $E\{\mb{u}\}=\mb{0}$, we have from \eqref{proj_ug}:
\begin{equation}\label{score s}
	\Pi(\mb{s}_{\bs{\mu}_0}|\mathcal{T}^u_{g_0}) = E\{\mb{s}_{\bs{\mu}_0}|\mathcal{Q}\}= \mb{0},
\end{equation}
then $\mb{s}_{\bs{\mu}_0} \perp \mathcal{T}^u_{g_0}$.

Let us now proceed to the calculation of the projection of $\mb{s}_{\vsVso}$ onto the finite-dimensional tangent space $\mathcal{T}_{s_0}$. From  \eqref{proj_T2_s}, \eqref{I_Vs}, \eqref{I_s} and \eqref{score s}, we get:
\begin{equation}\label{proj_Vs}
	\Pi(\mb{s}_{\vsVso}|\mathcal{T}_{s_0}) =  \mb{I}_{\vsVso, s_0}I_{s_0}^{-1}s_{s_0} = (2m)^{-1}\mb{M}_S^{\mb{V}_{S,0}} \vVsoinv \tonde{\mathcal{Q}\bar{\varphi}_0(\mathcal{Q}) - m}.
\end{equation}
Moreover, from the expression of the score $\mb{s}_{\vsVso}$ given by \eqref{s_V} and by using the independence between $\mathcal{Q}$ and $\mb{u}$ and $E\{\mb{u}\mb{u}^T\}=m^{-1}\mb{I}_m$, we have:
\begin{equation}\label{s_Sigma b}
	\begin{split}
		\Pi(\mb{s}_{\vsVso}|\mathcal{T}^u_{g_0})&=E\{\mb{s}_{\vsVso}|\mathcal{Q}\} = 2^{-1}\mb{M}_S^{\mb{V}_{S,0}}[\kronVsoinvm]\tonde{\mathcal{Q}\bar{\varphi}_0(\mathcal{Q}) \mathrm{vec}(E\{\mb{u}\mb{u}^T|\mathcal{Q}\}) - \cvec{\mb{I}_m}}\\
		&= (2m)^{-1}\mb{M}_S^{\mb{V}_{S,0}}[\kronVsoinvm]\cvec{\mb{I}_m}\tonde{\mathcal{Q}\bar{\varphi}_0(\mathcal{Q}) - m}\\
		& = (2m)^{-1}\mb{M}_S^{\mb{V}_{S,0}} \vVsoinv \tonde{\mathcal{Q}\bar{\varphi}_0(\mathcal{Q}) - m} = \Pi(\mb{s}_{\vsVso}|\mathcal{T}_{s_0}),
	\end{split}
\end{equation}
where we recall that $\mathcal{T}_{s_0} \subset \mathcal{T}^u_{g_0}$.
As a consequence of the previous outcomes and of the relation in \eqref{proj_s_Tc}, we have that the efficient semiparametric score vector $\bar{\mb{s}}_{\bs{\gamma}_0}$ can be expressed in three equivalent forms: 
\begin{equation}\label{bar_s_0}
	\begin{split}
		\bar{\mb{s}}_{\bs{\gamma}_0} & \triangleq  \mb{s}_{\bs{\gamma}_0}- \Pi( \mb{s}_{\bs{\gamma}_0}|\mathcal{T}_{s_0}) \perp \mathcal{T}_{s_0} \\
		& =  \mb{s}_{\bs{\gamma}_0} - \Pi( \mb{s}_{\bs{\gamma}_0}|\mathcal{T}^u_{g_0})  \perp \mathcal{T}^u_{g_0}, \\
		& =  \mb{s}_{\bs{\gamma}_0} - \Pi( \mb{s}_{\bs{\gamma}_0}|\mathcal{T}_{s_0} + \mathcal{T}^c_{\bar{g}_0}),
	\end{split}
\end{equation}
and then $\bar{\mb{I}}(\bs{\gamma}_0|s_0) = \bar{\mb{I}}(\bs{\gamma}_0|g_0) = \bar{\mb{I}}(\bs{\gamma}_0|s_0,\bar{g}_0)$.

We note in passing that the equality $\bar{\mb{I}}(\bs{\gamma}_0|s_0)  = \bar{\mb{I}}(\bs{\gamma}_0|s_0,\bar{g}_0)$ can be obtained by exploiting Lemma \ref{Lem2}. Specifically, in this context, Condition \eqref{main_cond} of Lemma \ref{Lem2} can be expressed as:
\begin{equation}\label{bar_t}
 [\mb{s}_{\bs{\gamma}_0} - \Pi(\mb{s}_{\bs{\gamma}_0}|\mathcal{T}_{s_0})] \perp \mathcal{T}^c_{\bar{g}_0},
\end{equation}
that follows immediately from \eqref{bar_s_0} and from the fact that $\mathcal{T}^c_{\bar{g}_0} \subset \mathcal{T}^u_{g_0}$.

To conclude, the expression of $\bar{\mb{I}}_{\vsVso}$ can be derived from the efficient score $\bar{\mb{s}}_{\vsVso}$ in
\eqref{s_V} and \eqref{s_Sigma b} as:
\begin{equation}
	 \bar{\mb{s}}_{\vsVso} = 2^{-1}\mathcal{Q}\bar{\varphi}_0(\mathcal{Q})\mb{M}_S^{\mb{V}_{S,0}}[\kronVsoinvm]\cvec{\mb{u}\mb{u}^T - m^{-1}\mb{I}_m},
\end{equation}
from which \eqref{bar_I_V} of $\bar{\mb{I}}_{\vsVso} \triangleq  E\{\bar{\mb{s}}_{\vsVso}\bar{\mb{s}}_{\vsVso}^T \}$ follows directly.\footnote{The calculation is based on the independence between $\mathcal{Q}$ and $\mb{u}$ along with the identities 
 	$E\{{\rm vec}(\mb{u}\mb{u}^T){\rm vec}^T(\mb{u}\mb{u}^T)\}=\frac{1}{m(m+2)}\left({\bf I}_{m^2}+{\bf K}_{m}+{\rm vec}(\mb{I}_m){\rm vec}^T(\mb{I}_m)\right)$ and	$\mb{K}_m\mb{D}_m=\mb{D}_m$ \cite[p. 49]{Magnus} in \eqref{M_mat}}
We note, in passing that similar expressions have been derived in \cite[eqs. (8) and (9)]{PAINDAVEINE}, by using the already-recalled asymptotic representation of the projection operator.
\end{IEEEproof} 

As a side result of this proof and of the fact that $\mb{I}_{\bs{\mu}_0,\vsVso} = \mb{0}$ and $\mb{I}_{\bs{\mu}_0,s_0} = \mb{0}$, it is immediate to verify that knowing or not knowing the location vector $\bs{\mu}_0$ has no impact on the efficiency losses when we do parametric or semiparametric inference on the shape matrix $\mb{V}_{S,0} \in \mathcal{S}_{m,S}^\mathbb{R}$. Moreover, since $\Pi(\mb{s}_{\bs{\mu}_0}|\mathcal{T}_{g_0}^u) =\Pi(\mb{s}_{\bs{\mu}_0}|\mathcal{T}_{g_0}^c)= \mb{0}$, the lack of knowledge of the density generator does not induce any loss of efficiency w.r.t. $\bs{\mu}_0$ (see e.g. the discussion in \cite{For_SCRB}).

Roughly speaking, Proposition \ref{Prop_sg} tells us that when performing inference on the shape matrix $\mb{V}_{S,0} \in \mathcal{S}_{m,S}^\mathbb{R}$, not knowing the scale $s_0$ (while knowing the constrained density generator $\bar{g}_0$) or not knowing the density generator at all $g_0$, leads to the same efficiency losses.

\textbf{Remark 1} (\textit{Full adaptivity property}): In \cite{PAINDAVEINE}, it has been show that, if the scale function $S_d(\bs{\Sigma}_0) \triangleq |\bs{\Sigma}_0|^{1/m}$ is adopted, then the parametric FIM $\mb{I}_{\bs{\gamma}_0}$ in the parametric model $\mathcal{P}_{\bs{\gamma}}$ in \eqref{P1_eta} is equal to the efficient SFIM $\bar{\mb{I}}(\bs{\gamma}_0|s_0,\bar{g}_0)$ for $\mathcal{P}_{\bs{\eta},\bar{g}}$. In other worlds, if the scale $S_d(\bs{\Sigma}_0) = |\bs{\Sigma}_0|^{1/m}$ is adopted, when performing inference on $\mb{V}_{S,0} \in \mathcal{S}_{m,S_d}^\mathbb{R}$, knowing or not knowing the scale $s_0$ and/or the density generator $g_0$, does not lead to any asymptotic efficiency loss w.r.t. the case in which $s_0$ and/or $g_0$ are perfectly known. For this reason, the scale $S_d(\bs{\Sigma}_0)$ is called \virg{canonical} in \cite{PAINDAVEINE}. A sketch of the proof of this result is recalled in Appendix \ref{App_det_based}.

Let us now focus our attention on the parametric model $\mathcal{P}^c_{\gamma,s,\bs{\zeta}}$, whose sub-model is given in \eqref{P_124}, in which the finite-dimensional nuisance vector is represented by the parameters $\bs{\zeta}_0 \in \Delta$ of the unconstrained density generator.

\begin{corollary}\label{Prop_szeta} 
	The efficient parametric FIM $\bar{\mb{I}}(\bs{\gamma}_0|s_0)$ in $\mathcal{P}^c_{\bs{\gamma},s}$ in \eqref{P12_eta} is equal to the parametric efficient FIM $\bar{\mb{I}}(\bs{\gamma}_0|s_0, \bs{\zeta}_0)$ in $\mathcal{P}^c_{\gamma,s,\bs{\zeta}}$ in \eqref{P_124}.
\end{corollary}
\begin{IEEEproof}
	In Proposition \ref{Prop_sg}, we showed that $\mb{s}_{\bs{\gamma}_0} - \Pi(\mb{s}_{\bs{\gamma}_0}|\mathcal{T}_{s_0}) \perp \mathcal{T}^c_{\bar{g}_0}$. Moreover since, as discussed in Appendix \ref{proof_tan_spa}, $\mathcal{T}^c_{\bs{\zeta}_0} \subset \mathcal{T}^c_{\bar{g}_0}$, we trivially have that $\mb{s}_{\bs{\gamma}_0} - \Pi(\mb{s}_{\bs{\gamma}_0}|\mathcal{T}_{s_0}) \perp \mathcal{T}^c_{\bs{\zeta}_0}$. Consequently, the parametric efficient FIM $\bar{\mb{I}}(\bs{\gamma}_0|s_0, \bs{\zeta}_0)$ in $\mathcal{P}^c_{\gamma,s,\bs{\zeta}}$ is equal to the parametric efficient FIM $\bar{\mb{I}}(\bs{\gamma}_0|s_0)$ in $\mathcal{P}^c_{\bs{\gamma},s}$ in \eqref{P12_eta}.
\end{IEEEproof}

Roughly speaking, Corollary \ref{Prop_szeta} tells us that, if the scale $s_0$ is an unknown parameter, when performing inference on the shape matrix $\mb{V}_{S,0} \in \mathcal{S}_{m,S}^\mathbb{R}$, not knowing the parameters $\bs{\zeta}_0$ of the density generator $\bar{g}_{\bs{\zeta}_0}$ or not knowing the whole functional form $\bar{g}$, leads to the same efficiency losses.

Proposition \ref{Prop_sg} and Corollary \ref{Prop_szeta} give us a clear picture of the efficiency relationships among the various sub-models of elliptical distributions in terms of Fisher information matrices. It would now be interesting to obtain explicit closed-form expressions for the inverses of these FIM matrices in order to obtain the related information bounds. This will be the objective of the next section. 

\section{Parametric and semiparametric information bounds in $\mathcal{P}_{\bs{\gamma},g}$ and $\mathcal{P}_{\bs{\eta},\bar{g}}$}
\label{sec:CRB}

The purpose of this section is to provide the counterpart of the Propositions \ref{Prop_sg} in terms of information bounds.  

A clarification is in order before continuing. In the theory of parametric estimation, the best known information bound is the Cram\'er-Rao Bound (CRB) which coincides with inverse of the related FIM as shown e.g. \cite[Chap. 5, Sects 5 and 6]{Lehmann}. Consequently, for the different parametric sub-models of the semiparametric model in $\mathcal{P}_{\bs{\eta},\bar{g}}$ in \eqref{RES_full_par_model_V}, we can derive the related information/CR bounds as:
\begin{itemize}
	\item CRB for the parameters of interest $\bs{\gamma}_0 =(\bs{\mu}_0^T,\ovecs{\mb{V}_{S,0}}^T)^T$ in the presence of the finite-dimensional nuisance $s_0$ in the parametric sub-model $\mathcal{P}^c_{\bs{\gamma},s}$ in \eqref{P12_eta}:
	\begin{equation}\label{efficient FIM}
		\mathrm{CRB}(\bs{\gamma}_0|s_0) = \bar{\mb{I}}(\bs{\gamma}_0|s_0)^{-1}.
	\end{equation}
	\item CRB for the parameters of interest $\bs{\gamma}_0$ in the presence of two finite-dimensional nuisance terms, $s_0$ and $\bs{\zeta}_0$ in the parametric sub-model $\mathcal{P}^c_{\bs{\gamma},s,\bs{\zeta}}$ in \eqref{P_124}:
	\begin{equation}
		\mathrm{CRB}(\bs{\gamma}_0|s_0, \bs{\zeta}_0) = \bar{\mb{I}}(\bs{\gamma}_0|s_0, \bs{\zeta}_0)^{-1} = \bar{\mb{I}}(\bs{\gamma}_0|s_0)^{-1}
	\end{equation}
	as shown in Corollary \ref{Prop_szeta}.
\end{itemize}

When we move to the semiparametric case, the CRB can no longer be defined as in the classical parametric case. Remarkably, the H\'ajek-Le Cam convolution theorem (see e.g. \cite[Sect. 3.3, Theo. 2]{BKRW}) provides the right theoretical framework to unify the concept of information bound in the parametric, semiparametric and non-parametric case. A formal presentation of this theorem would lead us too far from the main purpose of this article. Therefore, below we will simply define the \virg{Semiparametric CRB (SCRB)} as the information bound obtained as the inverse of the efficient semiparametric FIM (SFIM) (for an in-depth discussion about this point, we refer the readers to \cite[Chap. 3]{BKRW}). Specifically, for the semiparametric model $\mathcal{P}_{\bs{\eta},\bar{g}}$ in \eqref{RES_full_par_model_V}, we have that the SCRB for the parameters of interest $\bs{\gamma}_0$ in the presence of a finite-dimensional $s_0$ and of an infinite-dimensional $\bar{g}_0 \in \overline{\mathcal{G}}$ nuisance terms is given by: 
\begin{equation}
	\mathrm{SCRB}(\bs{\gamma}_0|s_0, \bar{g}_0) = \bar{\mb{I}}(\bs{\gamma}_0|s_0, \bar{g}_0)^{-1},
\end{equation}
where $\bar{\mb{I}}(\bs{\gamma}_0|s_0, \bar{g}_0)$ is the efficient Semiparametric FIM as defined in \eqref{eff_sem_FIM}.

The following two subsections are organized as follows: in subsect. \ref{ex_CRB_gamma}, an explicit expression for the parametric CRB of $\bs{\gamma}_0$ in $\mathcal{P}^c_{\bs{\gamma},s}$ in the presence of the nuisance parameter $s_0$ is provided, while the subsect. \ref{prop_bounds} proposes the counterpart of the Proposition \ref{Prop_sg} in terms of parametric and semiparametric information bounds. 

\subsection{Explicit expression for $\mathrm{CRB}(\bs{\gamma}_0|s_0)$ in  $\mathcal{P}^c_{\bs{\gamma},s}$}\label{ex_CRB_gamma}
Let us consider the inverse, say $\mb{w}^{-1}$, of the diffeomorphism $\mb{w}$ given in \eqref{diff_par}. It can be explicitly obtained as:
\begin{equation}\label{diff_par b}
	\begin{split}
		\mb{w}^{-1}:\;&  \Omega \rightarrow\Phi \\
		& \bs{\nu} = (\bs{\mu}^T,\vecs{\bs{\Sigma}}^T)^T  \mapsto \mb{w}^{-1}(\bs{\nu}) = (\bs{\mu}^T,S^{-1}(\bs{\Sigma}) \cdot \ovecs{\bs{\Sigma}}^T, S(\bs{\Sigma}))^T,
	\end{split}  
\end{equation}
whose Jacobian matrix $\mb{J}[\mb{w}^{-1}](\bs{\nu}_0)$ is given by:
\begin{equation}\label{Jac_e}
	\mb{J}[\mb{w}^{-1}](\bs{\nu}_0) = \tonde{
		\begin{array}{cc}
			\mb{I}_m & \mb{0}  \\
			\mb{0} & S^{-1}(\bs{\Sigma}_0)\underline{\mb{I}}_m\mb{D}_m^{\#} \left[\mb{I}_{m^2}-\cvec{\mb{V}_{S,0}}\bs{\nabla}_{\cvec{\bs{\Sigma}}}^T S(\bs{\Sigma}_0)\right]\mb{D}_m\\
			\mb{0} & \bs{\nabla}_{\cvec{\bs{\Sigma}}}^T S(\bs{\Sigma}_0)\mb{D}_m  \\
		\end{array}
	}.
\end{equation}
The following proposition is proved in Appendix \ref{App_crb2}:
\begin{proposition}\label{Prop CRB eta}
The parametric CRB for $\bs{\gamma}_0$ in $\mathcal{P}^c_{\bs{\gamma},s}$ is given by
\begin{equation}\label{CRB_gamma_s}
		\mathrm{CRB}(\bs{\gamma}_0|s_0)
		\triangleq \bar{\mb{I}}(\bs{\gamma}_0|s_0)^{-1} =
		\tonde{
			\begin{array}{cc}
				\mathrm{CRB}(\bs{\mu}_0) & \mb{0}   \\
				\mb{0} & \mathrm{CRB}(\ovecs{\mb{V}_{S,0}}|s_0)  \\
			\end{array}
		},
\end{equation}
where:
\begin{equation}\label{CRB mu}
		\mathrm{CRB}(\bs{\mu}_0) \triangleq \mb{I}_{\bs{\mu}_0}^{-1} = \frac{s_0}{\beta(\bar{g}_0)} \mb{V}_{S,0}
\end{equation} 
and
\begin{equation}\label{CRB V}	
		\mathrm{CRB}(\ovecs{\mb{V}_{S,0}}|s_0) \triangleq \bar{\mb{I}}(\ovecs{\mb{V}_{S,0}}|s_0)^{-1}
		=
		\alpha(\bar{g}_0)^{-1} 
		\underline{\mb{I}}_m
		\mb{D}_m^{\#}\mb{P}_S(\mb{V}_{S,0})
		(\mb{I}_{m^2}+\mb{K}_{m})(\mb{V}_{S,0}\otimes \mb{V}_{S,0})
		\mb{P}_S^T(\mb{V}_{S,0})\mb{D}_m^{\#T}
		\underline{\mb{I}}_m^T,
\end{equation}
with
$\mb{P}_S(\mb{V}_{S,0}) \triangleq \mb{I}_{m^2}-{\rm vec}(\mb{V}_{S,0}) \bs{\nabla}_{{\rm vec}({\bs{\Sigma})}}^T S(\bs{\Sigma}_0)$, which takes the following expressions:
	\begin{itemize}
		\item $\mb{I}_{m^2}-{\rm vec}(\mb{V}_{S,0}){\bf e}^T_{1,m^2}$, for $S(\bs{\Sigma}) = [\bs{\Sigma}]_{11}$,
		\item $\mb{I}_{m^2}-\frac{1}{m}{\rm vec}(\mb{V}_{S,0}){\rm vec}(\mb{I}_{m})^T$, for $S(\bs{\Sigma}) = \trace{\bs{\Sigma}}/m$,
		\item $\mb{I}_{m^2}-\frac{1}{m}{\rm vec}(\mb{V}_{S,0}){\rm vec}(\mb{V}_{S,0}^{-1})^T$, for $S(\bs{\Sigma}) = |\bs{\Sigma}|^{1/m}$.
	\end{itemize}
Furthermore, for the scale function $S_d(\bs{\Sigma}) \triangleq |\bs{\Sigma}|^{1/m}$, 
the parameters $\mb{V}_{S,0}$ and $s_0$ are decoupled in the CRB for the parameter pair $(\mb{V}_{S,0},s_0)$
and
the matrix $\mathrm{CRB}(\ovecs{\mb{V}_{S,0}}|s_0)$ reduces to 
	\begin{equation}
		\label{CRB V d}
		\mathrm{CRB}(\ovecs{\mb{V}_{S_d,0}}|s_{d,0})
		=
		\alpha(\bar{g}_0)^{-1} 
		\underline{\mb{I}}_m
		\mb{D}_m^{\#}
		\quadre{
		(\mb{I}_{m^2}+\mb{K}_{m})(\mb{V}_{S,0}\otimes \mb{V}_{S,0})
		-2m^{-1}{\rm vec}(\mb{V}_{S,0}){\rm vec}(\mb{V}_{S,0})^T
		}
		\mb{D}_m^{\#^T}
		\underline{\mb{I}}_m^T.
	\end{equation}
\end{proposition}
We note that the proof of  Proposition \ref{Prop CRB eta} given in  Appendix \ref{App_crb2} also provides us with closed-form expressions of  
CRB$(\vecs{\bs{\Sigma}_0})$ (see \eqref{CRB vecs Sigma b}) and $\mathrm{CRB}(s_{0}|\ovecs{\mb{V}_{S_d,0}})$ (see \eqref{CRB s} \eqref{CRB s d}), which are also new results.

\subsection{Equality chains for parametric and semiparametric information bounds in RES distributions} \label{prop_bounds}
The following proposition summarizes the key points from the previous sections and can serve as a take-away message for readers who are not interested in full mathematical details. 

\begin{proposition}\label{Prop_chains}
	Let $\mathbb{R}^m \ni \mb{x}  \sim RES_{m}(\bs{\mu}_0,\bs{\Sigma}_0,\bar{g}_0)=RES_{m}(\bs{\mu}_0,\mb{V}_{S,0},g_0)$ be a RES distributed vector with location vector $\bs{\mu}_0$, scatter matrix $\bs{\Sigma}_0=s_0\mb{V}_{S,0}$ and (constrained or unconstrained) density generator $\bar{g}_0 \in \overline{\mathcal{G}}$ or $g_0 \in \mathcal{G}$, such that $g_0(t) = s_0^{-m/2}\bar{g}_0(s_0^{-1}t)$, $t \in \mathbb{R}^{+}$. Let $S(\bs{\Sigma}_0) = s_0$ be a scale function satisfying A1, A2 and A3 and let $\mb{V}_{S,0} \in \mathcal{S}_{m,S}^\mathbb{R}$ the relevant shape matrix. Then, the following chain of equalities hold:
	\begin{equation}\label{CRB_chain_V}
		\begin{split}
			\alpha(\bar{g}_0)^{-1}& 
			\underline{\mb{I}}_m
			\mb{D}_m^{\#}\mb{P}_S(\mb{V}_{S,0})
			(\mb{I}_{m^2}+\mb{K}_{m})(\mb{V}_{S,0}\otimes \mb{V}_{S,0})
			\mb{P}_S^T(\mb{V}_{S,0})\mb{D}_m^{\#T}
			\underline{\mb{I}}_m^T,\\
			&= \mathrm{SCRB}(\ovecs{\mb{V}_{S,0}}|g_0),\quad \quadre{s_0 \; \text{unknown},\; g_0\;\text{functionally unknown}} \\
			&= \mathrm{SCRB}(\ovecs{\mb{V}_{S,0}}|s_0,\bar{g}_0),\quad \quadre{s_0 \; \text{unknown},\; \bar{g}_0\;\text{functionally unknown}} \\
			&= \mathrm{CRB}(\ovecs{\mb{V}_{S,0}}|s_0,\bs{\zeta}_0),\quad \quadre{s_0 \; \text{unknown},\; \bar{g}_0\;\text{functionally known up to its parameters}\; \bs{\zeta}_0} \\
			&= \mathrm{CRB}(\ovecs{\mb{V}_{S,0}}|s_0),\quad \quadre{s_0 \; \text{unknown},\; \bar{g}_0\;\text{fully known}} \\
			&= \mathrm{CRB}(\ovecs{\mb{V}_{S_d,0}}),\; \textit{iff}\; S_d(\bs{\Sigma}_0) \triangleq |\bs{\Sigma}_0|^{1/m}\; \textit{is adopted} \quad \quadre{s_0 \; \text{known},\; \bar{g}_0\;\text{fully known}},
		\end{split}
	\end{equation}
where $\mb{P}_S(\mb{V}_{S,0}) \triangleq \mb{I}_{m^2}-{\rm vec}(\mb{V}_{S,0}) \bs{\nabla}_{{\rm vec}({\bs{\Sigma})}}^T S(\bs{\Sigma}_0)$. A similar result holds true for the location vector $\bs{\mu}_0$.
\end{proposition}
\begin{IEEEproof}
	The proof of the Proposition \ref{Prop_chains} follows directly from the Propositions \ref{Prop_sg} and Corollary \ref{Prop_szeta} and in particular from the fact that:
	\begin{enumerate}
		\item $\bar{\mb{I}}(\bs{\gamma}_0|g_0) = \bar{\mb{I}}(\bs{\gamma}_0|s_0,\bar{g}_0)  = \bar{\mb{I}}(\bs{\gamma}_0|s_0, \bs{\zeta}_0)= \bar{\mb{I}}(\bs{\gamma}_0|s_0) $,
		\item If $S_d(\bs{\Sigma}_0) \triangleq |\bs{\Sigma}_0|^{1/m}$ is adopted, then $\bar{\mb{I}}_{\bs{\gamma}_0}=\bar{\mb{I}}(\bs{\gamma}_0|s_0,g_0)$.
	\end{enumerate}
	The proof is concluded by noticing that the SCRB and the different CRBs are defined as the inverse of the related FIM according to the H\'ajek-Le Cam convolution theorem for parametric models \cite[Sect. 2.3, Th. 1]{BKRW}, semiparametric models \cite[Sect. 3.3, Theo. 2]{BKRW}.   
\end{IEEEproof}

\subsection{Efficiency and robustness trade-off: the $R$-estimator} \label{R_est}

A natural question at this point concerns the attainability of the previously derived semiparametric bound on the shape matrix $\mb{V}_{S,0}$. Indeed, to address this question, in the seminal work \cite{Hallin_Annals_Stat_2}, a \textit{rank-based} $R$-estimator for the shape matrix has been proposed and its asymptotic performance derived. Specifically, in \cite{Hallin_Annals_Stat_2}, the authors showed that an \textit{oracle}, $g_0$-aware, version of the proposed $R$-estimator can attain the semiparametric bound, while, more importantly its \textit{robust}, $g_0$-unaware, version guarantees an excellent trade-off between robustness and (semiparametric) efficiency. In their work, the authors consider the scale function $S(\bs{\Sigma}) = [\bs{\Sigma}]_{11}$. In this subsection, we will show how to construct this estimator for \textit{any scale function} $S$. This generalization follows directly from the outcomes of the recent work \cite{Const_Sem_CRB} which provides a general framework on constrained information bounds in semiparametric models. Due to space limitations, we will not discuss here the properties or theoretical basis of this estimator. Interested readers are referred to \cite{Hallin_P_Annals, Hallin_Annals_Stat_2,Hallin_P_2006,Const_Sem_CRB} while a tutorial introduction more in line with the approach adopted in this paper can be found in \cite{Sem_eff_est_TSP,For_MLSP_SI}. According to the findings presented in the previous section, without loss of generality we assume that $\bs{\mu} = {\bf 0}$.

Let $\mb{x}_l  \sim RES_{m}(\mb{0},\mb{V}_{S,0},g_0)$, $l=1,\ldots,n$ be a set of \textit{independent, identically distributed} (iid), RES distributed, observations. The semiparametric $R$-estimator of $\mb{V}_{S,0}$ is given by:
\begin{equation}
	\label{R_estS}
	\vecs{\widehat{\mb{V}}_{S,R}}  = \vecs{\widehat{\mb{V}}_S^\star} + \frac{1}{\hat{\alpha}\sqrt{n}}\bs{\Xi}_{\widehat{\mb{V}}_S^\star}\bs{\Delta}_{\widehat{\mb{V}}_S^\star},
\end{equation}
where $\widehat{\mb{V}}_S^\star$ is a preliminary, $\sqrt{n}$-consistent, estimator of $\mb{V}_{S,0}$. The matrix $\bs{\Xi}_{\widehat{\mb{V}}_S^\star}$ and the vector $\bs{\Delta}_{\widehat{\mb{V}}_S^\star}$ are defined as:
\begin{equation}
	\bs{\Xi}_{\widehat{\mb{V}}_S^\star} \triangleq 2\mb{U}_{\widehat{\mb{V}}_S^\star}\quadre{\mb{U}_{\widehat{\mb{V}}_S^\star}^T\bs{\Upsilon}_{\widehat{\mb{V}}_S^\star} \bs{\Upsilon}_{\widehat{\mb{V}}_S^\star}^T\mb{U}_{\widehat{\mb{V}}_S^\star}}^{-1}\mb{U}_{\widehat{\mb{V}}_S^\star}^T,
\end{equation}
\begin{equation}
	\bs{\Delta}_{\widehat{\mb{V}}_S^\star} \triangleq \frac{1}{2\sqrt{n}}\bs{\Upsilon}_{\widehat{\mb{V}}_S^\star}\sum_{l=1}^{n}K\tonde{\frac{r^\star_l}{n+1}} \mathrm{vec}(\hat{\mb{u}}^\star_l(\hat{\mb{u}}^\star_l)^T),
\end{equation}
and the scalar $\hat{\alpha}$ can be obtained as \cite[Sect. 4]{Hallin_Annals_Stat_2}.

The function $K:(0,1)\rightarrow \mathbb{R}^+$ is generally called the \textit{score function} and belongs to the set $\mathcal{K}$ of continuous, square integrable functions that can be expressed as the difference of two monotone increasing functions \cite{Hallin_Annals_Stat_2,Hallin_PCA}. 

To define the $R$-estimator in \eqref{R_estS}, we also need the following terms:
\begin{itemize}
	\item $\hat{Q}^\star_l \triangleq \mb{x}_l^T[\widehat{\mb{V}}^\star_S]^{-1}\mb{x}_l$, 
	\item $\hat{\mb{u}}^\star_l \triangleq (\hat{Q}^\star_l)^{-1/2}[\widehat{\mb{V}}^\star_S]^{-1/2}\mb{x}_l$,
	\item $\{r^\star_l\}_{l=1}^n$ are the \textit{ranks} of the (continuous) real random variables $\{\hat{Q}^\star_l\}_{l=1}^n$,
	\item $\bs{\Upsilon}_{\widehat{\mb{V}}^\star_S} \triangleq \mb{D}_m^T\tonde{[\widehat{\mb{V}}^\star_S]^{-1/2}\otimes[\widehat{\mb{V}}^\star_S]^{-1/2}}\Pi^{\perp}_{\cvec{\mb{I}_m}}$, where $\Pi^{\perp}_{\cvec{\mb{I}_m}}=\mb{I}_{m^2} - m^{-1}\mathrm{vec}(\mb{I}_m)\mathrm{vec}(\mb{I}_m)^T$,
	\item $\mb{U}_{\widehat{\mb{V}}_S^\star}$ is the matrix defined in Appendix \ref{App6} in \eqref{U_mat}-\eqref{K_S_U}.
\end{itemize}

To show the effectiveness, in terms of efficiency and robustness, of the $R$-estimator in \eqref{R_estS}, we conducted a simulation study for finite $n$. Specifically, we generate $n$ iid observation vectors drawn form a $t$-distribution with $\nu \in (2,+\infty)$ degrees of freedom \cite[Sect. 1.5.2]{Chap_back_JP} \footnote{We note that for the $t$-distribution, to have a finite second-order moment $\e{\mathcal{Q}^2}<\infty$, we would need $\nu > 4$. However, given that this simulation is carried out within the model \eqref{RES_full_par_model_Sigma} and that the $R$-estimator does not require the second-order moment of $\mathcal{Q}$ to be finite, we can consider values of $\nu$ smaller than 4.} and we evaluated the following Mean Squared Error (MSE) index:
\begin{equation}
	\varepsilon \triangleq E\{\norm{\ovecs{\widehat{\mb{V}}_S-\mb{V}_{S,0}}}_2^2\}
\end{equation}
for the following three estimators:
\begin{enumerate}
	\item the constrained Sample Covariance Matrix $\widehat{\mb{V}}_{S,SCM}$,
	\item the constrained Tyler estimator $\widehat{\mb{V}}_{S,Tyler}$ \cite{Tyler1},
	\item the $R$-estimator in \eqref{R_estS} where the Tyler estimator is used as preliminary estimator, i.e. $\widehat{\mb{V}}_S^\star = \widehat{\mb{V}}_{S,Tyler}$. As score function $K$ we used:
	\begin{itemize}
		\item the \textit{van der Waerden} score $K_{vdW}(u) \triangleq \Psi^{-1}(u)$ where $\Psi^{-1}$ indicates the inverse function of the cdf of a central chi-squared-distributed random variable with $m$ degrees of freedom \cite{Hallin_P_Annals}. As discussed in \cite{Hallin_Annals_Stat_2, Sem_eff_est_TSP}, $K_{vdW}$ has been proved to provide a good compromise between robustness and efficiency. Moreover, the resulting $R$-estimator is asymptotically equivalent to the (constrained) maximum likelihood estimator of the shape matrix for Gaussian data.
		\item the $t_\nu$-based score $K_{t_\nu}(u) = \frac{m(m +\nu)F^{-1}_{m,\nu}(u)}{(\nu + mF^{-1}_{m,\nu}(u))}$ where $F_{m,\nu}(u)$ stands for the cdf of a Fisher random variable with $m$ and $\nu \in (0,\infty)$ degrees of freedom \cite{Hallin_P_Annals}. Here the parameter $\nu$ should be thought of as a knob to adjust the trade-off between robustness and efficiency (at Gaussian data). In fact, small $\nu$ values favor robustness while large values promote efficiency (at Gaussian data), and in particular $\lim_{\nu \rightarrow \infty} K_{t_\nu}(u) = K_{vdW}(u)$. Also in this case, the resulting $R$-estimator is asymptotically equivalent to the (constrained) maximum likelihood estimator of the shape matrix for $t_\nu$-distributed data.
	\end{itemize}
\end{enumerate}
The MSE of the above mentioned estimators is compared to the trace of the semiparametric CRB for $\ovecs{\mb{V}_{S,0}}$ obtained in \eqref{CRB_chain_V} and of the parametric CRB (numerically derived from the inversion of $\mb{I}_{\vsVso}$ in \eqref{I_V}) for the estimation of $\ovecs{\mb{V}_{S,0}}$ in the parametric sub-model $\mathcal{P}_{\bs{\gamma}}$ in \eqref{P1_eta}, i.e. when the scale $s$ and the density generator are perfectly known.

\begin{figure}[h]
	\centering
	\includegraphics[height=6cm]{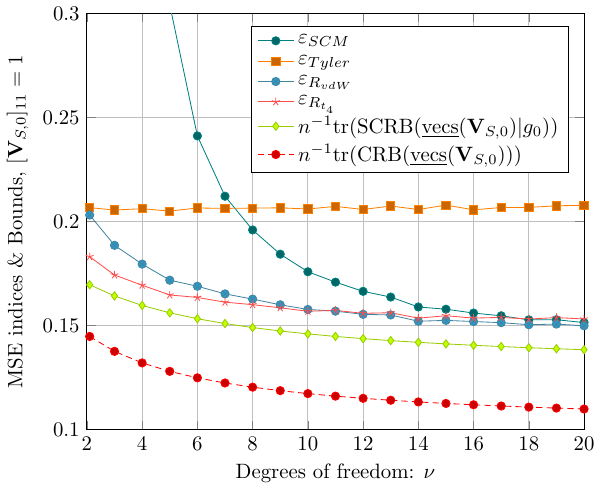}
	\caption{MSE indices and bounds for $S(\bs{\Sigma}) = [\bs{\Sigma}]_{11}$ vs $\nu$ ($n = 100$).}
	\label{fig:Fig1}
\end{figure}

\begin{figure}[h]
	\centering
	\includegraphics[height=6cm]{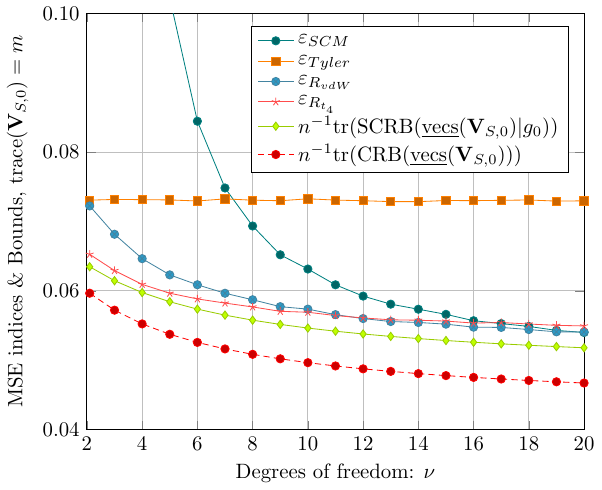}
	\caption{MSE indices and bounds for $S(\bs{\Sigma}) = \trace{\bs{\Sigma}}/m$ vs $\nu$ ($n = 100$).}
	\label{fig:Fig2}
\end{figure}

\begin{figure}[h]
	\centering
	\includegraphics[height=6cm]{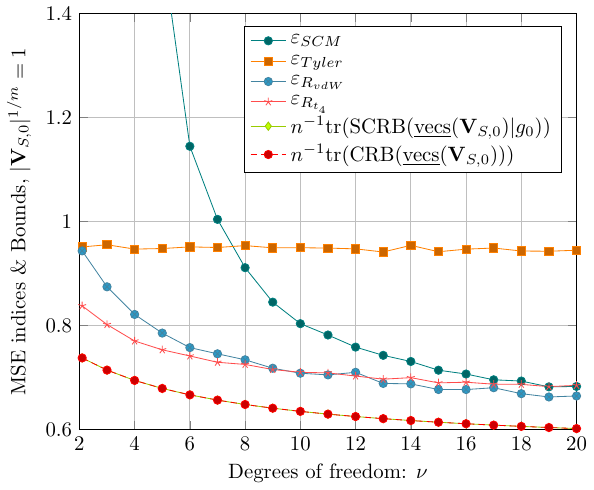}
	\caption{MSE indices and bounds for $S(\bs{\Sigma}) = |\bs{\Sigma}|^{1/m}$ vs $\nu$ ($n = 100$).}
	\label{fig:Fig3}
\end{figure}

In our simulations we generated $n = 100$ $t$-distributed observation vectors of dimension $m = 4$ with a degrees of freedom $\nu$ varying from $2.1$ to $20$. For the $t$-distribution, the coefficient in \eqref{def alpha} is equal to $\alpha(\bar{g}_0) = \frac{m + \nu}{2 + m + \nu}$. The scatter matrix used to generate the data has the Toeplitz form $[\bs{\Sigma}_0]_{ij} = \rho_{0}^{|i-j|}$ with $\rho_{0} = 0.8$. The MSE indices of the estimators have been evaluated by means of $10^5$ Monte-Carlo trials. A preliminary version of the code is available in \cite{R_Github}. \footnote{In \cite{R_Github}, readers can find an implementation of the algorithm for estimating the shape matrix for Circular CES data \cite{Esa}.} 
The simulation analysis is conducted for the three scale functions: $S(\bs{\Sigma}) = [\bs{\Sigma}]_{11}$  (Fig. \ref{fig:Fig1}), $S(\bs{\Sigma}) = \trace{\bs{\Sigma}}/m$ (Fig. \ref{fig:Fig2}) and $S(\bs{\Sigma}) = |\bs{\Sigma}|^{1/m}$ (Fig. \ref{fig:Fig3}). 

The simulations show that the $R$-estimator $\widehat{\mb{V}}_{S,R}$ presents better performance in terms of efficiency and robustness, compared to $\widehat{\mb{V}}_{S,SCM}$ (neither robust nor efficient) and $\widehat{\mb{V}}_{S,Tyler}$ (robust but not efficient). We note, however, that as $\nu \rightarrow \infty$, i.e.\, as the distribution tends to Gaussian, the $\widehat{\mb{V}}_{S,SCM}$ attains the lowest MSE, it is then the optimal estimator. As expected, one can note that the $K_{t_4}$ score function provides better performance, with respect to the van der Waerden score $K_{vdW}$ for heavy-tailed data ($\nu \approx 2$). On the other hand, $K_{vdW}$ is preferable when the data are closed to Gaussian ($\nu \gg 2$). Moreover, Fig. \ref{fig:Fig3} confirms that, as expected from Proposition \ref{Prop_chains}, when the determinant-based scale is adopted, the semiparametric CRB is equal to the parametric $\mathrm{CRB}(\ovecs{\mb{V}_{S_d,0}})$ in $\mathcal{P}_{\bs{\gamma}}$ in \eqref{P1_eta}.

A comment on the implementation of the $R$-estimator in \eqref{R_estS} is in order. Since, in general, the matrix $\mb{U}_{\widehat{\mb{V}}_S^\star}$ depends on the preliminary estimator, it may be of interest to iterate the estimation process. Specifically, as preliminary estimator at step $k$, we may use the $\sqrt{n}$-consistent $R$-estimate obtained at the step $k-1$. For example, we may implement the following iterations:
	\begin{equation}
		\left\lbrace 
		\begin{array}{l}
			\widehat{\mb{V}}_{S,R}^{(0)} = \widehat{\mb{V}}_{S,Tyler},\\
			\vecs{\widehat{\mb{V}}_{S,R}^{(k+1)}}  = \vecs{\widehat{\mb{V}}_{S,R}^{(k)}} + \frac{1}{\hat{\alpha}\sqrt{n}}\bs{\Xi}_{\widehat{\mb{V}}_{S,R}^{(k)}}\bs{\Delta}_{\widehat{\mb{V}}_{S,R}^{(k)}},
		\end{array}
		\right. 
	\end{equation}
	for $k = 0,\ldots,K$. While this iterative procedure does not have any impact on the efficiency of the resulting $R$-estimator, it may bring some improvement in the finite-sample regime. On the contrary, this procedure leads to no or negligible improvement when the matrix $\mb{U}_{\widehat{\mb{V}}_S^\star}$ does not depend on the preliminary estimator, as in the case of the scale functions $S(\bs{\Sigma}) = [\bs{\Sigma}]_{11}$ and $S(\bs{\Sigma}) = \trace{\bs{\Sigma}}/m$.

Clearly, this subsection does not resolve all questions about the possibility to derive \textit{semiparametric efficient} estimates of the shape matrix $\mb{V}_{S,0}$. Further in-depth analysis is needed, both from a theoretical and a simulation perspectives. Special attention should be devoted to the following two points:
\begin{enumerate}
	\item To derive $\hat{\alpha}$ in \eqref{R_estS}, we used the $\sqrt{n}$-consistent but not efficient estimator proposed in \cite[eq. (4.1)]{Hallin_Annals_Stat_2}. An efficient estimation method has been proposed in \cite[Sect. 4.2]{Hallin_Annals_Stat_2}, but presents a considerable computational (numerical) difficulty. Future works will then address its efficient numerical implementation.
	\item A limitation of the estimator in \eqref{R_estS} is that it satisfies the constraint $S(\widehat{\mb{V}}_{S,R})=1$ only asymptotically as $n \rightarrow \infty$ \cite[Theo. 2]{Const_Sem_CRB}. A possible way to \virg{force} the constraint \textit{for each $n$} is to derive $\widehat{\mb{V}}_{S,R}$ directly on the differential manifold induced by $S$, i.e. $\mathcal{S}_{m,S}^\mathbb{R}$ in \eqref{man_S}. This important issue may be addressed building upon the recent and seminal work \cite{asy_manifolds}. Again, this original line of research is left for future works.
\end{enumerate}

\section{Parameterization of the location vector and scatter matrix}\label{sec_para}
In this section, we focus our attention to the case where both the location vector and the scatter matrix can be parametrized by a \textit{real} $d$-dimensional parameter vector $\bs{\theta} = (\bs{\gamma}^T, \bs{\xi}^T)^T \in \Theta  \triangleq \Gamma \times \Psi \subset \mathbb{R}^d$, where $d = q+r$.

\subsection{Some preliminaries}\label{sect_par_pre}
Let $\mathcal{X} \ni \mb{x} \sim RES_m(\bs{\mu}_0,\bs{\Sigma}_0, \bar{g}_0)$ be a RES-distributed random vector whose location vector 
$\bs{\mu}_0 \triangleq \bs{\mu}(\bs{\theta}_0) \in \mathbb{R}^m$ and scatter matrix (which, under the constraint \eqref{set_G_c}, equates the covariance of $\mb{x}$ as shown in \eqref{scatter_cov}) $\bs{\Sigma}_0 \triangleq \bs{\Sigma}(\bs{\theta}_0) \in  \mathcal{S}_m^\mathbb{R}$ are parameterized by a $d$-dimensional parameter vector $\bs{\theta}_0 = (\bs{\gamma}_0^T, \bs{\xi}_0^T)^T \in \Theta \triangleq \Gamma \times \Psi$. As in Sect. \ref{sec_lem1_2}, $\Gamma \subseteq \mathbb{R}^q$ denotes the set of the (finite-dimensional) parameter vectors $\bs{\gamma}$ of interests, $\Psi \subseteq \mathbb{R}^r$ denotes the set of (finite-dimensional) nuisance parameter vectors $\bs{\xi}$ and $g_0 \in \mathcal{G}$ is the (infinite-dimensional) nuisance functions. By introducing the parameterization $\bs{\theta} \mapsto (\bs{\mu}(\bs{\theta}),\bs{\Sigma}(\bs{\theta}))$, we note that, from \eqref{Q_RES}, we have:
\begin{equation}
	\label{Q_alpha}
	(\mb{x}-\bs{\mu}(\bs{\theta}))^T\bs{\Sigma}(\bs{\theta})^{-1}(\mb{x}-\bs{\mu}(\bs{\theta})) =Q_{\bs{\mu}(\bs{\theta}),\bs{\Sigma}(\bs{\theta})}({\bf x}) = \mathcal{Q}\; \forall \bs{\theta} \in \Theta.
\end{equation}

In the rest of this Section, we always assume that the parameterization $\bs{\theta} \mapsto (\bs{\mu}(\bs{\theta}),\bs{\Sigma}(\bs{\theta}))$ following assumptions:
 \begin{itemize}
 	\item[P1)] it is continuous on $\Theta$,
 	\item[P2)] it is locally one-to-one (then \textit{locally identifiable}) in an open neighborhood of $\bs{\theta} \in \Theta$. To ensure this, the Jacobian matrix of the stacked map $\bs{\theta} \mapsto \mb{c}(\bs{\theta})\triangleq(\bs{\mu}(\bs{\theta})^T,\vecs{\bs{\Sigma}(\bs{\theta})}^T)^T$, i.e.\ $\mb{J}[\mb{c}](\bs{\theta})$, has full column rank.
 	\item[P3)] the inverse $\bs{\Sigma}(\bs{\theta})^{-1}$ exists for all $\bs{\theta} \in \Theta$.
 	\item[P4)] In addition, we assume that the second-order moment $\e{\mathcal{Q}^2}$ of $\mathcal{Q}$ in \eqref{Q_alpha}, $\alpha(\bar{g}_0)$ in \eqref{def alpha} and $\beta(\bar{g}_0)$ in \eqref{def beta} are all finite.
 \end{itemize}

As in Sect. \ref{RES_models_essentials}, we can then build the semiparametric model:
\begin{equation}\label{RES_full_par_model_theta_con}
	\mathcal{P}_{\bs{\theta},\bar{g}} = \graffe{p_X(\mb{x}|\bs{\theta},\bar{g}) = |\bs{\Sigma}(\bs{\theta})|^{-1/2}  \bar{g} \left(Q_{\bs{\mu}(\bs{\theta}),\bs{\Sigma}(\bs{\theta})}({\bf x}) \right); \bs{\theta} \in \Theta, \bar{g} \in \overline{\mathcal{G}} }.
\end{equation}
We note, in passing that, if we choose as parameter vector $\bs{\theta} = (\bs{\gamma}^T,s)^T$ in \eqref{eta_Sigma}, we immediately have that $\mathcal{P}_{\bs{\theta},\bar{g}}$ correspond to $\mathcal{P}_{\bs{\eta},\bar{g}}$ in \eqref{RES_full_par_model_V}.

Moreover, exactly as for the semiparametric models presented Sect. \ref{RES_models_essentials}, for the model $\mathcal{P}_{\bs{\theta},\bar{g}}$, we may identify the parametric sub-model:
\begin{equation}\label{P12_eta_t}
	\mathcal{P}_{\bs{\gamma},\bs{\xi}} = \graffe{p_X(\mb{x}|\bs{\gamma},\bs{\xi},\bar{g}_0): \bs{\gamma} \in \Gamma,  \bs{\xi} \in \Psi},
\end{equation}
and the non-parametric sub-model (that is the same of the one introduced in \eqref{con_non_par}) as:
\begin{equation}\label{con_non_par_t} 
	\mathcal{P}^c_{\bar{g}} = \graffe{p_X(\mb{x}|\bs{\gamma}_0,\bs{\xi}_0,\bar{g}): \bar{g} \in \overline{\mathcal{G}}}.
\end{equation}
Note that the (finite-dimensional) nuisance tangent space for the parametric sub-model $\mathcal{P}_{\bs{\gamma},\bs{\xi}}$ can be expressed as
\begin{equation}\label{T_xi}
	\mathcal{T}_{\bs{\xi}_0} = \mathrm{Span}\{[\mb{s}_{\bs{\xi}_0}]_1,\ldots,[\mb{s}_{\bs{\xi}_0}]_r\},
\end{equation}
where $\mb{s}_{\bs{\xi}_0}$ is the score vector of the nuisance parameters $\bs{\xi}_0 \in \Psi$. Moreover, the (infinite-dimensional) nuisance tangent space $\mathcal{T}^c_{\bar{g}_0}$ for the non-parametric model $\mathcal{P}^c_{\bar{g}}$, along with its relative projection operator $\Pi(h|\mathcal{T}^c_{\bar{g}_0})$, are the ones given in \eqref{Tc_g0_1} and \eqref{proj_cg}, respectively. 

To compute this projection operator, let us start by evaluating the score vectors in the parametric sub-model $\mathcal{P}_{\bs{\gamma},\bs{\xi}}$ in \eqref{P12_eta_t}, i.e. $\mb{s}_{\bs{\theta}_0} = (\mb{s}_{\bs{\gamma}_0}^T,\; \mb{s}_{\bs{\xi}_0}^T)^T$ and the related FIM $\mb{I}_{\bs{\theta}_0} = \ev{\mb{s}_{\bs{\theta}_0}\mb{s}_{\bs{\theta}_0}^T}$. 
Following the derivation in \cite[Sec. 3.1]{miss_sb} and \cite[Sec. III]{Bess}, each entry of $\mb{s}_{\bs{\theta}_0}$ can be expressed as:
\begin{equation}
	\label{score_vect_rep}
	\begin{split}
		[\mb{s}_{\bs{\theta}_0}]_i = - 2^{-1}\mathrm{tr}( \mb{P}_i^0) + \bar{\varphi}_0(\mathcal{Q}) \left( \sqrt{\mathcal{Q}} \mb{u}^T \bs{\Sigma}_0^{-1/2}\bs{\mu}_i^0  + 2^{-1}\mathcal{Q}\mb{u}^T \mb{P}_i^0 \mb{u} \right), \; i=1,\ldots,d. 
	\end{split}
\end{equation}
where $\mb{P}_i^0 \triangleq \bs{\Sigma}_0^{-1/2}\bs{\Sigma}_i^0\bs{\Sigma}_0^{-1/2}$, $\bs{\Sigma}_i^{0}\triangleq {\frac{\partial \bs{\Sigma}(\bs{\theta})}{\partial \theta_i}}|_{\bs{\theta}=\bs{\theta}_0}$ and $\bs{\mu}_i^0  \triangleq \frac{\partial \bs{\mu}(\bs{\theta})}{\partial \theta_i}|_{\bs{\theta}=\bs{\theta}_0}.$
Consequently, $[\mb{s}_{\bs{\gamma}_0}]_i = [\mb{s}_{\bs{\theta}_0}]_i$, $i=1,\ldots,q$ and $[\mb{s}_{\bs{\xi}_0}]_j = [\mb{s}_{\bs{\theta}_0}]_{q+j}$, $j=1,\ldots,r$.

To derive the FIM $\mb{I}_{\bs{\theta}_0}$, we may use the procedure in \cite{Bess} that leads to the following compact expression given in \cite[Sect. 1.6.5]{Chap_back_JP} as:
\begin{equation}\label{I_theta_par}
	\begin{split}
		\mb{I}_{\bs{\theta}_0}  = \beta(\bar{g}_0)&\mb{J}[\bs{\mu}_0]^T\bs{\Sigma}_0^{-1}\mb{J}[\bs{\mu}_0] \\
		&+2^{-1}\alpha(\bar{g}_0)\mb{J}[\cvec{\bs{\Sigma}_0}]^T\quadre{\bs{\Sigma}_0^{-1} \otimes \bs{\Sigma}_0^{-1} + 2^{-1}(1-\alpha^{-1}(\bar{g}_0)) \mathrm{vec}(\bs{\Sigma}_0^{-1})\mathrm{vec}(\bs{\Sigma}_0^{-1})^T}\mb{J}[\cvec{\bs{\Sigma}_0}],
	\end{split}
\end{equation}
where the scalars $\beta(\bar{g}_0)$ and $\alpha(\bar{g}_0)$ are given in \eqref{def beta} and \eqref{def alpha} respectively, while the Jacobian matrices, $\mb{J}[\bs{\mu}_0]= [\mb{J}_{\bs{\gamma}}[\bs{\mu}_0], \mb{J}_{\bs{\xi}}[\bs{\mu}_0]]$ and $\mb{J}[\cvec{\bs{\Sigma}_0}] = [\mb{J}_{\bs{\gamma}}[\cvec{\bs{\Sigma}_0}], \mb{J}_{\bs{\xi}}[\cvec{\bs{\Sigma}_0}]]$, are explicitly expressed as
\begin{equation}\label{def J}
	[\mb{J}[\bs{\mu}_0]]_{i,j} = \left. \parder{[\bs{\mu}(\bs{\theta})]_i}{\theta_j}\right|_{\bs{\theta}=\bs{\theta}_0}, \; [\mb{J}[\cvec{\bs{\Sigma}_0}]]_{i,j} = \left. \parder{[\cvec{\bs{\Sigma}(\bs{\theta})}]_i}{\theta_j}\right|_{\bs{\theta}=\bs{\theta}_0},
\end{equation}
where $\mb{J}[\bs{\mu}_0] \in \mathbb{R}^{m \times d}$, $\mb{J}_{\bs{\gamma}}[\bs{\mu}_0] \in \mathbb{R}^{m \times q}$ and $\mb{J}_{\bs{\xi}}[\bs{\mu}_0] \in \mathbb{R}^{m \times r}$. Moreover, $\mb{J}[\cvec{\bs{\Sigma}_0}] \in \mathbb{R}^{m^2 \times d}$, $\mb{J}_{\bs{\gamma}}[\cvec{\bs{\Sigma}_0}] \in \mathbb{R}^{m^2 \times q}$ and $\mb{J}_{\bs{\xi}}[\cvec{\bs{\Sigma}_0}] \in \mathbb{R}^{m^2 \times r}$. It is important to note that the Assumption P2 guarantees that $\mb{I}_{\bs{\theta}_0}$ is invertible.

We can now introduce the semiparametric efficient score vector $\bar{\mb{s}}_{\bs{\theta}_0} = \mb{s}_{\bs{\theta}_0} - \Pi(\mb{s}_{\bs{\theta}_0}|\mathcal{T}^c_{\bar{g}_0})$ for $\bs{\theta}_0 \in \Theta$ in the semiparametric model $\mathcal{P}_{\bs{\theta},\bar{g}}$ in \eqref{RES_full_par_model_theta_con}. From the expression of $\mb{s}_{\bs{\theta}_0}$ in \eqref{score_vect_rep} and the one of the projection operator $\Pi(\cdot|\mathcal{T}^c_{\bar{g}_0})$ in \eqref{proj_cg}, we have that:
\begin{equation}\label{eff_score_par}
	\begin{split}
		[\bar{\mb{s}}_{\bs{\theta}_0}]_i &= [\mb{s}_{\bs{\theta}_0}]_i - E\{[\mb{s}_{\bs{\theta}_0}]_i|\mathcal{Q}\} + \eg{\mathcal{Q}[\mb{s}_{\bs{\theta}_0}]_i}\sigma_{\mathcal{Q}}^{-2}(\mathcal{Q}-m)\\
		&=\sqrt{\mathcal{Q}}\bar{\varphi}_0(\mathcal{Q}) \mb{u}^T \bs{\Sigma}_0^{-1/2}\bs{\mu}_i^0 + 2^{-1}\mathcal{Q}\bar{\varphi}_0(\mathcal{Q})\tonde{\mb{u}^T \mb{P}_i^0 \mb{u} -  m^{-1}\mathrm{tr}( \mb{P}_i^0)} +\mathrm{tr}( \mb{P}_i^0)\sigma_{\mathcal{Q}}^{-2}(\mathcal{Q}-m),
	\end{split}
\end{equation}
where we have used the equality $E\{\mathcal{Q}^2\bar{\varphi}_0(\mathcal{Q})\}=m(m+2)$ given in footnote 3.

Moreover, by direct calculation,
the efficient semiparametric FIM for $\bs{\theta}_0$ in $\mathcal{P}_{\bs{\theta},\bar{g}}$ can be expressed as:
\begin{equation}\label{bar_I_theta_par}
	\begin{split}
		\bar{\mb{I}}_{\bs{\theta}_0}  &= \beta(\bar{g}_0)\mb{J}[\bs{\mu}_0]^T\bs{\Sigma}_0^{-1}\mb{J}[\bs{\mu}_0] + \\
		&\quad +2^{-1}\alpha(\bar{g}_0)\mb{J}[\cvec{\bs{\Sigma}_0}]^T \quadre{\bs{\Sigma}_0^{-1} \otimes \bs{\Sigma}_0^{-1}  +( 2\alpha^{-1}(\bar{g}_0)\sigma_{\mathcal{Q}}^{-2}-m^{-1})\cvec{\bs{\Sigma}_0^{-1}}\cvec{\bs{\Sigma}_0^{-1}}^T  } \mb{J}[\cvec{\bs{\Sigma}_0}],
	\end{split}
\end{equation} 
where we recall that $\sigma_{\mathcal{Q}}^2 \triangleq \e{(\mathcal{Q}-m)^2} = \e{\mathcal{Q}^2} - m^2$. \footnote{It is important to note that the expression of the efficient semiparametric score vector $\bar{\mb{s}}_{\bs{\theta}_0}$ and the related FIM $\bar{\mb{I}}_{\bs{\theta}_0}$ derived in eqs. (43) and (46) of \cite{For_SCRB_complex} are incorrect. In fact, in the derivation proposed in \cite{For_SCRB_complex}, the projection was incorrectly taken with respect to the \virg{unconstrained} nuisance tangent space $\mathcal{T}_{g_0}^u$, whereas $\mb{s}_{\bs{\theta}_0}$ must be projected onto the \virg{constrained} one, i.e. $\mathcal{T}_{\bar{g}_0}^c$.} 

\textit{Remark}: One can expect that, knowing the density generator $\bar{g}_0$ should lead to a larger (in the sense of the positive definite matrices ordering) FIM for $\bs{\theta}_0$. However, proving that ${\mb{I}}_{\bs{\theta}_0}\ge\bar{\mb{I}}_{\bs{\theta}_0}$ through the inequality $2^{-1}(1-\alpha^{-1}(\bar{g}_0))\ge 2\alpha^{-1}(\bar{g}_0)\sigma_{\mathcal{Q}}^{-2}-m^{-1}$ from \eqref{I_theta_par}, \eqref{bar_I_theta_par} for arbirtary RES distribution seems quite involved. Nevertheless, it is straightforward to prove that for the Gaussian distribution for which $\mathcal{Q}$ is $\chi_m^2$ distributed, $\alpha(\bar{g}_0)=1$ and $\sigma_{\mathcal{Q}}^{2}=2m$ imply ${\mb{I}}_{\bs{\theta}_0}=\bar{\mb{I}}_{\bs{\theta}_0}$ and for the $t$-distribution, $\alpha(\bar{g}_0)=\frac{m+\nu}{m+\nu+2}$ and $\sigma_{\mathcal{Q}}^{2}=(\frac{2m}{\nu-4})(m+\nu-2)$ imply
${\mb{I}}_{\bs{\theta}_0}>\bar{\mb{I}}_{\bs{\theta}_0}$.

%
%
%

The calculation of the efficient semiparametric score vector of the parameters of interest $\bar{\mb{s}}_{\bs{\gamma}_0} = \mb{s}_{\bs{\gamma}_0} - \Pi(\mb{s}_{\bs{\gamma}_0}|\mathcal{T}_{\bs{\xi}_0} + \mathcal{T}^c_{\bar{g}_0})$ and the related efficient semiparametric FIM $\bar{\mb{I}}(\bs{\gamma}_0|\bs{\xi}_0,\bar{g}_0)$ is more difficult task. In fact, even if the projection $\Pi(h|\mathcal{T}_{\bs{\xi}_0} + \mathcal{T}^c_{\bar{g}_0})$, $\forall h \in \mathcal{H}$ can be explicitly obtained as shown in Appendix \ref{App_Theo_sum} (see also \cite[Theo. 5, A.4]{BKRW}), the calculation may be quite involved depending on the parameterisation at hands. 

It may therefore be of practical interest to ask the following question: \textit{is it possible to characterize all parameterisations $\bs{\theta} \mapsto (\bs{\mu}(\bs{\theta}),\bs{\Sigma}(\bs{\theta}))$ that imply that not knowing the finite-dimensional nuisance vector $\bs{\xi}_0 \in \Psi$, while knowing $\bar{g}_0 \in \overline{\mathcal{G}}$, leads to the same loss of efficiency as not knowing $\bs{\xi}_0 \in \Psi$ \textit{nor} the density generator $\bar{g}_0$?}

As we will see below, many parameterisations used in applications satisfy this adaptivity property. 
We will therefore be able to easily compute the relevant $\bar{\mb{s}}_{\bs{\gamma}_0}$ and $\bar{\mb{I}}(\bs{\gamma}_0|\bs{\xi}_0,\bar{g}_0)$.

\subsection{An adaptivity condition for the parameterized case}
By collecting the previous results, we are ready to state the following original proposition:

\begin{proposition}\label{prop_corr_1_RES}
	Let $\mathbb{R}^m \ni \mb{x} \sim RES_m(\bs{\mu}_0,\bs{\Sigma}_0, \bar{g}_0)$ be a RES distributed vector whose location vector $\bs{\mu}_0 \triangleq \bs{\mu}(\bs{\theta}_0)$ and scatter matrix $\bs{\Sigma}_0 \triangleq \bs{\Sigma}(\bs{\theta}_0)$ are parameterized by $\bs{\theta}_0 = (\bs{\gamma}_0^T, \bs{\xi}_0^T)^T \in \Gamma \times \Psi$, where $\bs{\gamma}_0 \in \Gamma$ is the $q$-dimensional vector of interest and $\bs{\xi}_0 \in \Psi$ is the $r$-dimensional nuisance vector, such that $d = q+r$. Let $\bar{g}_0 \in \overline{\mathcal{G}}$ be the infinite-dimensional nuisance parameter.
	Then the following
	\begin{equation}\label{cond_corr_2_RES}
		\mathrm{Condition\; \ref{prop_corr_1_RES}:}\tonde{ \mb{J}_{\bs{\gamma}}^T[\cvec{\bs{\Sigma}_0}] 
			- \mb{I}_{\bs{\gamma}_0\bs{\xi}_0}\mb{I}_{\bs{\xi}_0}^{-1}\mb{J}_{\bs{\xi}}^T[\cvec{\bs{\Sigma}_0}]}
		\cvec{\bs{\Sigma}_0^{-1}} 
		= \mb{0},
	\end{equation}
	where $\mb{I}_{\bs{\gamma}_0\bs{\xi}_0}$ and $\mb{I}_{\bs{\xi}_0}$ are two sub-blocks of the FIM $\mb{I}_{\bs{\theta}_0} = \tonde{\begin{array}{cc}
			\mb{I}_{\bs{\gamma}_0} & \mb{I}_{\bs{\gamma}_0\bs{\xi}_0}\\
			\mb{I}_{\bs{\gamma}_0\bs{\xi}_0}^T & \mb{I}_{\bs{\xi}_0}
	\end{array}}$ given in \eqref{I_theta_par},
	\begin{itemize}
		\item is \textbf{sufficient} for the equality between the efficient SFIM $\bar{\mb{I}}(\bs{\gamma}_0|\bs{\xi}_0,\bar{g}_0)$ for the model $\mathcal{P}_{\bs{\theta},\bar{g}}$ in \eqref{RES_full_par_model_theta_con} and the parametric efficient FIM $\bar{\mb{I}}(\bs{\gamma}_0|\bs{\xi}_0)$ in $\mathcal{P}_{\bs{\gamma},\bs{\xi}}$ in \eqref{P12_eta_t}:
		\begin{equation}\label{cond_corr_2_RES_suf}
			\mathrm{Condition\; \ref{prop_corr_1_RES}} \; \Rightarrow \;\bar{\mb{I}}(\bs{\gamma}_0|\bs{\xi}_0,\bar{g}_0) = \bar{\mb{I}}(\bs{\gamma}_0|\bs{\xi}_0)\; \forall \bar{g}_0 \in \overline{\mathcal{G}}.
		\end{equation}
		\item is also \textbf{necessary} for all $\bar{g}_0 \in \overline{\mathcal{G}}$ except the Gaussian one, i.e. $\bar{g}_G(t) =(2\pi)^{-m/2}\exp(-t/2)$:
		\begin{equation}\label{cond_corr_2_RES_nec}
			\mathrm{Condition\; \ref{prop_corr_1_RES}} \; \Leftarrow \;\bar{\mb{I}}(\bs{\gamma}_0|\bs{\xi}_0,\bar{g}_0) = \bar{\mb{I}}(\bs{\gamma}_0|\bs{\xi}_0)\; \forall \bar{g}_0 \neq \bar{g}_G \in \overline{\mathcal{G}}.
		\end{equation}
		For the Gaussian case, we have that $\bar{\mb{I}}(\bs{\gamma}_0|\bs{\xi}_0,\bar{g}_G) = \bar{\mb{I}}(\bs{\gamma}_0|\bs{\xi}_0)$ irrespective of the parameterization.
	\end{itemize}

\end{proposition}	
\begin{IEEEproof}
	Let us apply Lemma \ref{Lem2} to the case of RES distributions, we can identify $\mathcal{T}_{2}$ with  $\mathcal{T}_{\bs{\xi}_0}$ in \eqref{T_xi} and $\mathcal{T}_{3}$ with $\mathcal{T}^c_{\bar{g}_0}$ in \eqref{Tc_g0_1}. 
	Then, Condition \eqref{main_cond} of Lemma  \ref{Lem2} can be expressed as:
		\begin{equation}\label{main_cond_RES1}
			[\mb{s}_{\bs{\gamma}_0} - \Pi(\mb{s}_{\bs{\gamma}_0}|\mathcal{T}_{\bs{\xi}_0})] \perp \mathcal{T}^c_{\bar{g}_0}.
		\end{equation}
		Since, as shown in Appendix \ref{proof_tan_spa}, $\mathcal{T}^u_{g_0} = \mathcal{T}^c_{\bar{g}_0} \oplus \mathrm{Span}\{\mathcal{Q}-m\}$, a necessary and sufficient condition for \eqref{main_cond_RES1} to hold is that: 
		\begin{equation}\label{main_cond_RES}
			[\Pi\tonde{\mb{s}_{\bs{\gamma}_0} - \Pi(\mb{s}_{\bs{\gamma}_0}|\mathcal{T}_{\bs{\xi}_0})|\mathcal{T}^u_{g_0}}]_i \in \mathrm{Span}\{\mathcal{Q}-m\},
		\end{equation}
		for $i=1,\ldots,q$, that is, from \eqref{proj_ug}, $[E\graffe{\mb{s}_{\bs{\gamma}_0} - \Pi(\mb{s}_{\bs{\gamma}_0}|\mathcal{T}_{\bs{\xi}_0})|\mathcal{Q}}]_i \in \mathrm{Span}\{\mathcal{Q}-m\}$

	
		By using the expression of $\Pi(\mb{s}_{\bs{\gamma}_0}|\mathcal{T}_{\bs{\xi}_0}) $ deduced from \eqref{proj_T2}, we have:
		\begin{equation}\label{cond_corr2_appo}
			E\graffe{\Pi(\mb{s}_{\bs{\gamma}_0}|\mathcal{T}_{\bs{\xi}_0})|\mathcal{Q}} = \mb{I}_{\bs{\gamma}_0\bs{\xi}_0}\mb{I}_{\bs{\xi}_0}^{-1}E\graffe{\mb{s}_{\bs{\xi}_0}|\mathcal{Q}}.
		\end{equation}
		By definition, $E\graffe{\mb{s}_{\bs{\gamma}_0}|\mathcal{Q}}$ and $E\graffe{\mb{s}_{\bs{\xi}_0}|\mathcal{Q}}$ are the two sub-vector of $E\{\mb{s}_{\bs{\theta}_0} |\mathcal{Q}\}$, then, to conclude the proof, we just need to evaluate $E\{\mb{s}_{\bs{\theta}_0} |\mathcal{Q}\}$. To this end, from the expression of $\mb{s}_{\bs{\theta}_0}$ given in \eqref{score_vect_rep} and by noticing that:
		\begin{equation}\label{trace_vec}
			\mathrm{tr}\left( \mb{P}_i^0 \right)  = \mathrm{tr}\left( \bs{\Sigma}_i^0 \bs{\Sigma}_0^{-1} \right)=  \cvec{ \bs{\Sigma}_0^{-1}}^T \cvec{\bs{\Sigma}_i^0},
		\end{equation}
		we immediately have that:
		\begin{equation}\label{stheta_T3}
			\begin{split}
				E\{\mb{s}_{\bs{\theta}_0} |\mathcal{Q}\} &=  2^{-1}(m^{-1} \mathcal{Q} \bar{\varphi}_0(\mathcal{Q}) - 1)\mb{J}[\cvec{\bs{\Sigma}_0}]^T\cvec{\bs{\Sigma}_0^{-1}}\\
				&= 2^{-1}(m^{-1} \mathcal{Q} \bar{\varphi}_0(\mathcal{Q}) - 1)[\mb{J}_{\bs{\gamma}}[\cvec{\bs{\Sigma}_0}], \mb{J}_{\bs{\xi}}[\cvec{\bs{\Sigma}_0}]]^T\cvec{\bs{\Sigma}_0^{-1}}. 
			\end{split}
		\end{equation}
		Consequently, from \eqref{stheta_T3} and \eqref{cond_corr2_appo}, the condition \eqref{main_cond_RES} holds true if and only if
		\begin{equation}\label{final_cond}
			a_i 2^{-1}(m^{-1} \mathcal{Q} \bar{\varphi}_0(\mathcal{Q}) - 1) \in \mathrm{Span}\{\mathcal{Q}-m\},
		\end{equation}
		for $i = 1,\ldots,q$ and where $a_i$ are the $q$ components of the vector in the left-hand side of \eqref{cond_corr_2_RES}. 
		
		To conclude the proof, we can trivially verified that \footnote{Note that $\mathrm{Pr}\tonde{\{m^{-1} \mathcal{Q} \bar{\varphi}_0(\mathcal{Q})=1\}}=0$, since $\mathcal{Q}$ is a continuous random variable.} $(m^{-1} \mathcal{Q} \bar{\varphi}_0(\mathcal{Q})-1)\notin \mathrm{Span}\{(\mathcal{Q}-m)\}$ for all $\bar{g}_0 \neq  \bar{g}_G$. Consequently, \eqref{final_cond} is verified \textit{if and only if} $a_i = 0$, for all $i=1,\ldots,q$ that is equivalent to the condition in \eqref{cond_corr_2_RES}. On the other hands, for the Gaussian case, we have that $\bar{\varphi}_0(t) = 1,\forall t \in \mathbb{R}^+$, then $(m^{-1} \mathcal{Q} \bar{\varphi}_0(\mathcal{Q})-1) = m^{-1}(\mathcal{Q}-m) \in \mathrm{Span}\{(\mathcal{Q}-m)$ irrespective of the values of $a_i$. The fact that the Gaussian case is somehow special can be seen directly by comparing the FIM $\mb{I}_{\bs{\theta}_0}$ in \eqref{I_theta_par} and its semiparametric counterpart $\bs{\theta}_0$ in \eqref{bar_I_theta_par}. In fact, as already noticed, for the Gaussian case, we have that $\mb{I}_{\bs{\theta}_0} = \bar{\mb{I}}_{\bs{\theta}_0}$ that trivially implies that $\bar{\mb{I}}(\bs{\gamma}_0|\bs{\xi}_0,\bar{g}_G) = \bar{\mb{I}}(\bs{\gamma}_0|\bs{\xi}_0)$ irrespective of the parameterization.
\end{IEEEproof}

As examples of how condition \eqref{cond_corr_2_RES} can be used, let us consider two quite common parametrizations in signal processing application: 
the case where the location vector and the scatter matrix have no parameters in common and  the \virg{low-rank}, parameterization model.

\subsection{The elliptical parameterized model where the location vector and the scatter matrix have no parameters in common}
Let $\mathbb{R}^m \ni \mb{x} \sim RES_{m}(\bs{\mu}_0, \bs{\Sigma}_0, \bar{g}_0)$ where $\bs{\mu}_0 \triangleq \bs{\mu}(\bs{\gamma}_0)$ is parameterized by the parameter of interest $\bs{\gamma}_0 \in \Gamma \subseteq  \mathbb{R}^{q}$ while the scatter matrix $\bs{\Sigma}_0 \triangleq \bs{\Sigma}(\bs{\xi}_0)$ is parameterized by the nuisance parameter
$\bs{\xi}_0  \in \Psi \subseteq \mathbb{R}^{r}$. Furthermore, we assume that such parameterisation satisfies the assumptions P1 and P2 given in Sect. \ref{sect_par_pre}.

Let us consider the three following models:
\begin{equation}\label{model location scatter semiparametrique}
\mathcal{P}_{\bs{\gamma},\bs{\xi},\bar{g}} \triangleq \graffe{p_X(\mb{x}|\bs{\gamma},\bs{\xi},\bar{g}) 
= |\bs{\Sigma}(\bs{\xi})|^{-1/2} 
 \bar{g}((\mb{x}-\bs{\mu}(\bs{\gamma}))^T\bs{\Sigma}(\bs{\xi})^{-1}(\mb{x}-\bs{\mu}(\bs{\gamma}))) : \bs{\gamma}\in \Gamma,  \bs{\xi}\in \Psi, \bar{g} \in \overline{\mathcal{G}}},
\end{equation}
\begin{equation}\label{model location scatter parametrique}
\mathcal{P}_{\bs{\gamma},\bs{\xi}}\triangleq \graffe{p_X(\mb{x}|\bs{\gamma},\bs{\xi},\bar{g}_0) 
= |\bs{\Sigma}(\bs{\xi})|^{-1/2} 
 \bar{g}_0((\mb{x}-\bs{\mu}(\bs{\gamma}))^T\bs{\Sigma}(\bs{\xi})^{-1}(\mb{x}-\bs{\mu}(\bs{\gamma}))) : \bs{\gamma}\in \Gamma,  \bs{\xi}\in \Psi},
\end{equation}
\begin{equation}\label{model location scatter parametrique decoupled}
\mathcal{P}_{\bs{\gamma}} \triangleq \graffe{p_X(\mb{x}|\bs{\gamma},\bs{\xi}_0,\bar{g}_0)
= |\bs{\Sigma}(\bs{\xi}_0)|^{-1/2} 
 \bar{g}_0((\mb{x}-\bs{\mu}(\bs{\gamma}))^T\bs{\Sigma}(\bs{\xi}_0)^{-1}(\mb{x}-\bs{\mu}(\bs{\gamma}))) : \bs{\gamma}\in \Gamma}.
\end{equation}
Then, the following proposition holds true.
\begin{proposition}\label{Prop_location scatter model} 
	Let $\mathbb{R}^m \ni \mb{x} \sim RES_{m}(\bs{\mu}_0, \bs{\Sigma}_0, \bar{g}_0)$ where $\bs{\mu}_0 \triangleq \bs{\mu}(\bs{\gamma}_0)$ and $\bs{\Sigma}_0 \triangleq \bs{\Sigma}(\bs{\xi}_0)$. Then, the efficient Semiparametric FIM (SFIM) $\bar{\mb{I}}(\bs{\gamma}_0|\bs{\xi}_0,\bar{g}_0)$ for the model \eqref{model location scatter semiparametrique}
	is equal to the parametric efficient FIM $\bar{\mb{I}}(\bs{\gamma}_0|\bs{\xi}_0)$ in \eqref{model location scatter parametrique}. Moreover this latter FIM is equal to the FIM  $\mb{I}_{\bs{\gamma}_0}$ for the model \eqref{model location scatter parametrique decoupled}.
\end{proposition}
\begin{IEEEproof}
For this specific parameterisation, we trivially have that  $\mb{J}_{\bs{\gamma}}[\cvec{\bs{\Sigma}_0}]=\mb{0}$ and $\mb{J}_{\bs{\xi}}[\bs{\mu}_0]=\mb{0}$.
Furthermore, since  $\mb{J}[\bs{\mu}_0]= [\mb{J}_{\bs{\gamma}}[\bs{\mu}_0], \mb{J}_{\bs{\xi}}[\bs{\mu}_0]]$ and $\mb{J}[\cvec{\bs{\Sigma}_0}]= [\mb{J}_{\bs{\gamma}}[\cvec{\bs{\Sigma}_0}], \mb{J}_{\bs{\xi}}[\cvec{\bs{\Sigma}_0}]]$, it follows from \eqref{I_theta_par} that the FIM $\mb{I}_{\bs{\theta}_0}$ is block-diagonal, i.e. $\mb{I}_{\bs{\gamma}_0\bs{\xi}_0}=\mb{0}$. Consequently the condition \eqref{cond_corr_2_RES} of Proposition \ref{prop_corr_1_RES} is satisfied and $\bar{\mb{I}}(\bs{\gamma}_0|\bs{\xi}_0)$ reduces to $\mb{I}_{\bs{\gamma}_0}$.
\end{IEEEproof}

\subsection{The elliptical parameterized \virg{low-rank} model}
Let $\mathbb{R}^m \ni \mb{x} \sim RES_{m}(\mb{0}, \bs{\Sigma}_0, \bar{g}_0)$ be a zero-mean, RES-distributed vector whose scatter matrix $\bs{\Sigma}_0$ is modeled as:
\begin{equation}
	\label{Cov_mat_DOA}
	\bs{\Sigma}_0 \equiv \bs{\Sigma}(\bs{\gamma}_0,\bs{\xi}_0) = \mb{A}_0\bs{\Xi}_0\mb{A}_0^T + \lambda_0\mb{I}_m,
\end{equation}
where:
\begin{itemize}
	\item $\mb{A}(\bs{\gamma}) \in \mathbb{R}^{m \times p}$ is differentiable and its Jacobian $\mb{J}[\cvec{\bs{\mb{A}}}](\bs{\gamma})$ has a full rank column for all $\bs{\gamma} \in \Gamma$ and then $\mb{A}(\bs{\gamma})$ is locally injective in a neighborhood of $\bs{\gamma}_0 \in \Gamma$, 
	\item $\bs{\Xi}_0 \in  \mathcal{S}_p^\mathbb{R}$ is a symmetric and \textit{positive definite} matrix, 
	\item $\lambda_0 \in \mathbb{R}^+$.
\end{itemize}
Consequently, the finite-dimensional nuisance vector can be defined as:
\begin{equation}
	\bs{\xi}_0 \triangleq (\vecs{\bs{\Xi}_0}^T, \lambda_0)^T \in \Psi \subset \vecs{\mathcal{S}_p^\mathbb{R}} \times \mathbb{R}^+ \subseteq \mathbb{R}^{r},
\end{equation}
with $r=p(p+1)/2+1$, while the infinite dimensional nuisance is the density generator $\bar{g}_0 \in \overline{\mathcal{G}}$. 

Let us consider the two following models related to the low-rank parameterization of the scatter matrix in \eqref{Cov_mat_DOA}:
\begin{equation}\label{model low rank semiparametrique}
\mathcal{P}_{\bs{\gamma},\bs{\xi},\bar{g}} \triangleq \graffe{p_X(\mb{x}|\bs{\gamma},\bs{\xi},\bar{g}) = |\bs{\Sigma}(\bs{\gamma},\bs{\xi})|^{-1/2}  \bar{g}(\mb{x}^T\bs{\Sigma}(\bs{\gamma},\bs{\xi})^{-1}\mb{x}) ; \bs{\gamma}\in \Gamma,  \bs{\xi}\in \Psi, \bar{g} \in \overline{\mathcal{G}}},
\end{equation}
\begin{equation}\label{model low rank parametrique}
\mathcal{P}_{\bs{\gamma},\bs{\xi}} \triangleq \graffe{p_X(\mb{x}|\bs{\gamma},\bs{\xi},\bar{g}_0) = 
 |\bs{\Sigma}(\bs{\gamma},\bs{\xi})|^{-1/2}  \bar{g}_0(\mb{x}^T\bs{\Sigma}(\bs{\gamma},\bs{\xi})^{-1}\mb{x}) ; \bs{\gamma}\in \Gamma,  \bs{\xi}\in \Psi}.
\end{equation}
The following result can be proved by a direct application of Condition \eqref{cond_corr_2_RES} of Proposition \ref{prop_corr_1_RES} as shown in Appendix \ref{proof_Prop_low-rank model}.
\begin{proposition}\label{Prop_low-rank model} 
For the low-rank scatter model \eqref{Cov_mat_DOA}, the efficient Semiparametric FIM (SFIM) $\bar{\mb{I}}(\bs{\gamma}_0|\bs{\xi}_0,\bar{g}_0)$ is equal to the parametric efficient FIM $\bar{\mb{I}}(\bs{\gamma}_0|\bs{\xi}_0)$.
\end{proposition}

\section{Applications to circular and noncircular CES distributions}
\label{Applications}
So far in this paper, we have only dealt with cases of real observation vectors. 
However complex-valued observations, i.e. $\mb{x} \in \mathbb{C}^m$ are an integral part of many science and engineering problems, including those
in communications, radar, biomedicine, geophysics, oceanography, electromagnetics, and optics, among
others. The complex field does not only provide a convenient representation for the observations but also
provides a natural way to capture their physical nature as well as the transformations
they go through (see e.g. \cite{Schreier,Mandic}).
In many studies it has often been (implicitly) assumed that complex random vectors are \emph {circular}, i.e. with invariant distribution under rotation around a center $\bs{\mu}$, that is $(\mb{x}-\bs{\mu}) = e^{j\theta}(\mb{x}-\bs{\mu}),\; \forall \theta \in \mathbb{R}$ (see e.g. in \cite{Esa}). This assumption however discards the information conveyed by the relationship between real and imaginary
parts of the observation vectors. Consequently, the \emph{noncircularity} may be an important feature that characterizes observation in many practical scenarios where the data are \textit{non-stationary}.
For this reason, in this section, we believe it is important to provide some evidence to ensure that all results obtained so far in the case of real observations remain entirely valid for Complex Elliptically Symmetric (CES)-distributed observations where the parameters are still real-valued.
Moreover, given their importance in well-known engineering applications, we devote particular attention to the \virg{complex} version of the Propositions \ref{Prop_location scatter model} and \ref{Prop_low-rank model}. Before moving on, it is important to note that, although there are many different notations in the literature, in the following sections we will use the notation introduced in \cite{Chap_back_JP}. 
%
\subsection{Real-complex representations}
A random observation vector $\mb{x} \in \mathcal{X} \subseteq \mathbb{C}^m$ is said to be CES distributed if  the associated real-valued vector
$\overline{\bf x} \in \mathbb{R}^{2m}$ with $\overline{\bf x} \triangleq ({\rm Re}({\bf x}^T),{\rm Im}({\bf x}^T))^T$ is RES distributed, i.e., with pdf given by \eqref{RES_pdf} when it exists. It follows that all the properties of RES distributions and propositions given for even $m$ in the previous sections apply for CES distributions.
However, it is more convenient to express these properties and propositions using notations suited to the complex representation that naturally arises when using the one-to-one mapping 
\begin{equation*}
	\overline{\bf x} \mapsto \widetilde{\bf x}\triangleq ({\bf x}^T,{\bf x}^H)^T= \sqrt{2}{\bf M}\overline{\bf x} 
\end{equation*}
where 
${\bf M}\triangleq
\tiny{\frac{1}{\sqrt{2}}\left(\begin{array}{cc}
	 	{\bf I}& i{\bf I}\\
		{\bf I} &  -i{\bf I}\\
	\end{array}
	\right)}$ is a unitary matrix \cite{Chap_back_JP}.
So, by indicating as $\bar{g}_r \in \overline{\mathcal{G}}$ (defined in \eqref{set_G_c}) the density generator of the $2m$-dimensional real vector $\overline{\bf x}$ and if 
its scatter matrix $\overline{\boldsymbol \Sigma} \in \mathcal{S}_{2m}^{\mathbb{R}}$ is positive definite, then the pdf \eqref{RES_pdf} of $\mb{x} \in \mathcal{X} \subseteq \mathbb{C}^m$ is generally rewritten in the following form:
\begin{equation}
	\label{non circulaire CES_pdf}
	p_X(\mb{x}|\bs{\mu},\bs{\Sigma},{\bf \Omega},\bar{g}_c)
	=
	|\widetilde{\bf \Sigma}|^{-1/2}\bar{g}_c
	\left[2^{-1}(\widetilde{\bf x}-\widetilde{\boldsymbol \mu})^H \widetilde{\bf \Sigma}^{-1}(\widetilde{\bf x}-\widetilde{\boldsymbol \mu})\right],
\end{equation}
where 
\begin{equation}\label{t_rc_mu}
	\widetilde{\boldsymbol \mu}\triangleq({\boldsymbol \mu}^T,{\boldsymbol \mu}^H)^T=\sqrt{2}{\bf M}\overline{\boldsymbol \mu}
\end{equation}
with
$\overline{\boldsymbol \mu} \in \mathbb{R}^{2m}$ denotes the  location vector of $\overline{\bf x}$,
\begin{equation}\label{t_rc_Sigma}
	\widetilde{\bf \Sigma}
	\triangleq
	\left(\begin{array}{cc}
			{\bf \Sigma}& {\bf \Omega}\\
			{\bf \Omega}^* &  {\bf \Sigma}^* \\
		\end{array}
		\right)=2{\bf M}\overline{\boldsymbol \Sigma}{\bf M}^H,\; \text{and} 
\end{equation}
\begin{equation}\label{t_rc_g}
	\bar{g}_c\triangleq 2^m \bar{g}_r(2t).
\end{equation}
We note that ${\bf \Sigma}=E\{({\bf x}-{\boldsymbol \mu})({\bf x}-{\boldsymbol \mu})^H\}\in {\cal M}_m^{\mathbb{C}}$ 
and ${\bf \Omega}=E\{({\bf x}-{\boldsymbol \mu})({\bf x}-{\boldsymbol \mu})^T\}\in {\cal S}_m^{\mathbb{C}}$ 
where ${\cal M}_m^{\mathbb{C}}$ and ${\cal S}_m^{\mathbb{C}}$ denote the sets of all Hermitian positive definite and complex symmetric matrices, respectively.

Depending on whether ${\bf \Omega}$ is a zero-matrix or not, the CES distribution is called \textit{circular} (C-CES) or \textit{non-circular} (NC-CES). 
Due to its widespread usage, let us have a closer look at the C-CES distributions. Such particular case is characterized by structured scatter matrices 
$\overline{\boldsymbol \Sigma}=\tiny{\left(\begin{array}{cc}
		\bs{\Sigma}_1& -\bs{\Sigma}_2\\
		\bs{\Sigma}_2& \bs{\Sigma}_1
	\end{array}
	\right)}$ where $\bs{\Sigma}_1$ and $\bs{\Sigma}_2$ are symmetric and skew-symmetric, respectively. Moreover, since by definition of circularity, $\bs{\Omega} = \mb{0}$, the C-CES pdf is a particular case of \eqref{non circulaire CES_pdf} that can be explicitly expressed as:
\begin{equation}
	\label{circulaire CES_pdf}
	p_X(\mb{x}|\bs{\mu},\bs{\Sigma},\bar{g}_c)
	=
	|{\bf \Sigma}|^{-1}\bar{g}_c
	\left[({\bf x}-{\boldsymbol \mu})^H {\bf \Sigma}^{-1}({\bf x}-{\boldsymbol \mu})\right].
\end{equation}
Using this complex representation of even-dimensional RES distributions, the Stochastic Representation Theorem \eqref{SRT_dec} can be extended to both circular and non-circular CES distributions as discussed in details in \cite[Sect. 1.3.2]{Chap_back_JP}. These \virg{complex} stochastic representations in fact use the mutually independent, random variable $\mathcal{Q}_c\triangleq \frac{1}{2}\mathcal{Q}_r$ (where 
$\mathcal{Q}_r$ is the 2nd-order modular variate associated with $\overline{\bf x}$) and the random vector ${\bf u}_c\sim \mathcal{U}(S_{\mathbb{C}}^{m-1})$. Moreover, the definition of the shape matrix \eqref{shape_m} can be straightforwardly extended while keeping 
definition of matrix scale function in \eqref{scale_f}.
In particular for NC-CES distributions, the shape matrix 
\begin{equation*}\label{V_comp}
	\widetilde{\bf V}_S \triangleq \widetilde{\bf \Sigma}/S(\widetilde{\bf \Sigma})=
	\left(\begin{array}{cc}
		{\bf \Sigma}_S& {\bf \Omega}_S\\		
		{\bf \Omega}_S^* &  {\bf \Sigma}_S^* \\
	\end{array} \right)
\end{equation*}
is structured like $\widetilde{\bf \Sigma}$.

Although certainly possible, directly rewriting all the results obtained in the previous sections using complex formalism is laborious and adds nothing to the statistical significance of our findings. Below, we therefore limit ourselves to providing two different, yet equivalent recipes that any practitioner can follow to obtain the desired complex-version of FIM and related information bounds:
\begin{description}
	\item[\textit{Recipe 1}]: Real to complex mapping.
	\begin{enumerate}
		\item Take the results obtained in the previous sections and consider the specific case of a real observed vector of dimension $2m$, i.e. $\mathbb{R}^{2m} \ni \overline{\bf x} \sim RES_{2m}(\overline{\boldsymbol \mu},\overline{\boldsymbol \Sigma},\bar{g}_r)$,
		\item Use the transformations given in \eqref{t_rc_mu}, \eqref{t_rc_Sigma} and \eqref{t_rc_g} to map the \virg{real-based} results to the \virg{complex-based} results.
	\end{enumerate}
	\item[\textit{Recipe 2}]: Wirtinger calculus.
	\begin{enumerate}
		\item Consider directly a complex observation vector, i.e. $\mathbb{C}^m \ni {\bf x} \sim \text{\textit{NC-CES}}_{m}({\boldsymbol \mu},{\boldsymbol \Sigma},{\boldsymbol \Omega}, \bar{g}_c)$ in the case of NC-CES distributions or $\mathbb{C}^m \ni {\bf x} \sim \text{\textit{C-CES}}_{m}({\boldsymbol \mu},{\boldsymbol \Sigma}, \bar{g}_c)$ in the case of C-CES distributions.
		\item Use the complex Hilbert space $(\mathcal{H},\innerprod{\cdot}{\cdot}_{\mathcal{H}})$ in \eqref{H_set}, the \virg{complex-aware} inner product $\innerprod{h_1}{h_2}_{\mathcal{H}}\triangleq \ev{h_1h_2^*}$ to express projection operators and tangent spaces,
		\item Use the Wirtinger calculus \cite{Remmert,Complex_M,Kreutz} to handle derivatives with respect to (real or/and complex) parameters.
\end{enumerate}
\end{description}
Extensive discussions and related examples about the above mentioned recipes can be found in \cite{Bos,OllCRB,Complex_MCRB,Menni,DELMAS2015,For_SCRB_complex}.

To conclude this section, let us now take a closer look at the \virg{complex version} of Propositions \ref{Prop_location scatter model} and \ref{Prop_low-rank model} in which the parameter vector $\bs{\theta}_0 = (\bs{\gamma}_0^T, \bs{\xi}_0^T)^T$ is still assumed to be real-valued.

\subsection{Applications of Proposition \ref{Prop_location scatter model}}

Proposition \ref{Prop_location scatter model} extends directly to CES distributions by considering the one-to-one mapping
$\overline{\boldsymbol \mu} \mapsto  \widetilde{\boldsymbol \mu}=\sqrt{2}{\bf M}\overline{\boldsymbol \mu}$. In fact, as shown in \cite{ABEIDA2023}, we have:
\begin{equation}\label{I_gamma_par nc}
	\bar{\mb{I}}(\bs{\gamma}_0|\bs{\xi}_0,\bar{g}_{c,0})
	=\bar{\mb{I}}(\bs{\gamma}_0|\bs{\xi}_0)
	=\mb{I}_{\bs{\gamma}_0}  = \beta(\bar{g}_{c,0})
	\mb{J}[\widetilde{\bs{\mu}}_0]^H\widetilde{\bs{\Sigma}}_0^{-1}\mb{J}[\widetilde{\bs{\mu}}_0],
\end{equation}
for NC-CES distributions where $\beta(\bar{g}_{c,0}) \triangleq \frac{E\{\mathcal{Q}_c\bar{\varphi}^2_{c,0}(\mathcal{Q}_c)\}}{m}$ with $\bar{\varphi}_{c,0}(t) \triangleq \frac{-1}{\bar{g}_{c,0}(t)}\frac{d\bar{g}_{c,0}(t)}{dt}$ and
$[\mb{J}[\widetilde{\bs{\mu}}_0]]_{i,j}  \triangleq  \left. \parder{[\widetilde{\bs{\mu}}(\bs{\theta})]_i}{\theta_j}\right|_{\bs{\theta}=\bs{\theta}_0}$. In the particular case of C-CES distributions, it is immediate to verify that the FIM in \eqref{I_gamma_par nc} reduces to 
\begin{equation}\label{I_gamma_par c}
	\bar{\mb{I}}(\bs{\gamma}_0|\bs{\xi}_0,\bar{g}_{c,0})
	=\bar{\mb{I}}(\bs{\gamma}_0|\bs{\xi}_0)
	=\mb{I}_{\bs{\gamma}_0}  =2 \beta(\bar{g}_{c,0})
	{\rm Re}\graffe{\mb{J}[\bs{\mu}_0]^H\bs{\Sigma}_0^{-1}\mb{J}[\bs{\mu}_0]}
\end{equation}
with $\mb{J}[\bs{\mu}_0]$ given by \eqref{def J}.

In signal processing, there are many examples of complex-valued observations ${\bf x}\in \mathbb{C}^m$ where the location vector $\bs{\mu}$ includes the parameters of interest $\bs{\gamma}$, while the scatter matrix $\bs{\Sigma}$ gathers the nuisance parameters $\bs{\xi}$.
This is the case where a deterministic signal of interest parameterized by $\bs{\gamma}$ is disturbed by a zero-mean C-CES distributed noise, whose density generator 
$\bar{g}$ is unspecified. A classic example is given by the statistical model  for time delay and Doppler estimation problems (see e.g. \cite{FORTUNATI2024,Ortega_For_SPL}).
%
\subsection{Applications of Proposition \ref{Prop_low-rank model}}
Proposition \ref{Prop_low-rank model} also extends to C- and NC-CES distributions by considering the one-to-one mapping
$\overline{\boldsymbol \Sigma} \mapsto \widetilde{\bf \Sigma}
=2{\bf M}\overline{\boldsymbol \Sigma}{\bf M}^H$.
More specifically for the C-CES distribution, \eqref{Cov_mat_DOA} can be rewritten in the following form \cite{Chap_Delmas}:
\begin{equation}\label{Cov_circulaire}
	\bs{\Sigma}_0  = \mb{A}_0\bs{\Xi}_0\mb{A}_0^H + \lambda_0\mb{I}_m,
\end{equation}
where $\bs{\Sigma}_0 \in {\cal S}_m^{\mathbb{C}}$, $\bs{\Xi}_0 \in {\cal S}_p^{\mathbb{C}}$ are two Hermitian and \textit{positive definite} matrices and $\mb{A}_0 \triangleq \mb{A}(\bs{\gamma}_0)  \in \mathbb{C}^{m \times p}$ is full column rank that collects the parameters of interest $\bs{\gamma}_0$ that characterize $\mb{A}_0$.  

For the C-CES distributions, Proposition \ref{Prop_low-rank model} then provides the following equality where the closed-form expression of the parametric efficient FIM has been given in \cite{ABEIDA2023}:
\begin{equation}\label{I_gamma_low-rank c}
	\bar{\mb{I}}(\bs{\gamma}_0|\bs{\xi}_0,\bar{g}_{c,0})
	=\bar{\mb{I}}(\bs{\gamma}_0|\bs{\xi}_0)
	=
	\frac{2\alpha(\bar{g}_{c,0})}{\lambda_0}
	{\rm Re}\graffe{
	{\bf J}[\mathrm{vec}({\bf A}_0)]^H
	({\bf H}_0^T \otimes \boldsymbol{\Pi}^{\bot}_{{\bf A}_0})
	{\bf J}[\mathrm{vec}({\bf A}_0)]},
\end{equation}
where:
\begin{equation}\label{def H0}
	{\bf H}_0 \triangleq \bs{\Xi}_0{\bf A}^H_0{\bf \Sigma}_0^{-1}{\bf A}_0\bs{\Xi}_0,
\end{equation}
\begin{equation}
	\boldsymbol{\Pi}^{\bot}_{{\bf A}_0}
	\triangleq {\bf I}_m-{\bf A}_0[{\bf A}^H_0{\bf A}_0]^{-1}{\bf A}^H_0,
\end{equation}
$[{\bf J}[\mathrm{vec}({\bf A}_0)]]_{i,j} \triangleq  \left. \parder{[\mathrm{vec}({\bf A}(\bs{\gamma}))]_i}{\gamma_j}\right|_{\bs{\gamma}=\bs{\gamma}_0}$ and where
here $\alpha(\bar{g}_{c,0}) \triangleq \frac{E\{\mathcal{Q}_c^2\bar{\varphi}^2_{c,0}(\mathcal{Q}_c)\}}{m(m+1)}$.
Note that the parametric efficient FIM \eqref{I_gamma_low-rank c} reduces to:
\begin{equation}
	\bar{\mb{I}}_\textrm{DOA}(\bs{\gamma}_0|\bs{\xi}_0,\bar{g}_{c,0})
	=\bar{\mb{I}}_\textrm{DOA}(\bs{\gamma}_0|\bs{\xi}_0) = \frac{2\alpha(\bar{g}_{c,0})}{\lambda_0}
	{\rm Re}\graffe{
	({\bf D}_0^H{\bf \Pi}^{\bot}_{{\bf A}_0}{\bf D}_0)
	\odot{\bf H}_0^T}
\end{equation}
for direction of arrival (DOA) modeling with one parameter per source where ${\bf A} \triangleq [{\bf a}_1,...,{\bf a}_p]$ and $({\bf a}_k)_{k=1,...,p}$ are the steering vectors parameterized by the DOA $\gamma_k$ with $\boldsymbol{\gamma}\triangleq (\gamma_1,...,\gamma_p)^T$ and ${\bf D}_0 \triangleq \left. [\frac{d{\bf a}_1}{d\gamma_1},...,\frac{d{\bf a}_p}{d\gamma_p}] \right|_{\bs{\gamma}=\bs{\gamma}_0}$ for $p$ sources.

For the NC-CES distributions, $\widetilde{\bs{\Sigma}}_0$ becomes:
\begin{equation}\label{Cov_noncirculaire}
	\widetilde{\bs{\Sigma}}_0  
	= \left(\begin{array}{cc}
		\mb{A}_0& \mb{0}\\
		\mb{0}& \mb{A}_0^*
	\end{array}
	\right)
	\widetilde{\bs{\Xi}}_0
	\left(\begin{array}{cc}
		\mb{A}^H_0& \mb{0}\\
		\mb{0}& \mb{A}_0^T
	\end{array}
	\right) + \lambda_0\mb{I}_{2m},
\end{equation}
where $\widetilde{\bs{\Xi}}_0=\left(\begin{array}{cc}
		\mb{\Sigma}_{c,0}& \mb{\Omega}_{c,0}\\
		\mb{\Omega}_{c,0}^*& \mb{\Sigma}_{c,0}^*
	\end{array}
	\right)$ is structured like $\widetilde{\bs{\Sigma}}_0$ and \eqref{I_gamma_low-rank c} is also valid, where now ${\bf H}_0$ is replaced by \cite{ABEIDA2023}:
\begin{equation}\label{def H0 r}
{\bf H}_0 \triangleq
\left(\mb{\Sigma}_{c,0}\mb{A}_0^H,\ \mb{\Omega}_{c,0}\mb{A}_0^T \right)
\widetilde{\bf \Sigma}_0^{-1}
\left(\begin{array}{c}
		\mb{A}_0\mb{\Sigma}_{c,0}\\
		\mb{A}_0^* \mb{\Omega}_{c,0}^*
	\end{array}
	\right).
\end{equation}

Note that \eqref{Cov_noncirculaire}  reduces
in the so-called rectilinear case to \cite{Chap_Delmas}
\begin{equation}\label{Cov_rectilinear}
	\widetilde{\bs{\Sigma}}_0  =\widetilde{\mb{A}}_{r,0}\bs{\Xi}_{r,0}\widetilde{\mb{A}}_{r,0}^H + \lambda_0\mb{I}_{2m},
\end{equation}
where $\bs{\Xi}_{r,0} \in {\cal S}_p^{\mathbb{R}}$ is symmetric positive definite and $\widetilde{\mb{A}}_{r,0} = \left(\begin{array}{c}
	\mb{A}_{r,0}\\
	\mb{A}_{r,0}^*
\end{array}
\right) \in \mathbb{C}^{2m \times p}$ is full column rank that collects the parameters of interest 
$\bs{\gamma}_0$ that characterize $\widetilde{\mb{A}}_{r,0}$.
Under the condition $2m>p$, \eqref{I_gamma_low-rank c} is written in the form  \cite{ABEIDA2019}:
\begin{equation}\label{I_gamma_low-rank r}
	\bar{\mb{I}}(\bs{\gamma}_0|\bs{\xi}_0,\bar{g}_{c,0})
	=\bar{\mb{I}}(\bs{\gamma}_0|\bs{\xi}_0)
	=
	\frac{\alpha(\bar{g}_{c,0})}{\lambda_0}
	{\bf J}[\mathrm{vec}(\widetilde{\mb{A}}_{r,0})]^H
	(\widetilde{\bf H}_0^T \otimes \boldsymbol{\Pi}^{\bot}_{\widetilde{\mb{A}}_{r,0}})
	{\bf J}[\mathrm{vec}(\widetilde{\mb{A}}_{r,0})],
\end{equation}
where
$\widetilde{\bf H}_0 \triangleq \bs{\Xi}_{r,0}\widetilde{\mb{A}}_{r,0}^H\widetilde{\bf \Sigma}_0^{-1}\widetilde{\mb{A}}_{r,0}\bs{\Xi}_{r,0}$ and
$\boldsymbol{\Pi}^{\bot}_{\widetilde{\mb{A}}_{r,0}}
	\triangleq {\bf I}_{2m}-\widetilde{\mb{A}}_{r,0}[\widetilde{\mb{A}}_{r,0}^H\widetilde{\mb{A}}_{r,0}]^{-1}\widetilde{\mb{A}}_{r,0}^H$.

This low-rank scatter model encompasses many far or near-field, narrow or wide-band  DOA models
with scalar or vector-sensors for an arbitrary
number of parameters per source and many other models such as the bandlimited SISO, SIMO
\cite{Moulines1995} and MIMO \cite{AbedMeraim1997} channel models.
For example, parametrization \eqref{Cov_rectilinear} can be applied  for DOA estimation modeling  with rectilinear or strictly second-order sources and  for SIMO channels estimation modeling with BPSK or MSK symbols \cite{Delmas2009} where $\bs{\gamma}_0$ represents both the localization parameters (azimuth, elevation, range) and the phase of the sources, and the real and imaginary parts of channel impulse response  coefficients, respectively.
Parametrization \eqref{Cov_noncirculaire} on the other hand  is used for DOA modeling with generally  non-circular and non-rectilinear complex sources.

\textit{Remark}: In all the above-mentioned applications, the observation vector ${\bf x}$ is generally assumed to be the sum of a low rank zero-mean signal of interest
and a zero-mean noise term, which is mutually uncorrelated with the signal. This approach requires that both the statistical models of the useful signal and the noise are chosen a-priori. It should be noted, however, that neither the useful signal nor the noise are observable, so the chosen model
could be completely misspecified. Furthermore, since they are not observable, their model cannot be estimated (in a non-parametric way) from the observed vector $\mb{x}$. Unlike what has been done generally in the literature, we adopt here a semi-parametric statistical model on the observed vector $\mb{x}$ only, without relying on any additional assumption on the statistical model on the unobservable signal of interest and noise. From a statistical point of view, we therefore believe that our approach is more valid and robust to the misspecification of the standard model generally assumed in the literature.
%
\section{Conclusion}\label{Conclusion}
The semiparametric statistical efficiency in estimation problems for elliptically symmetric distributed data was analyzed in this paper. In particular, we studied the impact of finite and infinite-dimensional nuisance parameters can have on the estimation of the parameters of interest which, in the case of elliptical distributions, are contained in the location vector $\bs{\mu}$ and/or in the covariance matrix $\bs{\Sigma}$. The profound and counter-intuitive result that emerged is that, in the presence of specific finite-dimensional nuisance parameters, semiparametric efficiency can be equivalent to parametric efficiency. Specifically, in the case of elliptical distributions, not knowing the density generator does not cause any loss of efficiency when estimating $\bs{\mu}$ or a scaled version of $\bs{\Sigma}$. This result had already been demonstrated by Hallin and Paindaveine in \cite{Hallin_P_2006,PAINDAVEINE} by leveraging an asymptotic approximation of the projection operator as a conditional expectation with respect to an appropriate invariant sub-$\sigma$ algebra. Unlike these works, in this paper we provide an explicit derivation of such projection operator. Furthermore, as an advancement over the state of the art, we analyzed the case, important in many applications, in which the parameters of interest and the finite-dimensional nuisance parameters are given by a parametrization of the location vector and of the covariance matrix. A general condition that the parameterization in question must satisfy in order for the semiparametric efficiency to be equal to the parametric efficiency has been derived in this paper. 
This condition therefore allows us to test this property for any particular parameterizations. Two examples were investigated here,
 including the well-known low-rank parameterization, often arising in many practical signal processing applications. The paper concluded with a section in which the results derived for RES distributions are extended to the case of C-CES and NC-CES distributions. The natural follow-up to this paper will be on the development of semiparametric estimators capable of achieving parametric efficiency. A promising approach for achieving this goal is that of rank-based ($R$-) estimators \cite{Hallin_Annals_Stat_2,Sem_eff_est_TSP,Ortega_For_SPL}.
\clearpage


\clearpage

{\appendix
\section*{Some technical results and their proofs}
\begin{enumerate}[align=left,label*=\thesection.\arabic*,labelindent=2em]
\item \textbf{Some useful results in $\mathcal{H}_q$}
\begin{lemma}\label{lem_p}
Let $\mathcal{H}$ be an Hilbert space and let $\mathcal{U}\subseteq \mathcal{H}$ be a closed linear subspace of $\mathcal{H}$. Let us now introduce the $q$-replicating versions of $\mathcal{H}$ and $\mathcal{U}$	as $\mathcal{H}_q = \mathcal{H} \times \cdots \times \mathcal{H}$ and $\mathcal{H}_q \supseteq \mathcal{U}_q = \mathcal{U} \times \cdots \times \mathcal{U}$. For $\mb{h} \in \mathcal{H}_q$ and $h_i \in \mathcal{H}$, $i=1,\ldots,q$, we have that:
\begin{equation}\label{eq_p_Hq}
	\quadre{ \Pi(\mb{h}|\mathcal{U}_q)}_i = \Pi(h_i|\mathcal{U}),\; i = 1,\ldots,q,
\end{equation}
where $\Pi(\mb{h}|\mathcal{U}_q) \in \mathcal{U}_q$ and $\Pi(h_i|\mathcal{U}) \in \mathcal{U}, \; i = 1,\ldots,q$.
\end{lemma}
\begin{IEEEproof}
	For a given $\mb{h} \in \mathcal{H}_q$, the projection $\Pi(\mb{h}|\mathcal{H}_q)$ is defined as the vector in $\mathcal{H}_q$ such that \cite[Theo. 5.1]{Chap_sem}:
	\begin{equation}
		\norm{\mb{h}-\Pi(\mb{h}|\mathcal{H}_q)}_{\mathcal{H}_q} \leq \norm{\mb{h}-\mb{u}}_{\mathcal{H}_q},\; \forall \mb{u} \in \mathcal{U}_q.
	\end{equation}
From the definition of the inner product in $\mathcal{H}_q$, induced by the one in $\mathcal{H}$ as $\innerprod{\mb{h}}{\mb{u}}_{\mathcal{H}_q} \triangleq \sum_{i=1}^q \innerprod{h_i}{u_i}_{\mathcal{H}}$, we have that:
 \begin{equation}
		\begin{split}
			\norm{\mb{h}-\mb{u}}_{\mathcal{H}_q} &= \sqrt{\innerprod{\mb{h}-\mb{u}}{\mb{h}-\mb{u}}}_{\mathcal{H}_q}\\
			&=\sqrt{\sum\nolimits_{i=1}^{q}\innerprod{h_i-u_i}{h_i-u_i}_{\mathcal{H}}}=\sqrt{\sum\nolimits_{i=1}^{q}\norm{h_i-u_i}_{\mathcal{H}}}.
		\end{split}
	\end{equation}
	As a consequence, minimizing $\norm{\mb{h}-\mb{u}}_{\mathcal{H}_q}$ is equivalent to minimize each term $\norm{h_i-u_i}_{\mathcal{H}}$. Then, the equality \eqref{eq_p_Hq} follows from the definition of orthogonal projection in $\mathcal{H}$ onto $\mathcal{U}$. 
\end{IEEEproof}
\vspace{0.5cm}

\begin{theorem}\label{Pyt_theo}
	The $q$-replicating Pythagorean theorem \cite[Theo. 3.3]{Tsiatis}: Let $\mathcal{H}_q$ and $\mathcal{U}_q \subseteq \mathcal{H}_q$ the $q$-replicating Hilbert space and subspace defined as in Lemma \ref{lem_p}. Let us take $\mb{h} \in \mathcal{H}_q$ and $\mb{u} \in \mathcal{U}_q$. If $\mb{h} \perp \mb{u}$, then:
	\begin{equation}\label{Pyt_theo_eq}
		\mb{G}(\mb{h} + \mb{u},\mb{h} + \mb{u}) = \mb{G}(\mb{h},\mb{h}) + \mb{G}(\mb{u},\mb{u}),
	\end{equation} 
	where $[\mb{G}(\mb{a},\mb{b})]_{i,j} \triangleq\innerprod{a_i}{b_j}_{\mathcal{H}}$. 
\end{theorem}
\begin{IEEEproof}
	Let us start by showing that:
	\begin{equation}\label{int_rel}
		\mathcal{H}_q \ni \mb{h} \perp \mb{u} \in \mathcal{U}_q \Leftrightarrow \mathcal{H} \ni h_i \perp u_j \in \mathcal{U},\; i,j=1,\ldots,q.
	\end{equation}
	\begin{enumerate}
		\item Proof of the implication $\Leftarrow$. Since $\mathcal{H} \ni h_i \perp u_j \in \mathcal{U}$, $i,j=1,\ldots,q$, we have that $\innerprod{h_i}{u_j}_{\mathcal{H}}=0$, $i,j=1,\ldots,q$. Then, $\innerprod{\mb{h}}{\mb{u}}_{\mathcal{H}_q}=\sum\nolimits_{i=1}^q\innerprod{h_i}{u_i}_{\mathcal{H}} =0$. Consequently, $\mb{h} \perp \mb{u}$.
		\item Proof of the implication $\Rightarrow$. We notice that $\mb{h} \perp \mb{u} \Leftrightarrow \Pi(\mb{h}|\mathcal{U}_q)=\mb{0}$. Then, from Lemma \ref{lem_p}, $\mb{h} \perp \mb{u} \Leftrightarrow \Pi(h_i|\mathcal{U})=0$, $i=1,\ldots,q$ and consequently $h_i \perp u_j \in \mathcal{U}$, $i,j=1,\ldots,q$.
	\end{enumerate}
	Now that we know that \eqref{int_rel} holds true, we can deduce that:
	\begin{equation}
		\innerprod{\mb{h}}{\mb{u}}_{\mathcal{H}_q}=0 \Leftrightarrow \innerprod{h_i}{u_j}_{\mathcal{H}}=0,\; i,j=1,\ldots,q,
	\end{equation}
	or equivalently that $\innerprod{\mb{h}}{\mb{u}}_{\mathcal{H}_q}=0 \Leftrightarrow\mb{G}(\mb{h},\mb{u}) = \mb{G}(\mb{u},\mb{h}) = \mb{0}$. Finally, \eqref{Pyt_theo_eq} follows from a simple calculation:
	\begin{equation}
		\begin{split}
			\mb{G}&(\mb{h} + \mb{u},\mb{h} + \mb{u})\\
			 &= \mb{G}(\mb{h},\mb{h}) + \mb{G}(\mb{h},\mb{u}) +\mb{G}(\mb{u},\mb{h})+ \mb{G}(\mb{u},\mb{u})\\
			 & = \mb{G}(\mb{h},\mb{h}) + \mb{G}(\mb{u},\mb{u}).
		\end{split}
	\end{equation}
\end{IEEEproof}
 \label{App1}
\vspace{0.5cm}

\item {\bf Proof of Lemma \ref{Lem1}.} From \eqref{eff_score_nu}, it follows that $\bar{\mb{t}}_{\bs{\gamma}_0} =  \bar{\mb{s}}_{\bs{\gamma}_0} + \mb{p}$ and therefore the covariance matrix of $\bar{\mb{t}}_{\bs{\gamma}_0}$ breaks down:
\begin{equation}\label{Pythag_C}
E_0\graffe{\bar{\mb{t}}_{\bs{\gamma}_0}\bar{\mb{t}}_{\bs{\gamma}_0}^T}
 =\bar{\mb{I}}(\bs{\gamma}_0|\bs{\xi}_0)
=
\bar{\mb{I}}(\bs{\gamma}_0|\bs{\xi}_0,g_0)
 +  \ev{\mb{p}\mb{p}^T} + \ev{\bar{\mb{s}}_{\bs{\gamma}_0}\mb{p}^T} + \quadre{\ev{\bar{\mb{s}}_{\bs{\gamma}_0}\mb{p}^T}}^T.
\end{equation}

Then the equality in \eqref{Sem_FIM_eff} follows immediately iff $E_0\graffe{\bar{\mb{s}}_{\bs{\gamma}_0}\mb{p}^T} = \mb{0}$. From the component-wise application of the inner product, we have that:
\begin{equation}\label{cond_orth}
E_0\graffe{\bar{\mb{s}}_{\bs{\gamma}_0}\mb{p}^T} = \mb{0}\Leftrightarrow  \innerprod{[\bar{\mb{s}}_{\bs{\gamma}_0}]_i}{p_j}_\mathcal{H} = 0 \Leftrightarrow [\bar{\mb{s}}_{\bs{\gamma}_0}]_i \perp p_j,\ \forall i,j\in \{1,\ldots,q\}.
\end{equation}
To show that $[\bar{\mb{s}}_{\bs{\gamma}_0}]_i$ is orthogonal to $p_j$, $\forall i,j$, we note that, according to its definition given in \eqref{bar_s_nu}, $[\bar{\mb{s}}_{\bs{\gamma}_0}]_i$ is the residual of $[\mb{s}_{\bs{\gamma}_0}]_i$ after projection onto $\mathcal{T}_{2} + \mathcal{T}_{3}$. As a direct consequence, we have that $[\bar{\mb{s}}_{\bs{\gamma}_0}]_i \perp (\mathcal{T}_{2} + \mathcal{T}_{3})$, or equivalently, $[\bar{\mb{s}}_{\bs{\gamma}_0}]_i \in (\mathcal{T}_{2} + \mathcal{T}_{3})^\perp$. Moreover, again by definition in \eqref{def_w}, $p_j\in (\mathcal{T}_{2} + \mathcal{T}_{3}) \cap \mathcal{T}_{2}^\perp \subseteq \mathcal{T}_{2} + \mathcal{T}_{3}$. Consequently, $\innerprod{[\bar{\mb{s}}_{\bs{\gamma}_0}]_i}{p_j} = 0$, $\forall i,j$ since $(\mathcal{T}_{2} + \mathcal{T}_{3}) \perp (\mathcal{T}_{2} + \mathcal{T}_{3})^\perp$.
\hfill
\QED
\label{App2}
\vspace{0.5cm}

\item {\bf Proof of Lemma \ref{Lem2}.}
Let us start by noticing that the closed subspaces $(\mathcal{T}_{2} + \mathcal{T}_{3}) \cap \mathcal{T}_{2}^\perp$ and 
$\mathcal{T}_{2}^\perp \cap \mathcal{T}_{3}^\perp$ are orthogonal. In fact, it is immediate to verify that, each element $h \in \mathcal{T}_{2} + \mathcal{T}_{3}$ can be written as $h=h_2+h_3$ with $h_2\in \mathcal{T}_{2}$ and $h_3\in \mathcal{T}_{3}$. Moreover, $\forall g \in \mathcal{T}_{2}^\perp \cap \mathcal{T}_{3}^\perp$ we have that $g \in \mathcal{T}_{2}^\perp$ and $g \in \mathcal{T}_{3}^\perp$ and consequently $h \perp g$. In addition, $\mathcal{T}_{2}^\perp$ can be expressed as the (direct) sum of the  and these two closed subspaces, i.e. $\mathcal{T}_{2}^\perp = (\mathcal{T}_{2} + \mathcal{T}_{3}) \cap \mathcal{T}_{2}^\perp + (\mathcal{T}_{2}^\perp \cap \mathcal{T}_{3}^\perp)$. Consequently, from the property \eqref{orthogonal projector}, we get:
\begin{equation}\label{proj_int}
		\Pi(\bar{\mb{t}}_{\bs{\gamma}_0}|(\mathcal{T}_{2} + \mathcal{T}_{3}) \cap \mathcal{T}_{2}^\perp) = \Pi(\bar{\mb{t}}_{\bs{\gamma}_0}|\mathcal{T}_{2}^\perp) -  \Pi(\bar{\mb{t}}_{\bs{\gamma}_0}|\mathcal{T}_{2}^\perp \cap \mathcal{T}_{3}^\perp).
\end{equation} 

By its definition, given in \eqref{t_def}, $[\bar{\mb{t}}_{\bs{\gamma}_0}]_i \in  \mathcal{T}_{2}^\perp$. As a consequence, we have that the first projection in the RHS of \eqref{proj_int} can be evaluated as $\Pi(\bar{\mb{t}}_{\bs{\gamma}_0}|\mathcal{T}_{2}^\perp) = \bar{\mb{t}}_{\bs{\gamma}_0}$.
And therefore $\mb{p}=\bs{0} \Leftrightarrow  \bar{\mb{t}}_{\bs{\gamma}_0}=\Pi(\bar{\mb{t}}_{\bs{\gamma}_0}|\mathcal{T}_{2}^\perp \cap \mathcal{T}_{3}^\perp)
\Leftrightarrow  [\bar{\mb{t}}_{\bs{\gamma}_0}]_i  \in \mathcal{T}_{2}^\perp \cap \mathcal{T}_{3}^\perp \Leftrightarrow [\bar{\mb{t}}_{\bs{\gamma}_0}]_i  \in \mathcal{T}_{2}^\perp$ and $[\bar{\mb{t}}_{\bs{\gamma}_0}]_i  \in  \mathcal{T}_{3}^\perp$. Then because $[\bar{\mb{t}}_{\bs{\gamma}_0}]_i  \in \mathcal{T}_{2}^\perp$,
 $\mb{p}=\bs{0} \Leftrightarrow [\bar{\mb{t}}_{\bs{\gamma}_0}]_i  \in  \mathcal{T}_{3}^\perp$ or equivalently, $[\bar{\mb{t}}_{\bs{\gamma}_0}]_i \perp \mathcal{T}_3$.
\hfill
\QED
\label{App31}
\vspace{0.5cm}

\item {\bf Projection onto two subspaces $\mathcal{D} + \mathcal{T}$, where $\mathcal{D}$ has finite dimension.} The following theorem has been proposed in \cite[Theo. 5, A.4]{BKRW}.
Let $(\mathcal{H},\innerprod{\cdot}{\cdot}_\mathcal{H})$ and Hilbert space. Let $\mathcal{T} \subset \mathcal{H}$ be an infinite-dimensional subspace while $\mathcal{H} \supset \mathcal{D} = \mathrm{Span}\{s_1, \ldots, s_d\}$ is a \textit{finite $d$-dimensional} subspace. Then, the orthogonal projection of an element $h \in \mathcal{H}$ onto $\mathcal{D} + \mathcal{T}$ can be expressed as:
\begin{equation}
	\Pi(h|\mathcal{D} + \mathcal{T}) = \Pi(h|\mathcal{T}) + \Pi(h|\mathrm{Span}\{r_1,\ldots, r_k\}),
\end{equation}
where $r_i \triangleq s_i-\Pi(s_i|\mathcal{T})$ and $k \triangleq \mathrm{dim}( \mathrm{Span}\{r_1,..,r_k\})\leq d$. Moreover $k = d$ \textit{if and only if} $\mathcal{D} \cap \mathcal{T} = \{0\}$. 

\begin{IEEEproof}
To prove this expression, let us start by noticing that any element in $\{s_i\}_{i=1}^d$ can be decomposed as:
\begin{equation}
	s_i  = [\Pi(s_i|\mathcal{T})] + [s_i-\Pi(s_i|\mathcal{T})] \triangleq t_i + r_i,
\end{equation} 
where, by definition, $t_i \in  \mathcal{T}$ and $r_i \in \mathcal{T}^\perp$, $i = 1,\ldots,d$. As a consequence, any element $l \in \mathcal{D}$ can be decomposed as:
\begin{equation}
	l = \sum\nolimits_{i=1}^d \alpha_i t_i + \sum\nolimits_{i=1}^d \alpha_i r_i, \; \alpha_i \in \mathbb{R},
\end{equation}
where the first term belongs to $\mathcal{T}$ while the second belongs to $\mathrm{Span}\{r_1, \ldots, r_d\} \subset \mathcal{T}^\perp$ and then:
\begin{equation}
	\mathcal{D} + \mathcal{T} = \mathrm{Span}\{r_1, \ldots, r_d\} + \mathcal{T}.
\end{equation}
In general, the elements $\{r_1, \ldots, r_d\}$ may be linearly dependent, so we can extract a basis $\{r_{i_1}, \ldots, r_{i_k}\}$ such that $\mathrm{Span}\{r_1, \ldots, r_d\} = \mathrm{Span}\{r_{i_1}, \ldots, r_{i_k}\}$. For the sake of simplicity and without loss of generality, let us suppose the $k$ linearly independent elements are the first ones, then we have:
\begin{equation}
	\mathcal{D} + \mathcal{T} = \mathrm{Span}\{r_1, \ldots, r_k\} + \mathcal{T}.
\end{equation}
Finally, according to \cite[eq. (5.7)]{Chap_sem} or \cite[Sect. 2.4, Ex. 1]{Tsiatis}, the following explicit expression holds:
\begin{equation}\label{proj_sum}
	\Pi(h|\mathcal{D} + \mathcal{T}) = \Pi(h|\mathcal{T}) + \Pi(h|\mathrm{Span}\{r_1, \ldots, r_k\}) 
	= \Pi(h|\mathcal{T}) + \mb{p}^T\mb{G}^{-1}\mb{r},
\end{equation}
where $\mb{p} \triangleq (p_1, \ldots,p_k)^T$ with $p_i \triangleq \innerprod{h}{r_i}_\mathcal{H}$, $[\mb{G}]_{ij} \triangleq  \innerprod{r_i}{r_j}_\mathcal{H}$, for $i = 1,\ldots,k$ and $\mb{r} \triangleq (r_1,\ldots, r_k)^T$. Let us now focus on finding a geometric (necessary and sufficient) condition to have $k=d$, that is verified if and only if $\{r_1, \ldots, r_d\}$ are linearly independents. By definition, we have that:
\begin{equation}
	\mathrm{Span}\{r_1, \ldots, r_d\} = \Pi(\mathcal{D}|\mathcal{T}^\perp).
\end{equation}
Specifically, we can define a linear map $P$ as:
\begin{equation}
	\begin{split}
		P :\;& \mathcal{D} \rightarrow \mathcal{T}^\perp\\
		& l \mapsto \Pi(l|\mathcal{T}^\perp).
	\end{split}
\end{equation}
Therefore $k = \mathrm{dim}(\mathrm{Span}\{r_1, \ldots, r_d\}) = \mathrm{dim}(\mathrm{Im}\{P\})$. Form the rank-nullity theorem, we also have that $\mathrm{dim}(\mathrm{Im}\{P\}) + \mathrm{dim}(\mathrm{ker}\{P\}) = d$. Consequently:
\begin{equation}
	k = d \Leftrightarrow \mathrm{dim}(\mathrm{ker}\{P\}) = 0.
\end{equation}
Since, by definition $\mathrm{ker}\{P\} = \{l \in \mathcal{D} | \Pi(l|\mathcal{T}^\perp) = 0\} = \{l \in \mathcal{D} | \Pi(l|\mathcal{T}) = l\}$, we immediately have that $\mathrm{ker}\{P\} = \mathcal{D} \cap \mathcal{T}$ and consequently:
\begin{equation}
	k = d \Leftrightarrow \mathcal{D} \cap \mathcal{T} =\{0\}.
\end{equation} 
\end{IEEEproof}
\label{App_Theo_sum}
\vspace{0.5cm}

\label{App3}
\item {\bf Implicit derivation of $\bs{\nabla}_{\ovecs{\mb{V}_S}}^T [\mb{V}_S]_{11}$.} The following calculation has been derived in \cite[Sect. 4]{PAINDAVEINE}.
Let us start by defining the mapping $v_{11}^S : \mathbb{R}^{m(m+1)/2-1} \mapsto\  \mathbb{R}$ implicitly defined by the constraint $S(\mb{V}_S) \equiv S(v_{11}^S(\ovecs{\mb{V}_S}),\ovecs{\mb{V}_S}) = 1$. Thanks to the implicit function theorem, under Assumptions A1, A2, A3, this mapping exists, is unique and continuously differentiable around a given $\vecs{\mb{V}_S}$. Then, we can differentiate both side of $S(v_{11}^S(\ovecs{\mb{V}_S}),\ovecs{\mb{V}_S}) = 1$ w.r.t. $\ovecs{\mb{V}_S}$ to get:
\begin{equation}
	\begin{split}
		\bs{\nabla}_{\ovecs{\mb{V}_S}}S(v_{11}^S,\ovecs{\mb{V}_S}) + \parder{S(v_{11}^S,\ovecs{\mb{V}_S})}{v_{11}^S}\nabla_{\ovecs{\mb{V}_S}}v_{11}^S(\ovecs{\mb{V}_S}) = \mb{0},
	\end{split}
\end{equation}  
then consequently
\begin{equation}
	\begin{split}
		\nabla_{\ovecs{\mb{V}_S}}v_{11}^S(\ovecs{\mb{V}_S}) = - \frac{\nabla_{\ovecs{\mb{V}_S}}S(v_{11}^S,\ovecs{\mb{V}_S})}{\partial S(v_{11}^S,\ovecs{\mb{V}_S})/\partial v_{11}^S}.
	\end{split}
\end{equation}
\hfill
\QED
\label{App4}
\vspace{0.5cm}

\item {\bf Explicit calculation of $\nabla_{\ovecs{\mb{V}_S}}^T [\mb{V}_S]_{11}$ for $S(\bs{\Sigma}) = [\bs{\Sigma}]_{1,1}$, $S(\bs{\Sigma}) = \trace{\bs{\Sigma}}/m$ and $S(\bs{\Sigma}) = |\bs{\Sigma}|^{1/m}$.} For the scale $S(\bs{\Sigma}) = [\bs{\Sigma}]_{1,1}$ it is trivial to verify that $\nabla_{\ovecs{\mb{V}_S}}^T [\mb{V}_S]_{11}=\mb{0}^T_{m(m+1)/2-1}$. For the scale $S(\bs{\Sigma}) = \trace{\bs{\Sigma}}/m$, let us start by noticing that it implies $\trace{\mb{V}_S}=m$. By taking the differential on both side of this equality, and by expliting the linearity of the trace, we get:
\begin{equation}\label{dtV}
	d\trace{\mb{V}_S}=\trace{\mb{I}_m d\mb{V}_S} = \cvec{\mb{I}_m}^T\cvec{d\mb{V}_S} =0.
\end{equation} 
Now, by definition of $\mb{D}_m$:
\begin{equation}\label{cvI}
	\cvec{\mb{V}_S}=\mb{D}_m\cvec{\mb{V}_S} = \mb{D}_m\left(\begin{array}{c}
		 [\mb{V}_S]_{11}\\
		\ovecs{\mb{V}_S}
	\end{array}\right).
\end{equation}
Then, by substituting \eqref{cvI} in \eqref{dtV}, we obtain:
\begin{equation}
	\begin{split}
		 \cvec{\mb{I}_m}^T&\mb{D}_m\left(\begin{array}{c}
		 	d[\mb{V}_S]_{11}\\
		 	d\ovecs{\mb{V}_S}
		 \end{array}\right) = (\cvec{\mb{I}_m}^T\mb{D}_m)_1d[\mb{V}_S]_{11}+ \cvec{\mb{I}_m}^T\mb{D}_m\underline{\mb{I}}_m^Td\ovecs{\mb{V}_S} =0,
	\end{split}
\end{equation}
from which we immediately get: 
\begin{equation}
	\tonde{\frac{d[\mb{V}_S]_{11}}{d\ovecs{\mb{V}_S}}}^T \triangleq \nabla_{\ovecs{\mb{V}_S}}^T [\mb{V}_S]_{11} =-\frac{\cvec{\mb{I}_m}^T\mb{D}_m\underline{\mb{I}}_m^T}{[\cvec{\mb{I}_m}^T\mb{D}_m]_1}.
\end{equation}
The scale $S(\bs{\Sigma}) = |\bs{\Sigma}|^{1/m}$ implies that $|\mb{V}_S|=1$. Taking as before the differential on both side, we have \cite[p. 149]{Magnus}:
\begin{equation}
	|\mb{V}_S|\trace{\mb{V}_S^{-1} d\mb{V}_S} = 0 \Rightarrow \cvec{\mb{V}_S^{-1}}^T\cvec{d\mb{V}_S} =0.
\end{equation} 
Consequently, by applying exactly the same procedure as before, we get:
\begin{equation}
	\nabla_{\ovecs{\mb{V}_S}}^T [\mb{V}_S]_{11} = -\frac{\cvec{\mb{V}_S^{-1}}^T\mb{D}_m\underline{\mb{I}}_m^T}{[\cvec{\mb{V}_S^{-1}}^T\mb{D}_m]_1}.
\end{equation}  
\hfill
\QED
\label{App5}
\vspace{0.5cm}

\item {\bf Properties of the matrix $\mb{M}_S^{\bs{\Sigma}}$.} The following properties have been derived in \cite{PAINDAVEINE,Hallin_P_2006}. Let us define the matrix $\mb{D}_S^{\bs{\Sigma}}$ as the matrix derivative of the scale function $S$ defined in \eqref{scale_f}:
\begin{equation}\label{D_mat}
	\mb{D}_S^{\bs{\Sigma}} \triangleq \parder{S(\bs{\Sigma})}{\bs{\Sigma}} \in  \mathbb{R}^{m \times  m}.
\end{equation}
We note that, for $S(\bs{\Sigma}) = [\bs{\Sigma}]_{11}$, $S(\bs{\Sigma}) = \trace{\bs{\Sigma}}/m$ and $S(\bs{\Sigma}) = |\bs{\Sigma}|^{1/m}$, we have $\mb{D}_S^{\bs{\Sigma}} = \mb{e}_{1,m}\mb{e}_{1,m}^T$, $\mb{D}_S^{\bs{\Sigma}} = m^{-1}\mb{I}_m$ and $\mb{D}_S^{\bs{\Sigma}} = m^{-1}|\bs{\Sigma}|^{1/m}\bs{\Sigma}^{-1}$ respectively.

The matrices $\mb{D}_S^{\bs{\Sigma}}$ and $\mb{M}_S^{\bs{\Sigma}}$ have the following properties \cite{PAINDAVEINE}:
\begin{itemize}
	\item[P1] From the 1-homogeneity of $S$, we have that $\mb{D}_S^{c\bs{\Sigma}} = \mb{D}_S^{\bs{\Sigma}}$ for all $c>0$. Moreover $S(\bs{\Sigma}) = \trace{\mb{D}_S^{\bs{\Sigma}}\bs{\Sigma}}$.
	\item[P2] Let $\mb{A}$ be a $m \times m$ symmetric, real matrix. If $\trace{\mb{D}_S^{\bs{\Sigma}}\mb{A}}=0$, then 
	\begin{equation}\label{M_prop}
		[\mb{M}_S^{\bs{\Sigma}}]^T\ovecs{\mb{A}} = \mb{D}_m\mb{K}_{\mb{V}_{S,0}} \ovecs{\mb{A}} = \mb{D}_m\vecs{\mb{A}} = \cvec{\mb{A}}.
	\end{equation} 
	\item[P3] $\mb{M}_S^{\bs{\Sigma}}$ has full row rank equal to $m(m+1)/2-1$.
	\item[P4] $\text{Ker}\;\mb{M}_S^{\bs{\Sigma}} \cap \cvec{\mathcal{S}_m^\mathbb{R}} = \{a \cdot \cvec{\mb{D}_S^{\bs{\Sigma}}}, a \in \mathbb{R}\}$, that is the null space of the restriction to the vectorized space $\mathcal{S}_m^\mathbb{R}$ of the real, symmetric matrices of dimension $m \times m$ of the linear application defined by $\mb{M}_S^{\bs{\Sigma}}$ is the one-dimensional space generated by $\cvec{\mb{D}_S^{\bs{\Sigma}}}$.
\end{itemize}
The proof of these properties can be found in \cite[Lemma 4.2]{PAINDAVEINE}. Let us now focus on the property P4 to add some insight on the image of $\mb{M}_S^{\bs{\Sigma}}$. From the Property P3, we have that the dimension of $\text{Im}\;\mb{M}_S^{\bs{\Sigma}} \cap \cvec{\mathcal{S}_m^\mathbb{R}}$, i.e. the number of the linearly independent rows of $\mb{M}_S^{\bs{\Sigma}}$ is equal to $m(m+1)/2-1$. Moreover, as a direct consequence of the rank-nullity theorem, we have the following additional property:
\begin{itemize}
	\item[P5] An orthonormal basis of the $m(m+1)/2-1$-dimensional image $\text{Im}\;\mb{M}_S^{\bs{\Sigma}} \cap \cvec{\mathcal{S}_m^\mathbb{R}}$ of the restriction to $\cvec{\mathcal{S}_m^\mathbb{R}}$ of the linear application defined by $\mb{M}_S^{\bs{\Sigma}}$ is given by the columns of the matrix $\mb{V}_{\bs{\Sigma}} = \mb{D}_m \mb{U}_{\bs{\Sigma}} \in \mathbb{R}^{m^2 \times  (m(m+1)/2-1)}$ such that :
	\begin{equation}\label{U_mat}
		\vecs{\mb{D}_S^{\bs{\Sigma}}}^T \mb{U}_{\bs{\Sigma}} = \mb{0},\quad \mb{U}_{\bs{\Sigma}}^T\mb{U}_{\bs{\Sigma}} = \mb{I}_{m(m+1)/2-1}.
	\end{equation} 
	From its definition given in \eqref{M_mat}, we have that $\quadre{\mb{M}_S^{\bs{\Sigma}}}^T = \mb{D}_m\mb{K}_{\bs{\Sigma}}$. Consequently the column of the matrix $\mb{K}_{\bs{\Sigma}}$ and the one of $\mb{U}_{\bs{\Sigma}}$ span the same subspace. In other word, for each $\bs{\Sigma} \in \mathcal{S}_m^\mathbb{R}$, there exists an invertible transformation matrix $\mb{S}_{\bs{\Sigma}}$ such that:
	\begin{equation}\label{change_coor}
		\mb{K}_{\bs{\Sigma}}\mb{S}_{\bs{\Sigma}} = \mb{U}_{\bs{\Sigma}}.
	\end{equation}
	We note, for further reference, that for $\mb{V}=s^{-1}\bs{\Sigma}$ from the Property P1, we have $\vecs{\mb{D}_S^{\bs{\Sigma}}}^T \mb{U}_{\bs{\Sigma}} = 	\vecs{\mb{D}_S^{\mb{V}}}^T \mb{U}_{\mb{V}}$ and then $\mb{U}_{\mb{V}} = \mb{U}_{\bs{\Sigma}}$ and consequently: 
	\begin{equation}\label{K_S_U}
		\mb{K}_{\mb{V}}\mb{S}_{\mb{V}} = \mb{U}_{\mb{V}}.
	\end{equation}
\end{itemize}  
\hfill
\QED
\label{App6}
\vspace{0.5cm}


\item \textbf{Derivation of the nuisance tangent spaces $\mathcal{T}^u_{g_0}$ and $\mathcal{T}^c_{\bar{g}_0}$.}
Let us start by considering the \textit{unconstrained} non-parametric sub-model $\mathcal{P}^u_{g}$ given in \eqref{unc_non_par} and derive the relevant tangent space $\mathcal{T}^u_{g_0}$. To this end, following \cite[Sec. 4.2]{Tsiatis}, \cite[Sec. 3.1]{BKRW},\cite[Sec. 2.2]{Hallin_Werker}, we need to introduce a set of \textit{parametric sub-models} of $\mathcal{P}^u_{g}$. Formally, an $i$-\textit{th} parametric sub-model of $\mathcal{P}^u_{g}$ is defined as: 
\begin{equation}
	\label{sub_par}
	\mathcal{P}^u_{\tau_{i,\bs{\rho}}}  \triangleq \graffe{p_X(\mb{x}|\bs{\gamma}_0,\tau_{i,\bs{\rho}}), \bs{\rho} \in \Upsilon_i \subseteq \mathbb{R}^{r_i}},
\end{equation}
where:
\begin{equation}
	\begin{split}
		\tau_{i,\bs{\rho}}: \, & \mathcal{X} \times \Upsilon_i \rightarrow \mathcal{G} \\
		& \bs{\rho} \mapsto \tau_i(\mb{x},\bs{\rho}),
	\end{split}
\end{equation}
is a \textit{known} function parametrized by an (artificial) \textit{unknown} finite-dimensional vector $\bs{\rho}$. In particular, for every smooth parametric map $\tau_i(\cdot,\bs{\rho}) : \Upsilon_i \rightarrow \mathcal{G}$, $\mathcal{\mathcal{P}}^u_{\tau_{i,\bs{\rho}}}$ in \eqref{sub_par} is a parametric model satisfying the following two conditions \cite[Sec. 4.2]{Tsiatis}:
\begin{description}
	\item[{\normalfont C1)}] $\mathcal{P}^u_{\tau_{i,\bs{\rho}}} \subseteq \mathcal{P}^u_{g}, \; \forall i \in \mathbb{N}$,
	\item[{\normalfont C2)}] $p_0(\mb{x}) \in \mathcal{P}^u_{\tau_{i,\bs{\rho}}}$, i.e. $\forall i \in \mathbb{N}$ there exists a vector $\bs{\rho}_{0} \in \Upsilon_i$ such that $p_X(\mb{x}|\bs{\gamma}_0,\tau_{i,\bs{\rho}_0}) = p_X(\mb{x}|\bs{\gamma}_0,g_0)\triangleq p_0(\mb{x})$.
\end{description}
Intuitively, a parametric sub-model $\mathcal{P}^u_{\tau_{i,\bs{\rho}}}$ can be thought as a finite-dimensional approximation of the non-parametric model $\mathcal{P}^u_{g}$. The purpose of using a parametric sub-model lies in the fact that its tangent space is well-defined as shown in \eqref{ts_2} as:
\begin{equation}
	\label{ts_sub}
	\mathcal{H}  \supseteq  \mathcal{T}^u_{i,\bs{\rho}_0} \triangleq \mathrm{Span}\{[\mb{s}_{\bs{\rho}_0,i}]_1, \ldots, [\mb{s}_{\bs{\rho}_0,i}]_{r_i}\},
\end{equation} 
where $\mb{s}_{\bs{\rho}_0,i} =\nabla_{\bs{\rho}} \ln p_X(\mb{x}|\bs{\gamma}_0,\tau_{i,\bs{\rho}_0})$ is the score vector of $\bs{\rho}_{0} \in \Upsilon_i$ in the $i$-th parametric sub-model $\mathcal{P}_{\tau_i}^u$. Consequently, according to \cite{Newey} and \cite[Sec. 4.4]{Tsiatis}, the tangent space $\mathcal{T}^u_{g_0}$ can be defined as the closure in $\mathcal{H}$ \footnote{The closure $\overline{\mathcal{A}}$ of a set $\mathcal{A}$ is defined as the smallest closed set that contains $\mathcal{A}$, or equivalently, as the set of all elements in $\mathcal{A}$ together with all the limit points of $\mathcal{A}$.} of the union of all the (parametric) tangent spaces $\mathcal{T}^u_{i,\bs{\rho}_{0}}$:
\begin{equation}
	\label{semi_tangent_space_def}
	\mathcal{T}^u_{g_0} =  \overline{\bigcup\nolimits_{i \in \mathbb{N}}\mathcal{T}^u_{i,\bs{\rho}_{0}}} \subseteq \mathcal{H}.
\end{equation}
Equivalently, $\mathcal{T}^u_{g_0} \subseteq \mathcal{H}$ is the subspace \footnote{The closure of a linear spaces doesn't need to be linear, in general. However, as discussed in \cite[Assumption S]{Begun} and \cite[Sec. 4.4, Remark 5]{Tsiatis}), $\mathcal{T}^u_{g_0}$ is a linear space in the vast majority of the non-pathological statistical models.} of $\mathcal{H}$ composed by all the functions $h \in  \mathcal{H}$ for which there exists a sequence 
$\{\mb{c}_i^T\mb{s}_{\bs{\rho}_0,i}, \mb{c}_i \in \mathbb{R}^{r_i}\}_{i \in \mathbb{N}}$  such that $\norm{h-\mb{c}_i^T\mb{s}_{\bs{\rho}_{0,i}}}^2 = E_0\graffe{(h-\mb{c}_i^T\mb{s}_{\bs{\rho}_{0,i}})^2} \rightarrow 0$. 

Now that we have the theoretical and formal framework, let us go back to the application at hand. Specifically, we have to show that $\mathcal{T}^u_{g_0}$ can actually be expressed as in \eqref{Tu_g0_1}:
\begin{equation}
	\begin{split}
		\mathcal{T}^u_{g_0} & = \graffe{h \in \mathcal{H}| h(\mb{x}) = h((\mb{x}-\bs{\mu}_0)^T\mb{V}_{S,0}^{-1}(\mb{x}-\bs{\mu}_0)), a.s. \; \mb{x} \in \mathcal{X}}\\
		& = \graffe{h \in \mathcal{H}| h(\mb{x}) = h(Q_{\bs{\mu}_0,\mb{V}_{S,0}}(\mb{x})), a.s. \; \mb{x} \in \mathcal{X}}\\
		&= \graffe{h \in \mathcal{H}| h\;\text{is $\sigma(\mathcal{Q})$-measurable}}.
	\end{split}
\end{equation}

\begin{IEEEproof}
	The proof follows from some small modifications of the one in \cite[Theo. 4.4]{Tsiatis}. Specifically, we need to show that:
	\begin{itemize}
		\item[\textit{i})] Any element of $\mathcal{T}^u_{i,\bs{\rho}_0}, \forall i \in \mathbb{N}$ is an element of $\mathcal{T}^u_{g_0}$, i.e. $\mathcal{T}^u_{i,\bs{\rho}_0} \subset \mathcal{T}^u_{g_0},\forall i \in \mathbb{N}$,
		\item[\textit{ii})] Any element of $\mathcal{T}^u_{g_0}$ can be expressed as an element of a given $\mathcal{T}^u_{\tilde{i},\bs{\rho}_0}$, for some $\tilde{i} \in \mathbb{N}$, or as a converging sequence of such elements.
	\end{itemize} 
	Let us start by showing \textit{i}). Each element of $\mathcal{T}^u_{i,\bs{\rho}_0}, \forall i \in \mathbb{N}$ is of the form $\mb{c}_i^T\mb{s}_{\bs{\rho}_0,i}$, where $\mb{c}_i \in \mathbb{R}^{r_i}$ and:
	\begin{equation}
		\mb{s}_{\bs{\rho}_0,i} =\nabla_{\bs{\rho}} \ln \quadre{|\mb{V}_{S,0}|^{-1/2}  \tau_{i,\bs{\rho}_0} \tonde{Q_{\bs{\mu}_0,\mb{V}_{S,0}}(\mb{x})}}.
	\end{equation}
	Consequently, $\mb{c}_i^T\mb{s}_{\bs{\rho}_0,i}$ is an element of $\mathcal{T}^u_{g_0}$, since, by inverting the order of differentiation and integration (under the regularity conditions discussed in \cite[Sects. 6.2, 6.3]{Lehmann}), it is immediate to verify that $E_0\{[\mb{s}_{\bs{\rho}_0,i}]_j\}=0$, for $j=1,\ldots,r_i$. Then $[\mb{s}_{\bs{\rho}_0,i}]_j \in \mathcal{H}$ for $j=1,\ldots,r_i$ and it is measurable w.r.t. $Q_{\bs{\mu}_0,\mb{V}_{S,0}}(\mb{x})$, then $\mb{c}_i^T\mb{s}_{\bs{\rho}_0,i} \in \mathcal{T}^u_{g_0}, \forall i \in \mathbb{N}$.
	
	Let us now move to \textit{ii}). Let start by choosing $r_{\tilde{i}}$ elements $\{\tilde{h}_j(Q_{\bs{\mu}_0,\mb{V}_{S,0}}(\mb{x}))\}_{j=1}^{r_{\tilde{i}}} \in \mathcal{T}^u_{g_0}$ such that $\tilde{h}_j$ are \textit{bounded} functions. Then, as parametric sub-model of the form in \eqref{sub_par}, we may choose the following one:
	\begin{equation}
		\label{sub_par_SS_p}
		\mathcal{P}^u_{\tau_{\bs{\rho},\tilde{i}}}  = \graffe{ p_X(\mb{x}|\bs{\rho}) = p_0(\mb{x})\quadre{1+\sum\nolimits_{j=1}^{r_{\tilde{i}}}\rho_j \tilde{h}_j(Q_{\bs{\mu}_0,\mb{V}_{S,0}}(\mb{x}))}},
	\end{equation}
where $\bs{\rho} \in \Upsilon_{\tilde{i}}$ is sufficiently small to guarantee that 
	\begin{equation}
		1+\sum\nolimits_{j=1}^{r_{\tilde{i}}}\rho_j \tilde{h}_j(Q_{\bs{\mu}_0,\mb{V}_{S,0}}(\mb{x}))\geq 0, \forall \mb{x} \in \mathcal{X},
	\end{equation}
	such then $p_X(\mb{x}|\bs{\rho}) \geq 0$. Note that, such \virg{small $\bs{\rho}$} exists since we are working with bounded functions $\tilde{h}_j$. Moreover, for each $p_X(\mb{x}|\bs{\rho}) \in \mathcal{P}_{\tau_{\bs{\rho},\tilde{i}}}$, we have that:
	\begin{equation}
	\label{rel}
		\begin{split}
			\int p_X(\mb{x}|\bs{\rho}) d\mb{x} & = \int p_0(\mb{x})d\mb{x} + \sum\nolimits_{j=1}^{r_{\tilde{i}}}\rho_j \int \tilde{h}_j(Q_{\bs{\mu}_0,\mb{V}_{S,0}}(\mb{x}))p_0(\mb{x})d\mb{x}\\
			&=1 + \sum\nolimits_{j=1}^{r_{\tilde{i}}}\rho_j \ev{\tilde{h}_j(Q_{\bs{\mu}_0,\mb{V}_{S,0}}(\mb{x}))} = 1 + 0 =1,
		\end{split}
	\end{equation}
	since $\tilde{h}_j \in \mathcal{H}$, then $p_X(\mb{x}|\bs{\rho})$ is a proper pdf and consequently $\mathcal{P}_{\tau_{\bs{\rho},\tilde{i}}}$ is a proper parametric sub-model that satisfies the conditions C1 and C2. Specifically, C2 is verified for $\bs{\rho}_0 = \mb{0}$. Now, a score vector of this specific parametric sub-model is of the form:
	\begin{equation}
		\label{pro}
		\begin{split}
			\mb{s}_{\bs{\rho}_0,\tilde{i}} &=\nabla_{\bs{\rho}} \ln p_X(\mb{x}|\bs{\gamma}_0,\tau_{i,\bs{\rho}_0}) =\\
			&= \left. \nabla_{\bs{\rho}} \ln\quadre{p_0(\mb{x})\tonde{1+\sum\nolimits_{j=1}^{r_{\tilde{i}}}\rho_j \tilde{h}_j(Q_{\bs{\mu}_0,\mb{V}_{S,0}}(\mb{x})}}\right|_{\bs{\rho}=\mb{0}} \\
			& = \left. \nabla_{\bs{\rho}}\ln\left[1+\sum\nolimits_{j=1}^{r_{\tilde{i}}}\rho_j \tilde{h}_j(Q_{\bs{\mu}_0,\mb{V}_{S,0}}(\mb{x})\right]\right|_{\bs{\rho}=\mb{0}}
			=(\tilde{h}_1,\ldots, \tilde{h}_{r_{\tilde{i}}})^T.
		\end{split}
	\end{equation}
Consequently, since, as said before, any element of $\mathcal{T}^u_{\tilde{i},\bs{\rho}_0}$ is of the form $\mb{c}^T\mb{s}_{\bs{\rho}_0,\tilde{i}}$ for some $\mb{c} \in \mathbb{R}^{r_{\tilde{i}}}$, we just need to choose $\mb{c} = \mb{e}_{j,r_{\tilde{i}}}$ to prove that $\mathcal{T}^u_{g_0} \ni \tilde{h}_j(Q_{\bs{\mu}_0,\mb{V}_{S,0}}(\mb{x})) \in \mathcal{T}^u_{\tilde{i},\bs{\rho}_0}$, where $\tilde{h}_j$ is a bounded function. The proof is completed by noticing that the set of bounded functions is dense in $\mathcal{H}$ and consequently, any element $h \in \mathcal{T}^u_{g_0}$ can be obtained as a converging sequence of bounded functions. This concludes the proof of the result in \eqref{Tu_g0_1}. We conclude this part by noticing that, $\Pi(h|\mathcal{T}^u_{g_0}) = E\{h|\mathcal{Q}\}, \; \forall h \in \mathcal{H}$ \cite[Ch. 23, Def. 4]{Jacod}.\\

Let us now move to the \textit{constrained} non-parametric sub-model $\mathcal{P}^c_{\bar{g}}$ in \eqref{con_non_par}. To derive the relevant tangent space $\mathcal{T}^c_{\bar{g}_0}$, we can reuse the previous proof except for the fact that we need to impose such constraint to the class of parametric sub-model given in \eqref{sub_par_SS_p}. Specifically, since form \eqref{Q_RES} $\mathcal{Q} = Q_{\bs{\mu}_0,\mb{V}_{S,0}}(\mb{x})$, we want that:
\begin{equation}
	\label{constraint}
	\begin{split}
		m = &\e{Q_{\bs{\mu}_0,\mb{V}_{S,0}}(\mb{x})} = \int Q_{\bs{\mu}_0,\mb{V}_{S,0}}(\mb{x}) p_0(\mb{x})\quadre{1+\sum\nolimits_{j=1}^{r_{\tilde{i}}}\rho_j \tilde{h}_j(Q_{\bs{\mu}_0,\mb{V}_{S,0}}(\mb{x}))}d\mb{x}\\
		& = \e{\mathcal{Q}} +\sum\nolimits_{j=1}^{r_{\tilde{i}}}\rho_j \e{\mathcal{Q}\tilde{h}_j(\mathcal{Q})},
	\end{split}
\end{equation} 
for all the possible $\{\rho_j\}_{j=1}^{r_{\tilde{i}}}$ and then we should impose that $\e{\mathcal{Q}\tilde{h}_j(\mathcal{Q})}=\e{(\mathcal{Q}-m)\tilde{h}_j(\mathcal{Q})}=0$ for all $j=1,\ldots,r_{\tilde{i}}$. Consequently, the tangent space $\mathcal{T}^c_{\bar{g}_0}$ of $\mathcal{P}^c_{\bar{g}}$ in \eqref{con_non_par} is given by:
\begin{equation}\label{D_c}
	\mathcal{T}^c_{\bar{g}_0} = \graffe{l \in \mathcal{T}^u_{g_0}| \e{(\mathcal{Q}-m)l(\mathcal{Q})} = 0} = \graffe{l \in \mathcal{T}^u_{g_0}| \e{\mathcal{Q}l(\mathcal{Q})} = 0} \subseteq \mathcal{H}.
\end{equation}

We now need to derive the orthogonal projection operator onto $\mathcal{T}^c_{\bar{g}_0}$, i.e. $\Pi(\cdot|\mathcal{T}^c_{\bar{g}_0})$. This can be readily done by noticing that $\mathcal{T}^u_{g_0}$ can be decomposed as the following orthogonal direct sum $\mathcal{T}^u_{g_0} = \mathcal{T}^c_{\bar{g}_0} \oplus \mathrm{Span}\{\mathcal{Q}-m\}$. This follows by the fact that: $\e{(\mathcal{Q}-m)l(\mathcal{Q})} = 0 \Leftrightarrow (\mathcal{Q}-m) \perp l(\mathcal{Q})$ where we have used the canonical inner product of $\mathcal{H}$ defined in \eqref{H_inner_prod} as $\innerprod{h_1}{h_2}_{\mathcal{H}}\triangleq \ev{h_1h_2},\; \forall h_1,h_2 \in \mathcal{H}$.\footnote{Note that we need to assume here the additional assumption $E\{\mathcal{Q}^2\} <\infty$. We note that this condition implies the assumption required in footnote 3. It imposes for example for the $t$-distribution that the degree of freedom satisfies $\nu>4$.}  Moreover, the previous decomposition of $\mathcal{T}^u_{g_0}$ implies that
\begin{equation}
	\Pi(h|\mathcal{T}^u_{g_0})=\Pi(h|\mathcal{T}^c_{g_0})+\Pi(h|\mathrm{Span}\{(\mathcal{Q}-m)\}),  \quad\forall h\in \mathcal{H},
\end{equation}
where 
\begin{equation}
	\Pi(h|\mathrm{Span}\{(\mathcal{Q}-m)\})= \frac{\innerprod{h}{(\mathcal{Q}-m)}_\mathcal{H}}{\norm{\mathcal{Q}-m}^2_\mathcal{H}} (\mathcal{Q}-m).
\end{equation}
Finally, since $\Pi(h|\mathcal{T}^u_{g_0}) = E\{h|\mathcal{Q}\}$, we have:
\begin{equation}
	\Pi(h|\mathcal{T}^c_{\bar{g}_0}) = E\{h|\mathcal{Q}\} - \frac{\eg{\mathcal{Q}h}}{\eg{(\mathcal{Q}-m)^2}}(\mathcal{Q}-m), \quad \forall h\in \mathcal{H},
\end{equation}
where $\e{\mathcal{Q}h} = \e{\mathcal{Q}\e{h|\mathcal{Q}}} = \innerprod{\Pi(h|\mathcal{T}^u_{g_0})}{(\mathcal{Q}-m)}_\mathcal{H}$. To conclude, let us have a look at the nuisance tangent space of the parametric sub-model $\mathcal{P}^c_{\bs{\gamma},s,\bs{\zeta}} \subset \mathcal{P}_{\bs{\eta},\bar{g}}$ in \eqref{P_124}. It is immediate to verify that its parametric sub-model $\mathcal{P}^c_{\bs{\gamma}_0,s_0,\bs{\zeta}}$ is exactly of the form given in \eqref{sub_par}. Consequently, its (finite-dimensional) tangent space $\mathcal{T}^c_{\bs{\zeta}_0}$ is included in $\mathcal{T}^c_{\bar{g}_0}$, i.e. $\mathcal{T}^c_{\bs{\zeta}_0} \subset \mathcal{T}^c_{\bar{g}_0}$.
\end{IEEEproof}
\label{proof_tan_spa}
\vspace{0.5cm}

\item {\bf A property of the scale function $S(\bs{\Sigma}_0) = |\bs{\Sigma}_0|^{1/m}$}. The aim of this appendix is to recall the proof, already derived in \cite{PAINDAVEINE}, of the determinant-based scale. Specifically, as corollary of the Proposition \ref{Prop_sg}, in \cite{PAINDAVEINE}, it has been show that, if the scale function $S_d(\bs{\Sigma}_0) \triangleq |\bs{\Sigma}_0|^{1/m}$ is adopted, then the parametric FIM $\mb{I}_{\bs{\gamma}_0}$ in the parametric model $\mathcal{P}_{\bs{\gamma}}$ in \eqref{P1_eta} is equal to the efficient SFIM $\bar{\mb{I}}(\bs{\gamma}_0|s_0,\bar{g}_0)$ for $\mathcal{P}_{\bs{\eta},\bar{g}}$. To prove this, we need to show that, if $S(\bs{\Sigma}_0) = |\bs{\Sigma}_0|^{1/m}$, then $\Pi(\mb{s}_{\vsVso}|\mathcal{T}_{s_0})=\mb{0}$. From \eqref{proj_Vs}, we have that:
\begin{equation}
	\Pi(\mb{s}_{\vsVso}|\mathcal{T}_{s_0})=\mb{0} \Leftrightarrow \mb{M}_S^{\mb{V}_{S,0}}\vVsoinv = \mb{0}
\end{equation}
or equivalently $\Pi(\mb{s}_{\vsVso}|\mathcal{T}_{s_0})=\mb{0}$ iff $\vVsoinv \in \text{Ker}\;\mb{M}_S^{\mb{V}_{S,0}} \cap \cvec{\mathcal{S}_m^\mathbb{R}} = \{a \cdot \cvec{\mb{D}_S^{\mb{V}_{S,0}}},a\in \mathbb{R} \}$ where the matrix $\mb{D}_S^{\mb{V}_{S,0}}$ is formally defined in \eqref{D_mat} of \ref{App6}. As already noted in \cite[Th. 3.1]{PAINDAVEINE}, this is a direct consequence of the Property P4 of the matrix $\mb{M}_S^{\mb{V}_{S,0}}$ (see appendix \ref{App6}). In fact, for $S\equiv S_d$, we have that $\mb{D}_{S_d}^{\bs{\Sigma}_0} = m^{-1}|\bs{\Sigma}_0|^{1/m}\bs{\Sigma}_0^{-1} = m^{-1}\mb{V}_{S,0}^{-1}$. Consequently, we immediately have that $\mathrm{vec}(\mb{V}_{S_d,0}^{-1}) \in \text{Ker}\;\mb{M}_{S_d}^{\mb{V}_{S,0}} \cap \cvec{\mathcal{S}_m^\mathbb{R}}$ and then $\Pi(\mb{s}_{\vsVso}|\mathcal{T}_{s_0})=\mb{0}$.
\hfill
\QED
\label{App_det_based}
\vspace{0.5cm}

\item \textbf{Explicit expressions of the CRB in the model $\mathcal{P}^c_{\bs{\gamma},s}$}.
Let us first evaluate the CRB on $\bs{\nu}_0 \triangleq (\bs{\mu}_0^T,\vecs{\bs{\Sigma}_0}^T)^T \in \Omega$. To this end, we note that the two components of the related score vector $\mb{s}_{\bs{\nu}_0} = (\mb{s}_{\bs{\mu}_0}^T, \mb{s}_{\vsSs}^T)^T$ have been already introduced in \eqref{s_mu} and \eqref{s_Sigma}. Consequently, from standard calculation and by using the independence between $\mathcal{Q}$ and $\mb{u}$ (along with the properties of $\mb{u}$), the FIM for $\bs{\nu}_0 \in \Omega$ is given by:
\begin{equation}\label{Par_model_nu_block}
	\mb{I}_{\bs{\nu}_0} \triangleq E_0\{\mb{s}_{\bs{\nu}_0}\mb{s}_{\bs{\nu}_0}^T\} = \tonde{
		\begin{array}{cc}
			\mb{I}_{\bs{\mu}_0} & \mb{I}_{\bs{\mu}_0,\vsSs} \\
			\mb{I}_{\bs{\mu}_0,\vsSs}^T	& \mb{I}_{\vsSs} \\
		\end{array}
	},
\end{equation}
where
\begin{equation}
	\mb{I}_{\bs{\mu}_0} \triangleq E_0\graffe{\mb{s}_{\bs{\mu}_0}\mb{s}_{\bs{\mu}_0}^T} = \beta(\bar{g}_0) \bs{\Sigma}_0^{-1},
	\tag{\ref{I_mu}}
\end{equation}
\begin{equation}\label{I_m_S}
	\mb{I}_{\bs{\mu}_0,\vsSs} \triangleq E_0\graffe{\mb{s}_{\bs{\mu}_0}\mb{s}_{\vsSs}^T}=\mb{0},
\end{equation}
\begin{equation}\label{I_vecs Sigma}
\mb{I}_{\vsSs}\triangleq E_0\graffe{\mb{s}_{\vsSs}\mb{s}_{\vsSs}^T}
= \mb{D}_m^T\quadre{\frac{1}{2}\alpha(\bar{g}_0)(\bs{\Sigma}_0^{-1}\otimes \bs{\Sigma}_0^{-1})
+ \frac{1}{4}({\alpha(\bar{g}_0) -1})\cvec{\bs{\Sigma}_0^{-1}}\cvec{\bs{\Sigma}_0^{-1}}^T}\mb{D}_m.
\end{equation}
Then the CRB on $\bs{\nu}_0$ can be derived from the block diagonal structured FIM in \eqref{Par_model_nu_block} as:
\begin{equation}\label{CRB nu}
	\mathrm{CRB}(\bs{\nu}_0) 
	\triangleq \mb{I}_{\bs{\nu}_0}^{-1}
=\tonde{
		\begin{array}{cc}
			\mathrm{CRB}(\bs{\mu}_0) & \mb{0} \\
			\mb{0} & \mathrm{CRB}(\vsSs) \\
		\end{array}
	},
\end{equation}
where:
\begin{equation}\label{CRB vecs Sigma}
	\mathrm{CRB}(\bs{\mu}_0) \triangleq \mb{I}_{\bs{\mu}_0}^{-1} = \beta(\bar{g}_0)^{-1} \bs{\Sigma}_0\ \
	\mbox{and}\ \
	\mathrm{CRB}(\vsSs) \triangleq \mb{I}_{\vsSs}^{-1}.
\end{equation}
To calculate $\mb{I}_{\vsSs}^{-1}$, 
let us rewrite \eqref{I_vecs Sigma} in the form $\mb{I}_{\vsSs}= \mb{D}_m^T[ \mb{A}+ \mb{a} \mb{a}^T] \mb{D}_m$. The inverse of the middle term of 
\eqref{I_vecs Sigma} can be derived from the inversion matrix lemma giving:
\begin{equation}\label{inversion matrix lemma}
	[ \mb{A}+ \mb{a} \mb{a}^T]^{-1}
	=
	2\alpha^{-1}(\bar{g}_0)(\bs{\Sigma}_0\otimes \bs{\Sigma}_0)
	-\frac{2\alpha^{-1}(\bar{g}_0)(\alpha(\bar{g}_0)-1)}{(m+2)\alpha(\bar{g}_0)-m}\cvec{\bs{\Sigma}_0}\cvec{\bs{\Sigma}_0}^T.
\end{equation}
Then using $\mb{D}_m \mb{D}_m^{\#}= \frac{1}{2}(\mb{I}_{m^2}+\mb{K}_m)$ and $\mb{K}_m (\bs{\Sigma}_0\otimes  \bs{\Sigma}_0)= (\bs{\Sigma}_0\otimes  \bs{\Sigma}_0)\mb{K}_m$ \cite[Ch. 3]{Magnus}, we straightforwardly get
\begin{equation}\label{inversion matrix lemma b}
	\mb{D}_m \mb{D}_m^{\#}[ \mb{A}+ \mb{a} \mb{a}^T]^{-1} \mb{D}_m^{\#T}\mb{D}^T_m
	=
	\alpha^{-1}(\bar{g}_0)(\mb{I}_{m^2}+\mb{K}_{m})(\bs{\Sigma}_0\otimes \bs{\Sigma}_0)
	-\frac{2\alpha^{-1}(\bar{g}_0)(\alpha(\bar{g}_0)-1)}{(m+2)\alpha(\bar{g}_0)-m}\cvec{\bs{\Sigma}_0}\cvec{\bs{\Sigma}_0}^T.
\end{equation}
Finally, using $\mb{K}_m \mb{D}_m = \mb{D}_m$, $\mb{D}_m^{\#}\mb{D}_m =\mb{I}_{m(m+1)/2}$, 
\cite[Ch. 3]{Magnus}, it is easy to check that
\begin{equation}\label{FIM times CRB=I}
	\underbrace{\left(\mb{D}_m^{T}[ \mb{A}+ \mb{a} \mb{a}^T] \mb{D}_m \right)}_{\mb{I}_{\vsSs}}
	\mb{D}_m^{\#}\underbrace{\left(\mb{D}_m \mb{D}_m^{\#}[ \mb{A}+ \mb{a} \mb{a}^T]^{-1} \mb{D}_m^{\#T}\mb{D}^T_m\right)}_{\eqref{inversion matrix lemma b}}\mb{D}_m^{\#T}
	=
	\mb{I}_{m(m+1)/2},
\end{equation}
and thus, since $\mb{D}_m^{\#}{\mb K}_m=\mb{D}_m^{\#}$ \cite[Ch. 3]{Magnus}
\begin{equation}\label{CRB vecs Sigma b}
	\mathrm{CRB}(\vsSs) =
	2\alpha^{-1}(\bar{g}_0)\mb{D}_m^{\#}\quadre{(\bs{\Sigma}_0\otimes \bs{\Sigma}_0)
		-\frac{(\alpha(\bar{g}_0)-1)}{(m+2)\alpha(\bar{g}_0)-m}\cvec{\bs{\Sigma}_0}\cvec{\bs{\Sigma}_0}^T}\mb{D}_m^{\# T}.
\end{equation}

In the second step of the proof, the CRB on 
 $\bs{\eta}_0$ is given by 
\begin{equation}\label{CRB eta}
	\mathrm{CRB}(\bs{\eta}_0)
	=
	\mb{J}[\mb{w}^{-1}](\bs{\nu}_0)
	\mathrm{CRB}(\bs{\nu}_0)
	[\mb{J}[\mb{w}^{-1}](\bs{\nu}_0)]^T,
\end{equation}
 from the inverse diffeomorphism of $\mb{w}$ \eqref{diff_par}, whose 
 Jacobian matrix $\mb{J}[\mb{w}^{-1}](\bs{\nu}_0)$ is given by \eqref{Jac_e}.

Let us explicitly evaluate the term $\nabla_{\cvec{\bs{\Sigma}}}^T S(\bs{\Sigma}_0)$ in \eqref{Jac_e}. From \eqref{shape_m}, we have
\begin{equation}
	\nabla_{\cvec{\bs{\Sigma}}}^T\cvec{\mb{V}_{S}}
	=
	S^{-1}(\bs{\Sigma})\left[\mb{I}_{m^2}-\cvec{\mb{V}_{S}}\nabla_{\cvec{\bs{\Sigma}}}^T S(\bs{\Sigma})\right].
\end{equation}
Then from $\ovecs{\mb{V}_S}=\underline{\mb{I}}_m\vecs{\mb{V}_S}$ and $\cvec{\bs{\Sigma}}=\mb{D}_m\vecs{\bs{\Sigma}}$, we get:
\begin{equation}\label{Jacobien proof}
	\nabla_{\vecs{\bs{\Sigma}}}^T\ovecs{\mb{V}_{S}}
	=
	S^{-1}(\bs{\Sigma})\underline{\mb{I}}_m\mb{D}_m^{\#} \left[\mb{I}_{m^2}-\cvec{\mb{V}_{S}}\nabla_{\cvec{\bs{\Sigma}}}^T S(\bs{\Sigma})\right]\mb{D}_m
\end{equation}
Consequently, the CRB on $\bs{\eta}_0$ is given by 
\begin{equation}\label{CRB eta b}
	\mathrm{CRB}(\bs{\eta}_0)
	=
	\tonde{
		\begin{array}{ccc}
			\mathrm{CRB}(\bs{\mu}_0) & \mb{0}  & \mb{0} \\
			\mb{0} & \mathrm{CRB}(\ovecs{\mb{V}_{S,0}}|s_0) & \bs{\Psi} \\
			\mb{0} & \bs{\Psi}^T &  \mathrm{CRB}(s_0|\ovecs{\mb{V}_{S,0}}) \\
		\end{array}
	},
\end{equation}
where $\mathrm{CRB}(\bs{\mu}_0)$ is given by \eqref{CRB mu} and:
\begin{eqnarray}
\label{CRB V proof}
\mathrm{CRB}(\ovecs{\mb{V}_{S,0}}|s_0)
&=&
 \nabla_{\vecs{\bs{\Sigma}}}^T\ovecs{\mb{V}_{S,0}}
\mathrm{CRB}(\vsSs)
[\nabla_{\vecs{\bs{\Sigma}}}\ovecs{\mb{V}_{S,0}}],
\\
\label{CRB s proof}
\mathrm{CRB}(s_0|\ovecs{\mb{V}_{S,0}})
&=&
 \nabla_{\vecs{\bs{\Sigma}}}^TS(\bs{\Sigma}_0)
\mathrm{CRB}(\vsSs)
 \nabla_{\vecs{\bs{\Sigma}}}S(\bs{\Sigma}_0),
\\
\label{CRB V,s proof}	
 \bs{\Psi}
&=&
 \nabla_{\vecs{\bs{\Sigma}}}^T\ovecs{\mb{V}_{S,0}}
\mathrm{CRB}(\vsSs)
 \nabla_{\vecs{\bs{\Sigma}}}S(\bs{\Sigma}_0).
\end{eqnarray}
Then, applying Euler's homogeneous function theorem to the score function $S(\bs{\Sigma})$ which is homogeneous of order one:
\begin{equation}\label{identity}
 \nabla_{\cvec{\bs{\Sigma}}}^T S(\bs{\Sigma})\cvec{\bs{\Sigma}}
 =S(\bs{\Sigma}),
\end{equation}
the relation  $\mb{D}_m \mb{D}_m^{\#}= \frac{1}{2}(\mb{I}_{m^2}+\mb{K}_m)$ and the fact that $\mb{K}_m \nabla_{{\rm vec}({\bs{\Sigma})}} S(\bs{\Sigma}_0)= \nabla_{{\rm vec}({\bs{\Sigma})}} S(\bs{\Sigma}_0)$ the relations \eqref{CRB V proof}, \eqref{CRB s proof} and \eqref{CRB V,s proof} can be explicitly expressed as:
\begin{equation}
 	\mathrm{CRB}(\ovecs{\mb{V}_{S,0}}|s_0) 
 	=
 	\alpha^{-1}(\bar{g}_0) 
 	\underline{\mb{I}}_m
 	\mb{D}_m^{\#}\mb{P}_S(\mb{V}_{S,0})
 	(\mb{I}_{m^2}+\mb{K}_{m})(\mb{V}_{S,0}\otimes \mb{V}_{S,0})
 	\mb{P}_S^T(\mb{V}_{S,0})\mb{D}_m^{\#T}
 	\underline{\mb{I}}_m^T,
 \end{equation}
 that correspond to the expression reported in \eqref{CRB V} and:
 \begin{equation}\label{CRB s}
 	\mathrm{CRB}(s_0|\ovecs{\mb{V}_{S,0}})
 	=
 	\frac{2s_0^2}{\alpha(\bar{g}_0)} \quadre{\nabla_{{\rm vec}({\bs{\Sigma})}}^T S(\bs{\Sigma}_0)(\mb{V}_{S,0}\otimes \mb{V}_{S,0})
 		\nabla_{{\rm vec}({\bs{\Sigma})}} S(\bs{\Sigma}_0)-\tonde{\frac{\alpha(\bar{g}_0)-1}{(m+2)\alpha(\bar{g}_0)-m} }},
 \end{equation}
 \begin{equation}\label{CRB Psi}
 	\bs{\Psi}
 	=
 	2\alpha^{-1}(\bar{g}_0)s_0
 	\underline{\mb{I}}_m
 	\mb{D}_m^{\#}\mb{P}_S(\mb{V}_{S,0})
 	(\mb{V}_{S,0}\otimes \mb{V}_{S,0})\nabla_{{\rm vec}({\bs{\Sigma})}}S(\bs{\Sigma}_0),
 \end{equation}
 
Remarkably, it can be shown that, for the scale functional $S_d(\bs{\Sigma}) \triangleq |\bs{\Sigma}|^{1/m}$ the term $\bs{\Psi}$ cancels out, in accordance with Proposition \ref{Prop_sg}. In fact, for $S_d(\bs{\Sigma})$, we get
$\nabla_{{\rm vec}({\bs{\Sigma})}} S_d(\bs{\Sigma}_0)=\frac{s_0}{m}{\rm vec}(\bs{\Sigma}_0^{-1})$ from \cite[Ch. 8, Th. 1]{Magnus}
and 
$\mb{P}_{S_d}(\mb{V}_{S,0}) =\mb{I}_{m^2}-\frac{1}{m}{\rm vec}(\bs{\Sigma}_0){\rm vec}(\bs{\Sigma}_0^{-1})^T$. 
Consequently $(\bs{\Sigma}_0\otimes \bs{\Sigma}_0)\nabla_{{\rm vec}({\bs{\Sigma})}} S_d(\bs{\Sigma}_0)=\frac{s_0}{m}{\rm vec}(\bs{\Sigma}_0)$ and thus

$\mb{P}_{S_d}(\mb{V}_{S,0})(\bs{\Sigma}_0\otimes \bs{\Sigma}_0)\nabla_{{\rm vec}({\bs{\Sigma})}} S_d(\bs{\Sigma}_0)
=\frac{s_0}{m}{\rm vec}(\bs{\Sigma}_0)-\frac{s_0}{m^2} \trace{{\mb I}_m}{\rm vec}(\bs{\Sigma}_0)=\mb{0}$ and therefore $\bs{\Psi}=\mb{0}$.

Moreover, $\mb{P}_{S_d}(\mb{V}_{S,0})
(\mb{I}_{m^2}+\mb{K}_{m})(\mb{V}_{S,0}\otimes \mb{V}_{S,0})
\mb{P}_{S_d}^T(\mb{V}_{S,0})=(\mb{I}_{m^2}+\mb{K}_{m})(\mb{V}_{S,0}\otimes \mb{V}_{S,0})
-\frac{2}{m}{\rm vec}(\mb{V}_{S,0}){\rm vec}(\mb{V}_{S,0})^T$ using in particular
${\rm vec}(\mb{V}_{S,0}^{-1})^T(\mb{V}_{S,0}\otimes \mb{V}_{S,0}){\rm vec}(\mb{V}_{S,0}^{-1})
=\trace{\mb{V}_{S,0}^{-1} \mb{V}_{S,0} \mb{V}_{S,0}^{-1} \mb{V}_{S,0}}=m$, which proves  \eqref{CRB V d}.

Finally, $\nabla_{{\rm vec}({\bs{\Sigma})}}^T {S_d}(\bs{\Sigma}_0)(\mb{V}_{S,0}\otimes \mb{V}_{S,0})
\nabla_{{\rm vec}({\bs{\Sigma})}} S_{d}(\bs{\Sigma}_0)
=\frac{s_0^2}{m^2}{\rm vec}(\mb{V}_{S,0}^{-1})^T(\mb{V}_{S,0}\otimes \mb{V}_{S,0}){\rm vec}(\mb{V}_{S,0}^{-1})=\frac{s_0^2}{m}$, which proves that \eqref{CRB s} reduces to:
\begin{equation}
	\label{CRB s d}
	\mathrm{CRB}(s_{d,0}|\ovecs{\mb{V}_{S_d,0}})
	=
	\frac{4|\bs{\Sigma}|^{2/m}}{m[m(\alpha(\bar{g}_0)-1)+2\alpha(\bar{g}_0)]}.
\end{equation}
\hfill
\QED
\label{App_crb2}

\vspace{0.5cm}
\item \textbf{Proof of Proposition \ref{Prop_low-rank model}}. The different steps of the proof are based on some reasoning and notations of \cite[Sec. VI of supplement material]{ABEIDA2019} for C-CES distributions and the general "low-rank" model \eqref{I_gamma_low-rank c},
that itself takes up the steps of the proof presented in \cite{Stoica_CRB}.

We deduce from the FIM in $\mb{I}_{\bs{\theta}_0}$ \eqref{I_theta_par}, the following expressions of the sub-blocks:
\begin{eqnarray}
\label{I xi}
\mb{I}_{\bs{\xi}_0} 
&=& 
\mb{J}_{\bs{\xi}}[\cvec{\bs{\Sigma}_0}]^T\quadre{(\bs{\Sigma}_0^{-1/2} \otimes \bs{\Sigma}_0^{-1/2}){\bf T}^{1/2}}
\quadre{{\bf T}^{1/2}(\bs{\Sigma}_0^{-1/2} \otimes \bs{\Sigma}_0^{-1/2})}
\mb{J}_{\bs{\xi}}[\cvec{\bs{\Sigma}_0}]
\\
\label{I gamma xi}
\mb{I}_{\bs{\gamma}_0 \bs{\xi}_0} 
&=& 
\mb{J}_{\bs{\gamma}}[\cvec{\bs{\Sigma}_0}]^T\quadre{(\bs{\Sigma}_0^{-1/2} \otimes \bs{\Sigma}_0^{-1/2}){\bf T}^{1/2}}
\quadre{{\bf T}^{1/2}(\bs{\Sigma}_0^{-1/2} \otimes \bs{\Sigma}_0^{-1/2})}
\mb{J}_{\bs{\xi}}[\cvec{\bs{\Sigma}_0}]
\end{eqnarray}
with
\begin{equation}\label{def T}
{\bf T}\triangleq 
2^{-1}\alpha(\bar{g}_0){\bf I}_{m^2}+4^{-1}(\alpha(\bar{g}_0)-1) \mathrm{vec}({\bf I}_{m})\mathrm{vec}({\bf I}_{m})^T,
\end{equation}
which can be written, by means of the notation used in \cite[Sec; V of supplement material]{ABEIDA2019}, in the following compact form
\begin{equation}\label{I gamma xi xi compact}
\mb{I}_{\bs{\xi}_0} 
={\boldsymbol \Delta}^T{\boldsymbol \Delta}
\ \ 
\mbox{and}
\ \ 
\mb{I}_{\bs{\gamma}_0 \bs{\xi}_0} 
={\bf G}^T{\boldsymbol \Delta},
\end{equation} 
with
\begin{equation}\label{def Delta G}
{\boldsymbol \Delta}\triangleq
\quadre{{\bf T}^{1/2}(\bs{\Sigma}_0^{-1/2} \otimes \bs{\Sigma}_0^{-1/2})}
\mb{J}_{\bs{\xi}}[\cvec{\bs{\Sigma}_0}]
\ \ \mbox{and}\ \ 
{\bf G}\triangleq \quadre{{\bf T}^{1/2}(\bs{\Sigma}_0^{-1/2} \otimes \bs{\Sigma}_0^{-1/2})}
\mb{J}_{\bs{\gamma}}[\cvec{\bs{\Sigma}_0}].
\end{equation} 
Consequently the left hand side of \eqref{cond_corr_2_RES} can be expressed as:
\begin{equation}\label{LHS}
{\bf G}^T\left[
{\bf I}_{m^2}-{\boldsymbol \Delta}
({\boldsymbol \Delta}^T{\boldsymbol \Delta})^{-1}{\boldsymbol \Delta}^T\right]
{\bf T}^{-1/2}(\bs{\Sigma}_0^{1/2} \otimes \bs{\Sigma}_0^{1/2})\cvec{\bs{\Sigma}_0^{-1}}
=
{\bf G}^T{\boldsymbol \Pi}^{\bot}_{\boldsymbol \Delta}
{\bf T}^{-1/2}\cvec{{\bf I}_m},
\end{equation} 
with
${\boldsymbol \Pi}^{\bot}_{\boldsymbol \Delta} \triangleq
{\bf I}_{m^2}-{\boldsymbol \Delta}
({\boldsymbol \Delta}^T{\boldsymbol \Delta})^{-1}{\boldsymbol \Delta}^T$.
This implies that condition \eqref{cond_corr_2_RES} of Proposition \ref{prop_corr_1_RES} is satisfied iif 
\begin{equation}\label{condition VI.1 a}
\cvec{{\bf I}_m}^T{\bf T}^{-1/2}{\boldsymbol \Pi}^{\bot}_{\boldsymbol \Delta}
{\bf g}_k=0,\ \ k=1,..,q,
\end{equation}
where ${\bf g}_k, k=1,..,q$ denotes the $k$th column of ${\bf G}$.
Let's further partition the matrix ${\boldsymbol \Delta}$ as 
\begin{equation}\label{partition Delta}
{\boldsymbol \Delta}={\bf T}^{1/2}(\bs{\Sigma}_0^{-1/2} \otimes \bs{\Sigma}_0^{-1/2})
\left[
\mb{J}_{\vecs{\bs{\Xi}}}[\cvec{\bs{\Sigma}_0}]\ |\ \mb{J}_{\lambda}[\cvec{\bs{\Sigma}_0}]
\right]
\triangleq [{\bf V}\ |\ {\bf u}_n],
\end{equation}
with ${\bf u}_n={\bf T}^{1/2}(\bs{\Sigma}_0^{-1/2} \otimes \bs{\Sigma}_0^{-1/2})\cvec{{\bf I}_m}={\bf T}^{1/2}\cvec{\bs{\Sigma}_0^{-1}}$.
It follows from \cite[rel. (14)]{Stoica_CRB} that 
\begin{equation}\label{projector Delta}
{\boldsymbol \Pi}^{\bot}_{\boldsymbol \Delta}
=
{\boldsymbol \Pi}^{\bot}_{\bf V}
-\frac{{\boldsymbol \Pi}^{\bot}_{\bf V}{\bf u}_n{\bf u}_n^T{\boldsymbol \Pi}^{\bot}_{\bf V}}{{\bf u}_n^T{\boldsymbol \Pi}^{\bot}_{\bf V}{\bf u}_n}.
\end{equation}
Reporting expression \eqref{projector Delta} of ${\boldsymbol \Pi}^{\bot}_{\boldsymbol \Delta}$ in \eqref{condition VI.1 a}, condition \eqref{cond_corr_2_RES} of Proposition \ref{prop_corr_1_RES} is satisfied iif 
\begin{equation}\label{condition VI.1 b}
\cvec{{\bf I}_m}^T{\bf T}^{-1/2}{\boldsymbol \Pi}^{\bot}_{\bf V}{\bf g}_k
-\frac{(\cvec{{\bf I}_m}^T{\bf T}^{-1/2}{\boldsymbol \Pi}^{\bot}_{\bf V}{\bf u}_n)({\bf u}_n^T{\boldsymbol \Pi}^{\bot}_{\bf V}{\bf g}_k)}
{{\bf u}_n^T{\boldsymbol \Pi}^{\bot}_{\bf V}{\bf u}_n}
=0,\ \ k=1,..,q.
\end{equation}
Consequently to conclude the proof, it is sufficient to prove the two equalities:
\begin{equation}\label{condition VI.1 c}
\cvec{{\bf I}_m}^T{\bf T}^{-1/2}{\boldsymbol \Pi}^{\bot}_{\bf V}{\bf g}_k=0
\ \ \mbox{and}\ \ 
{\bf u}_n^H{\boldsymbol \Pi}^{\bot}_{\bf V}{\bf g}_k=0,\ \ 
k=1,..,q.
\end{equation}
From the definition of ${\bf G}$ \eqref{def Delta G} and the derivative
$ \boldsymbol{\Sigma}_k^0={\bf A}_k^0 \boldsymbol{\Xi}_0{\bf A}_0^T
+{\bf A}_0 \boldsymbol{\Xi}_0{({\bf A}_k^0})^T$, we straightforwardly deduce that
\begin{equation} \label{expression gk}
{\bf g}_k
=
{\bf T}^{1/2}{\rm vec}({\bf Z}_k+{\bf Z}_k^T)
\end{equation}
with
\begin{equation} \label{def Zk}
{\bf Z}_k
\triangleq
\boldsymbol{\Sigma}_0^{-1/2}{\bf A}_0\boldsymbol{\Xi}_0{({\bf A}_k^0})^T\boldsymbol{\Sigma}_0^{-1/2}.
\end{equation}
Likewise from the definition of ${\bf V}$ \eqref{partition Delta} and 
${\bf J}_{{\rm vecs}( \boldsymbol{\Xi})}[{\rm vec}({\bf \Sigma}_0)]=({\bf A}_0\otimes{\bf A}_0){\bf D}_p$, we straightforwardly get:
\begin{equation}
\label{expression V}
{\bf V}
=
{\bf T}^{1/2}{\bf W}_0{\bf D}_p
\end{equation}
with
\begin{equation} \label{definition W}
{\bf W}_0
\triangleq
{\bf \Sigma}_0^{-1/2}{\bf A}_0 \otimes {\bf \Sigma}_0^{-1/2}{\bf A}_0.
\end{equation}
Consequently, ${\boldsymbol \Pi}^{\bot}_{\bf V}$ is written in the form:
\begin{equation}\label{Pi V orthogonal}
{\boldsymbol \Pi}^{\bot}_{\bf V}
={\bf I}_{m^2}-{\bf T}^{1/2}{\bf W}_0{\bf D}_p[{\bf D}_p^T({\bf W}_0^T{\bf T}{\bf W}_0){\bf D}_p]^{-1}{\bf D}_p^T{\bf W}_0^T{\bf T}^{1/2},
\end{equation}
where ${\bf W}_0^T{\bf T}{\bf W}_0
=\frac{\alpha(\bar{g}_0)}{2}({\bf A}_0^T{\bf \Sigma}_0^{-1}{\bf A}_0 \otimes {\bf A}_0^T{\bf \Sigma}_0^{-1}{\bf A}_0)
+\frac{\alpha(\bar{g}_0)-1}{4}{\rm vec}({\bf A}_0^T{\bf \Sigma}_0^{-1}{\bf A}_0){\rm vec}({\bf A}_0^T{\bf \Sigma}_0^{-1}{\bf A}_0)^T$.
By noticing that ${\bf W}_0^T{\bf T}{\bf W}_0$ is structured in the form  ${\bf A}+ {\bf a} {\bf a}^T$  and applying the inversion matrix lemma with the trick in \eqref{FIM times CRB=I}, we get:
\begin{equation}\label{Pi V orthogonal b}
{\boldsymbol \Pi}^{\bot}_{\bf V}
={\bf I}_{m^2}-{\bf T}^{1/2}{\bf B}_0{\bf T}^{1/2},
\end{equation}
with 
\begin{equation}\label{definition B}
{\bf B}_0
\triangleq
\left(
{\bf I}_{m^2}+{\bf K}_m
\right)
\left[
\frac{1}{\alpha(\bar{g}_0)}({\bf H}_{1,0}\otimes {\bf H}_{1,0})
-\frac{\alpha(\bar{g}_0)-1}{\alpha(\bar{g}_0)(2\alpha(\bar{g}_0)+(\alpha(\bar{g}_0)-1)p}{\rm vec}({\bf H}_{1,0}){\rm vec}({\bf H}_{1,0})^T
\right]
\end{equation}
with ${\bf H}_{1,0} \triangleq {\bf \Sigma}_0^{-1/2}{\bf A}_0({\bf A}_0^T{\bf \Sigma}_0^{-1}{\bf A}_0)^{-1}{\bf A}_0^T{\bf \Sigma}_0^{-1/2}$.
Finally, both left hand sides of expressions \eqref{condition VI.1 c} follow
\begin{eqnarray}
\label{condition VI.1 d}
\cvec{{\bf I}_m}^T{\bf T}^{-1/2}{\boldsymbol \Pi}^{\bot}_{\bf V}{\bf g}_k
&=&
\cvec{{\bf I}_m}^T\cvec{{\bf Z}_k+{\bf Z}_k^T}
-\cvec{{\bf I}_m}^T{\bf B}_0{\bf T}\cvec{{\bf Z}_k+{\bf Z}_k^T}
\\
{\bf u}_n^H{\boldsymbol \Pi}^{\bot}_{\bf V}{\bf g}_k
\label{condition VI.1 e}
&=&
\cvec{{\boldsymbol \Sigma}_0^{-1}}^T{\bf T}\cvec{{\bf Z}_k+{\bf Z}_k^T}
-\cvec{{\boldsymbol \Sigma}_0^{-1}}^T{\bf T}{\bf B}_0{\bf T}\cvec{{\bf Z}_k+{\bf Z}_k^T}.
\end{eqnarray}
Reporting the expressions of ${\bf T}$ \eqref{def T}, ${\bf Z}_k$ \eqref{def Zk} and ${\bf B}_0$ \eqref{definition B} in the right hand side of \eqref{condition VI.1 d} and \eqref{condition VI.1 e}, and using ${\bf H}_{1,0}^2={\bf H}_{1,0}$ and ${\rm tr}({\bf H}_{1,0})=p$,
terms \eqref{condition VI.1 d} and \eqref{condition VI.1 e} are proven after cumbersome calculations  to be equal to zero.
\hfill
\QED
\label{proof_Prop_low-rank model}

\end{enumerate}

\bibliographystyle{IEEEtran}
\bibliography{ref_semipar_eff_estim}

\end{document}